\documentclass{article}

\usepackage[main, final]{neurips_2026}

\usepackage[utf8]{inputenc} 
\usepackage[T1]{fontenc}    
\usepackage{hyperref}       
\usepackage{url}            
\usepackage{booktabs}       
\usepackage{amsfonts}       
\usepackage{nicefrac}       
\usepackage{microtype}      
\usepackage{xcolor}         

\usepackage{amsmath}
\usepackage{amssymb}
\usepackage{mathtools}
\usepackage{amsthm}
\usepackage{microtype}
\usepackage{graphicx}
\usepackage{subcaption}
\usepackage{booktabs} 
\usepackage{algorithm}
\usepackage{algorithmic}
\usepackage{enumitem}

\usepackage{amsmath,amsfonts,bm}
\usepackage{xparse}

\DeclareDocumentCommand\W{ g g }{%
        \IfNoValueTF {#1} {\mathbf{W}} {
            \IfNoValueTF {#2} {\mathbf{W}^{(#1)}}{\mathbf{W}^{(#1)}_{#2}}
        }
}

\DeclareDocumentCommand\bias{ g g }{%
        \IfNoValueTF {#1} {\mathbf{b}} {
            \IfNoValueTF {#2} {\mathbf{b}^{(#1)}}{\mathbf{b}^{(#1)}_{#2}}
        }
}

\DeclareDocumentCommand\betavar{ g g }{%
        \IfNoValueTF {#1} {\bm{\beta}} {
            \IfNoValueTF {#2} {{\bm{\beta}^{(#1)}}{}}{\bm{\beta}^{(#1)}_{#2}}
        }
}

\DeclareDocumentCommand\xivar{ g g }{%
        \IfNoValueTF {#1} {\bm{\xi}} {
            \IfNoValueTF {#2} {{\bm{\xi}^{(#1)}}{}}{\bm{\xi}^{(#1)}_{#2}}
        }
}

\DeclareDocumentCommand\xivarn{ g g }{%
        \IfNoValueTF {#1} {\bm{\xi^-}} {
            \IfNoValueTF {#2} {\bm{\xi^-}^{+(#1)}}{\bm{\xi^-}^{+(#1)}_{#2}}
        }
}

\DeclareDocumentCommand\xivarp{ g g }{%
        \IfNoValueTF {#1} {\bm{\xi^+}} {
            \IfNoValueTF {#2} {\bm{\xi^+}^{+(#1)}}{\bm{\xi^+}^{+(#1)}_{#2}}
        }
}

\DeclareDocumentCommand\nuvar{ g g }{%
        \IfNoValueTF {#1} {{\bm{\nu}}} {
            \IfNoValueTF {#2} {{\bm{\nu}^{(#1)}}{}}{\nu^{(#1)}_{#2}{}}
        }
}

\DeclareDocumentCommand\hnuvar{ g g }{%
        \IfNoValueTF {#1} {\bm{\hat{\nu}}} {
            \IfNoValueTF {#2} {{\bm{\hat{\nu}}^{(#1)}}{}}{\hat{\nu}^{(#1)}_{#2}{}}
        }
}

\DeclareDocumentCommand\muvar{ g g }{%
        \IfNoValueTF {#1} {\bm{\mu}} {
            \IfNoValueTF {#2} {{\bm{\mu}^{(#1)}}{}}{\mu^{(#1)}_{#2}}
        }
}

\DeclareDocumentCommand\tauvar{ g g }{%
        \IfNoValueTF {#1} {\bm{\tau}} {
            \IfNoValueTF {#2} {{\bm{\tau}^{(#1)}}{}}{\tau^{(#1)}_{#2}}
        }
}

\DeclareDocumentCommand\pivar{ g g }{%
        \IfNoValueTF {#1} {\bm{\pi}} {
            \IfNoValueTF {#2} {{\bm{\pi}^{(#1)}}{}}{\pi^{(#1)}_{#2}}
        }
}

\DeclareDocumentCommand\gammavar{ g g }{%
        \IfNoValueTF {#1} {\bm{\gamma}} {
            \IfNoValueTF {#2} {{\bm{\gamma}^{(#1)}}{}}{\gamma^{(#1)}_{#2}}
        }
}

\DeclareDocumentCommand\lambdavar{ g g }{%
        \IfNoValueTF {#1} {\bm{\lambda}} {
            \IfNoValueTF {#2} {{\bm{\lambda}^{(#1)}}{}}{\lambda^{(#1)}_{#2}}
        }
}

\DeclareDocumentCommand\tbetavar{ g g }{%
        \IfNoValueTF {#1} {{\bm{\tilde{\beta}}}} {
            \IfNoValueTF {#2} {{{\bm{\tilde{\beta}}}^{(#1)}}{}}{{{\tilde{\beta}}^{(#1)}_{#2}}}
        }
}

\DeclareDocumentCommand\alphavar{ g g }{%
        \IfNoValueTF {#1} {\bm{\alpha}} {
            \IfNoValueTF {#2} {{\bm{\alpha}^{(#1)}}}{\alpha^{(#1)}_{#2}}
        }
}

\DeclareDocumentCommand\hfunc{ g g }{%
        \IfNoValueTF {#1} {h} {
            \IfNoValueTF {#2} {{{h}^{(#1)}}}{h^{(#1)}_{#2}}
        }
}

\DeclareDocumentCommand\zcut{ g g }{%
        \IfNoValueTF {#1} {\bm{q}} {
            \IfNoValueTF {#2} {{\bm{q}^{(#1)}}}{q^{(#1)}_{#2}}
        }
}

\DeclareDocumentCommand\Zcut{ g g }{%
        \IfNoValueTF {#1} {\bm{Q}} {
            \IfNoValueTF {#2} {{\bm{Q}^{(#1)}}}{\bm{Q}^{(#1)}_{#2}}
        }
}

\DeclareDocumentCommand\xcut{ g g }{%
        \IfNoValueTF {#1} {\bm{h}} {
            \IfNoValueTF {#2} {{\bm{h}^{(#1)}}}{h^{(#1)}_{#2}}
        }
}

\DeclareDocumentCommand\Xcut{ g g }{%
        \IfNoValueTF {#1} {\bm{H}} {
            \IfNoValueTF {#2} {{\bm{H}^{(#1)}}}{\bm{H}^{(#1)}_{#2}}
        }
}

\DeclareDocumentCommand\hxcut{ g g }{%
        \IfNoValueTF {#1} {\bm{g}} {
            \IfNoValueTF {#2} {{\bm{g}^{(#1)}}}{g^{(#1)}_{#2}}
        }
}

\DeclareDocumentCommand\hXcut{ g g }{%
        \IfNoValueTF {#1} {\bm{G}} {
            \IfNoValueTF {#2} {{\bm{G}^{(#1)}}}{\bm{G}^{(#1)}_{#2}}
        }
}

\DeclareDocumentCommand\D{ g g }{%
        \IfNoValueTF {#1} {\mathbf{D}} {
            \IfNoValueTF {#2} {\mathbf{D}^{(#1)}}{\mathbf{D}^{(#1)}_{#2}}
        }
}

\DeclareDocumentCommand\A{ g g }{%
        \IfNoValueTF {#1} {\mathbf{A}} {
            \IfNoValueTF {#2} {\mathbf{A}^{(#1)}}{\mathbf{A}^{(#1)}_{#2}}
        }
}

\DeclareDocumentCommand\a{ g g }{%
        \IfNoValueTF {#1} {\mathbf{a}} {
            \IfNoValueTF {#2} {\mathbf{a}^{(#1)}}{{a}^{(#1)}_{#2}}
        }
}

\DeclareDocumentCommand\al{ g g }{%
        \IfNoValueTF {#1} {\underline{\mathbf{a}}} {
            \IfNoValueTF {#2} {\underline{\mathbf{a}}^{(#1)}}{\underline{a}^{(#1)}_{#2}}
        }
}

\DeclareDocumentCommand\au{ g g }{%
        \IfNoValueTF {#1} {\overline{\mathbf{a}}} {
            \IfNoValueTF {#2} {\overline{\mathbf{a}}^{(#1)}}{\overline{a}^{(#1)}_{#2}}
        }
}

\DeclareDocumentCommand\c{ g }{%
        \IfNoValueTF {#1} {\bm{c}} {
            {c^{({#1})}}
        }
}

\DeclareDocumentCommand\cl{ g }{%
        \IfNoValueTF {#1} {\underline{c}} {
            {\underline{c}^{({#1})}}
        }
}

\DeclareDocumentCommand\cu{ g }{%
        \IfNoValueTF {#1} {\overline{c}} {
            {\overline{c}^{({#1})}}
        }
}

\DeclareDocumentCommand\AA{ g g }{
        \IfNoValueTF {#1} {\mathbf{\Omega}} {
            \IfNoValueTF {#2} {\mathbf{\Omega}(#1, #1)}{\mathbf{\Omega}(#1, #2)}
        }
}

\DeclareDocumentCommand\S{ g g }{%
        \IfNoValueTF {#1} {\mathbf{S}} {
            \IfNoValueTF {#2} {\mathbf{S}^{(#1)}}{\mathbf{S}^{(#1)}_{#2}}
        }
}

\DeclareDocumentCommand\K{ g g }{%
        \IfNoValueTF {#1} {\mathbf{K}} {
            \IfNoValueTF {#2} {\mathbf{K}^{(#1)}}{\mathbf{K}^{(#1)}_{#2}}
        }
}

\DeclareDocumentCommand\B{ g g }{%
        \IfNoValueTF {#1} {\mathbf{B}} {
            \IfNoValueTF {#2} {\mathbf{B}^{(#1)}}{\mathbf{B}^{(#1)}_{#2}}
        }
}

\DeclareDocumentCommand\lowerb{ g g }{%
        \IfNoValueTF {#1} {{\mathbf{\underline{b}}}} {
            \IfNoValueTF {#2} {{\mathbf{\underline{b}}}^{(#1)}}{{\mathbf{\underline{b}}}^{(#1)}_{#2}}
        }
}

\DeclareDocumentCommand\z{ g g }{%
        \IfNoValueTF {#1} {\mathbf{z}} {
            \IfNoValueTF {#2} {\mathbf{z}^{(#1)}}{z^{(#1)}_{#2}}
        }
}

\DeclareDocumentCommand\hz{ g g }{%
        \IfNoValueTF {#1} {\hat{\mathbf{z}}} {
            \IfNoValueTF {#2} {\hat{\mathbf{z}}^{(#1)}}{\hat{z}^{(#1)}_{#2}}
        }
}

\DeclareDocumentCommand\hx{ g g }{%
        \IfNoValueTF {#1} {\hat{\boldsymbol{x}}} {
            \IfNoValueTF {#2} {\hat{\boldsymbol{x}}^{(#1)}}{\hat{x}^{(#1)}_{#2}}
        }
}

\DeclareDocumentCommand\x{ g g }{%
        \IfNoValueTF {#1} {\boldsymbol{x}} {
            \IfNoValueTF {#2} {\boldsymbol{x}^{(#1)}}{x^{(#1)}_{#2}}
        }
}

\DeclareDocumentCommand\bu{ g g }{%
        \IfNoValueTF {#1} {\boldsymbol{u}} {
            \IfNoValueTF {#2} {\boldsymbol{u}^{(#1)}}{{u}^{(#1)}_{#2}}
        }
}

\DeclareDocumentCommand\buh{ g g }{%
        \IfNoValueTF {#1} {\boldsymbol{u}^*} {
            \IfNoValueTF {#2} {\boldsymbol{{u}}^{*(#1)}}{{{u}}^{*(#1)}_{#2}}
        }
}

\DeclareDocumentCommand\bl{ g g }{%
        \IfNoValueTF {#1} {\boldsymbol{l}} {
            \IfNoValueTF {#2} {\boldsymbol{l}^{(#1)}}{{l}^{(#1)}_{#2}}
        }
}

\DeclareDocumentCommand\blh{ g g }{%
        \IfNoValueTF {#1} {\boldsymbol{l}^*} {
            \IfNoValueTF {#2} {\boldsymbol{l}^{*(#1)}}{{{l}}^{*(#1)}_{#2}}
        }
}

\DeclareDocumentCommand\aaa{ g }{%
        \IfNoValueTF {#1} {\bm{a}} {
            {\bm{a}^{({#1})}}
        }
}

\DeclareDocumentCommand\haaa{ g }{%
        \IfNoValueTF {#1} {\bm{\hat{a}}} {
            {\bm{\hat{a}}^{({#1})}}
        }
}

\DeclareDocumentCommand\bbb{ g g }{%
        \IfNoValueTF {#1} {\mathbf{P}} {
            \IfNoValueTF {#2} {{\mathbf{P}_{#1}}}{{\mathbf{P}_{#1}^{({#2})}}}
        }
}

\DeclareDocumentCommand\hbbb{ g g }{%
        \IfNoValueTF {#1} {\mathbf{\hat{P}}} {
            \IfNoValueTF {#2} {{\mathbf{\hat{P}}_{#1}}}{{\mathbf{\hat{P}}_{#1}^{({#2})}}}
        }
}

\DeclareDocumentCommand\ccc{ g g }{%
        \IfNoValueTF {#1} {\mathbf{q}} {
            \IfNoValueTF {#2} {{\mathbf{q}_{#1}}}{{\mathbf{q}_{#1}^{(#2)}}{}}
        }
}

\DeclareDocumentCommand\constc{ g }{%
        \IfNoValueTF {#1} {c} {
            {c^{({#1})}}
        }
}

\DeclareDocumentCommand\setz{ g g }{%
        \IfNoValueTF {#1} {\mathcal{Z}} {
            \IfNoValueTF {#2} {\mathcal{Z}^{(#1)}}{\mathcal{Z}^{(#1)}_{#2}}
        }
}

\DeclareDocumentCommand\setzp{ g g }{%
        \IfNoValueTF {#1} {\mathcal{Z^+}} {
            \IfNoValueTF {#2} {\mathcal{Z}^{+(#1)}}{\mathcal{Z}^{+(#1)}_{#2}}
        }
}

\DeclareDocumentCommand\setzn{ g g }{%
        \IfNoValueTF {#1} {\mathcal{Z^-}} {
            \IfNoValueTF {#2} {\mathcal{Z}^{-(#1)}}{\mathcal{Z}^{-(#1)}_{#2}}
        }
}

\DeclareDocumentCommand\tsetz{ g g }{%
        \IfNoValueTF {#1} {\tilde{\mathcal{Z}}} {
            \IfNoValueTF {#2} {\tilde{\mathcal{Z}}^{(#1)}}{\tilde{\mathcal{Z}}^{(#1)}_{#2}}
        }
}

\DeclareDocumentCommand\seti{ g g }{%
        \IfNoValueTF {#1} {\mathcal{I}} {
            \IfNoValueTF {#2} {\mathcal{I}^{(#1)}}{\mathcal{I}^{(#1)}_{#2}}
        }
}

\DeclareDocumentCommand\setip{ g g }{%
        \IfNoValueTF {#1} {\mathcal{I}^{+}} {
            \IfNoValueTF {#2} {\mathcal{I}^{+(#1)}}{\mathcal{I}^{+(#1)}_{#2}}
        }
}

\DeclareDocumentCommand\setin{ g g }{%
        \IfNoValueTF {#1} {\mathcal{I}^{-}} {
            \IfNoValueTF {#2} {\mathcal{I}^{-(#1)}}{\mathcal{I}^{-(#1)}_{#2}}
        }
}

\DeclareDocumentCommand\tseti{ g g }{%
        \IfNoValueTF {#1} {\tilde{\mathcal{I}}} {
            \IfNoValueTF {#2} {\tilde{\mathcal{I}}^{(#1)}}{\tilde{\mathcal{I}}^{(#1)}_{#2}}
        }
}

\DeclareDocumentCommand\tz{ g g }{%
        \IfNoValueTF {#1} {\tilde{z}} {
            \IfNoValueTF {#2} {\tilde{z}^{(#1)}}{\tilde{z}^{(#1)}_{#2}}
        }
}

\DeclareDocumentCommand\f{ g g }{%
        \IfNoValueTF {#1} {f} {
            \IfNoValueTF {#2} {f^{(#1)}}{f^{(#1)}_{#2}}
        }
}

\DeclareDocumentCommand\lf{ g g }{%
        \IfNoValueTF {#1} {\underline{f}} {
            \IfNoValueTF {#2} {\underline{f}^{(#1)}}{\underline{f}^{(#1)}_{#2}}
        }
}

\def\eqref#1{(\ref{#1})}

\def\1{\bm{1}}

\DeclareMathAlphabet{\mathsfit}{\encodingdefault}{\sfdefault}{m}{sl}
\SetMathAlphabet{\mathsfit}{bold}{\encodingdefault}{\sfdefault}{bx}{n}

\usepackage{hyperref}
\usepackage[T1]{fontenc}
\usepackage[capitalize,noabbrev]{cleveref}
\hypersetup{hidelinks}

\theoremstyle{plain}
\newtheorem{theorem}{Theorem}[section]
\newtheorem{proposition}[theorem]{Proposition}
\newtheorem{lemma}[theorem]{Lemma}

\theoremstyle{definition}

\newtheorem{assumption}[theorem]{Assumption}
\theoremstyle{remark}

\crefname{assumption}{Assumption}{Assumptions}
\Crefname{assumption}{Assumption}{Assumptions}

\newcommand{\GEAR}{\textsc{GEAR}}

\newcommand{\SBReweight}{$\lambda$-\textsc{SB}}

\title{GEAR: A GPU-Accelerated Global Solver for Bounded Continuous Nonlinear Programs \\ via Linear Bound Propagation}

\author{%
  Duo Zhou, Hesun Chen, Xiangru Zhong, Grani A. Hanasusanto, Huan Zhang
   \\
  University of Illinois Urbana-Champaign\\
  \texttt{duozhou2@illinois.edu, huan@huan-zhang.com} \\
}

\begin{document}
\raggedbottom
\maketitle
\addtocontents{toc}{\protect\setcounter{tocdepth}{-1}}

\setlength{\abovedisplayskip}{4pt plus 2pt minus 2pt}
\setlength{\belowdisplayskip}{4pt plus 2pt minus 2pt}
\makeatletter\@secpenalty=0
\let\GEARparagraph\paragraph
\renewcommand{\paragraph}{\@startsection{paragraph}{4}{\z@}{0.8ex \@plus 0.3ex \@minus 0.2ex}{-1em}{\normalsize\bf}}
\makeatother
\begin{abstract}
    We present a novel GPU-accelerated global solver for constrained nonlinear programming (NLP) inspired by the efficient linear bound propagation framework in neural network (NN) verification.
    Existing work on NN verification, a problem that can be cast as feasibility detection for a nonlinear objective, demonstrated orders-of-magnitude speedups with the linear bound propagation-based framework when compared to traditional solvers. 
    Extending this paradigm to global nonlinear programming requires addressing many challenges: a global constrained NLP solver must effectively exploit all, possibly unstructured, nonlinear constraints to improve the objective, whereas NN verification typically handles simple input box constraints over a well-structured NN and does not directly optimize an objective function. 
    We address these challenges by fully exploiting the linear bounds produced by bound propagation: these linear bounds can be viewed as relaxations of nonlinear constraints and used to provide lower bounds, check infeasibility, and guide the search for a primal solution.
    In addition, we propose tighter linear relaxations for bilinear, quotient, and other commonly used nonlinear functions in NLP problems, and a projected augmented Lagrangian method tightly coupled to our dual solving procedure to produce primal solutions.
    Our method solves 421 of 505 bounded and continuous NLP problems in the GAMS and MINLPLib benchmarks under a 180-second limit, more than each of the six global solvers we compare with, including  Gurobi and BARON (commercial), and SCIP (open-source). Sixteen of these problems are certified by none of the six global solvers; all of them contain exponential, logarithmic, or trigonometric terms, which are often challenging for traditional approaches. 
    Our method demonstrates a distinct, highly complementary regime for GPU-accelerated global optimization of bounded continuous nonlinear programs. Project webpage is at \url{https://verified-intelligence.github.io/GEAR}.
\end{abstract}

\section{Introduction}
\label{sec:intro}

Neural-network (NN) verification asks whether all outputs of a computation graph
satisfy a prescribed property over a bounded input set, and can be cast as
feasibility detection for a nonlinear objective~\cite{brown2022unifiedviewsdpbasedneural, fazlyab2020safety, raghunathan2018semidefinite, lan2022tight, chiu2025sdpcrown, dathathri2020enabling}. Many modern NN verifiers, such as $\alpha,\!\beta$-CROWN, solve this
through \emph{linear bound propagation}: linear coefficients are propagated
backward through the graph to obtain valid bounds on the output, and the
coefficient updates batch as dense tensor operations on GPUs. This pipeline has
delivered orders-of-magnitude speedups over traditional verification approaches
built on LP, SDP, or MILP solvers
\citep{zhang2018crown,xu2020automatic,wang2021betaCROWN,xu2021fast,zhang2022general,wong2020neural, venzke2020verification, sun2022romax, chen2024verification, yang2024lyapunov, wu2023neural}.
The same paradigm could in principle drive a general global solver for nonlinear
programs (NLPs) with continuous, box-bounded variables, where linear bound
propagation through the computation graph would supply valid bounds for thousands of problems in a single batched workload~\cite{bunel2018unified,bunel2020branch,singh2018fast,wang2018efficient,bak2021nfm}.

Extending linear bound propagation from NN verification to global NLP is,
however, far from automatic. A constrained NLP carries an objective together
with many, possibly unstructured, nonlinear equality and inequality constraints
over operators such as products, quotients, powers, logarithms, square roots,
and trigonometric terms; a global solver must exploit \emph{all} of these
constraints simultaneously to drive the objective down, certify optimality, and
recover feasible incumbents~\cite{biegler2010nonlinear}. NN verification, by contrast, typically operates
over a well-structured neural network and a simple input box, and never
directly optimizes an objective. Classical global NLP solvers
\citep{zhang2025baron,belotti2009branching,vigerske2018scip,bestuzheva2023scip,misener2014antigone}
implement spatial branch-and-bound on top of per-node relaxations,
McCormick envelopes, optimization-based bound tightening, cut pools, and
primal heuristics. These relaxations are mostly polyhedral, although some methods use convex
nonlinear relaxations~\citep{adjiman1998global,khajavirad2018hybrid} and, for nonconvex quadratic
programs, convex quadratic or semidefinite relaxations~\citep{nohra2021spectral,burer2008finite}.
Each subdomain may demand its own relaxation solves and
nonlinear bookkeeping, which makes thousands of subdomains and constraint
outputs hard to evaluate as one GPU-batched workload.

We introduce \GEAR{}, a GPU-accelerated global solver for bounded continuous
NLPs that closes this gap by \emph{fully exploiting} the linear bounds produced
by bound propagation. The same valid linear coefficients are read
simultaneously as: (i) linear relaxations of the original nonlinear
constraints, which also clip the input and intermediate domains; (ii) a Lagrangian combination of these linear bounds that yields a
constraint-aware lower bound on the objective; (iii) an infeasibility certificate when some nonnegative
combination of the relaxations is violated everywhere on the current box; and (iv) a
primal seed and a branching score that guide the spatial branch-and-bound
search. The objective and all constraints are vectorized as outputs of one
computation graph, so linear bound propagation returns coefficient tensors for
every constraint across thousands of branch-and-bound subproblems, and the
downstream lower bounding, infeasibility testing, primal initialization, and
branching all execute as batched tensor operations rather than as a sequence of
per-node LP solves.  The key observation is that the linear bounds themselves relax every constraint, and their Lagrangian combination is an LP
dual that \GEAR{} evaluates in closed form per box.  The price, a looser bound per box than an LP with
all McCormick planes, is offset in the tree by bounding many boxes per GPU step, which helps on highly
nonlinear problems of low and moderate dimension rather than on large structured models (\cref{sec:exp}).\looseness=-1

Two further ingredients make this pipeline effective. First, since the
tightness of linear bound propagation is determined largely by the per-operator
linear relaxations, we develop tighter ones for the operators that
dominate NLP but are rare in NN verification. In particular, for the bilinear
products that pervade NLP benchmarks, we choose among the McCormick planes by
their value on the current box, and we use secant-tangent bounds that follow
each operator's curvature and direct bounds for quotients, reciprocals, powers,
and trigonometric terms, all designed to keep bound propagation tensorized.
Second, we propose a projected augmented Lagrangian method (ALM) for primal
recovery that is tightly coupled to the dual procedure: the minimizer of the
Lagrangian bound starts the primal search, its multipliers initialize the ALM
multipliers, and first-order primal updates batch across thousands of
nodes instead of running an SQP or interior-point subproblem sequence.

On 505 bounded continuous NLP instances drawn from the \texttt{gams\_world} and \texttt{minlplib}
benchmarks under a 180-second wall-clock limit, feasibility tolerance $10^{-6}$ and relative gap
$10^{-4}$, \GEAR{} solves 421 instances, more than each of the six global
solvers we compare with, including the commercial Gurobi (401) and BARON and
the open-source SCIP. Against the strongest baseline, Gurobi, \GEAR{} solves
38 instances that Gurobi does not (Gurobi solves 18 that \GEAR{} does not), and
16 instances are certified by none of the six global solvers; all of them
contain exponential, logarithmic or trigonometric terms, whose compositions
are often out of reach for classical LP-relaxation-based pipelines. \GEAR{} therefore
complements CPU global solvers with a GPU-accelerated regime for global optimization of bounded
continuous nonlinear programs; mixed-integer and unbounded programs are outside the scope of this paper.
Our contributions:
\begin{itemize}[leftmargin=10pt,labelsep=0.5em,topsep=0pt,itemsep=1pt,parsep=0pt]
    \item \textbf{Linear bound propagation as a unified interface for global NLP}: Inspired by the success of linear bound propagation in neural network verification, we demonstrate its success as a new GPU-accelerated regime for bounded continuous NLP. One set of linear
    bounds serves every component, and computing them requires no LP solve (\cref{subsec:representation,subsec:dual,subsec:branch}).
    \item \textbf{Tighter operator-specific linear relaxations} for the bilinear, quotient, power and
    trigonometric operators of NLP (\cref{sec:operators}), essential to make linear bound propagation successful for NLPs.
    \item \textbf{A projected ALM primal solver tightly coupled to the dual solver} and guided by the linear bounds (\cref{subsec:primal}).
    \item \textbf{An open-source complementary GPU-batched solver for bounded continuous global NLP}, evaluated on 505 GAMS and MINLPLib instances with promising results against six global solvers. 
\end{itemize}


\vspace{-6pt}
\section{Preliminaries}
\label{sec:setup}
\nobreak\vspace{-5pt}

\paragraph{Bounded Continuous Nonlinear Programs}
\label{subsec:graph_nlp}
\label{subsec:problem_subproblem}

We consider bounded continuous nonlinear programs (NLPs)
\begin{equation}
p^\star=\min_{x\in B_0} f_0(x)
\quad\text{s.t.}\quad f_i(x)\le 0,\quad i=1,\ldots,m,
\label{eq:nlp_standard_form}
\end{equation}
where $f_0$ is the objective, $f_1,\ldots,f_m$ are the constraint functions,
and $B_0\subset\mathbb R^n$ is a box with finite bounds. Equalities and
two-sided inequalities of the source model are written as pairs of
one-sided inequalities, and a feasibility tolerance $\tau_{\mathrm{feas}}\ge0$
is absorbed into every $f_i$; for example, $g(x)=0$ becomes
$g(x)-\tau_{\mathrm{feas}}\le0$ and $-g(x)-\tau_{\mathrm{feas}}\le0$. Lower
bounds and feasibility checks therefore refer to the same
problem~\eqref{eq:nlp_standard_form}.
Each $f_i$ is factorable over the operator library of
\cref{subsec:operators}; integer variables, unbounded variables, and
unsupported operators are outside our scope.

The objective and all constraints can be written as the outputs of one
computation graph,
\begin{equation}
F(x)=\bigl(f_0(x),f_1(x),\ldots,f_m(x)\bigr),
\label{eq:nlp_graph}
\end{equation}
in which every intermediate variable $z$ is produced by one operator,
$z=h(y)$, applied to the input $x$ or to earlier intermediate variables $y$
(one variable, or two for binary operators such as products). For a box
$B=[\underline x,\overline x]\subseteq B_0$, the feasible set and optimal value of the
subproblem on $B$ are
\begin{equation}
\mathcal F(B)=\{x\in B:\ f_i(x)\le0,\ i=1,\ldots,m\},
\qquad
p^\star(B)=\min_{x\in\mathcal F(B)} f_0(x),
\label{eq:local_feasible_value}
\end{equation}
with $p^\star(B)=+\infty$ if $\mathcal F(B)=\emptyset$.

\paragraph{Neural Network Verification with Linear Bound Propagation}
\label{subsec:lbp}

Given a scalar output $f$ of a computation graph and an input box
$B=[\underline x,\overline x]$, neural network (NN) verification asks
whether $\min_{x\in B}f(x)>0$
\citep{wang2018formal,wang2018efficient,bunel2018unified}; in its basic
form this is a lower-bounding problem over a box, and constraint-aware
verifiers additionally incorporate linear constraints such as activation
splits and cutting planes
\citep{wang2021betaCROWN,zhang2022general,zhou2025clip}. Linear bound
propagation (LBP)
\citep{zhang2018crown,xu2020automatic,wang2021betaCROWN,xu2021fast,zhang2022general}
answers the basic problem as follows; we follow the notation of
\cite{wang2021betaCROWN} and \cite{zhou2025clip}, with underlines and
overlines for lower and upper bounds. For each nonlinear operator $z=h(y)$,
let $Y=[\underline y,\overline y]$ contain every value of $y$ attained for
$x\in B$ (the \emph{intermediate bounds}), and let $\underline h,\overline h$
be a linear relaxation of $h$ on $Y$:
\begin{equation}
\underline h(y)\le h(y)\le\overline h(y)\qquad\forall y\in Y.
\label{eq:operator_relaxation}
\end{equation}
Starting from $f$, LBP maintains a linear expression in the intermediate
variables that bounds $f$ from below and substitutes these relaxations from
the output toward the input: a term $\gamma z$ with $z=h(y)$ becomes
$\gamma\,\underline h(y)$ if $\gamma\ge0$ and $\gamma\,\overline h(y)$
otherwise, which preserves the inequality; linear operators are substituted
exactly. When only the input remains, the result is a linear lower bound
(an upper bound can be obtained symmetrically)
\begin{equation}
\underline f(x)=a^\top x+c\;\le\;f(x)\qquad\forall x\in B.
\label{eq:row_per_output}
\end{equation}
Intermediate bounds come from running LBP with each intermediate variable as the
output, intersected with interval bounds. A linear function attains its minimum over a box at a vertex
selected coordinate-wise by the signs of its coefficients $a_j$, so the
scalar bound, called the \emph{concretization} of~\eqref{eq:row_per_output}, is
closed-form:
\begin{equation}
\min_{x\in B}\underline f(x)=c+\sum_{j=1}^n\min\{a_j\underline x_j,\ a_j\overline x_j\}.
\label{eq:lbp_corner}
\end{equation}
Each propagation step is a product of coefficient tensors, so propagation from
many outputs over many boxes runs as one batched GPU computation.
Besides producing this scalar, the coefficients $(a,c)$ describe how the
bound varies across $B$; \GEAR{} keeps them and uses them in all of
its components (\cref{subsec:representation}).
Complete verifiers wrap LBP in branch-and-bound that splits either the
input box (\emph{input BaB}) or activations (\emph{activation BaB})
\citep{wang2021betaCROWN,zhou2025clip}.

\paragraph{Spatial Branch-and-Bound}
\label{subsec:sbab}

In its usual rectangular form, spatial branch-and-bound (sBaB) maintains a queue $\mathcal Q$ of boxes,
initialized with $B_0$, and an incumbent value $\mathrm{UB}=f_0(\bar x)$ of
the best feasible point $\bar x$ found so far by a \emph{primal} search
($\mathrm{UB}=+\infty$ initially). For each box $B$ it computes a
\emph{dual} (lower) bound $\mathrm{LB}(B)\le p^\star(B)$ and removes $B$ if $\mathcal F(B)=\emptyset$
is certified (\emph{infeasibility pruning}) or if
$\mathrm{LB}(B)\ge\mathrm{UB}$ (\emph{dominance pruning}). A box may first be
replaced by a smaller box that still contains $\mathcal F(B)$
(\emph{bound tightening}); a box that survives is split along one
coordinate into two children that cover it. At all times
$\min\{\mathrm{UB},\min_{B\in\mathcal Q}\mathrm{LB}(B)\}\le p^\star$; the
search stops when this bound is within $\varepsilon$ of $\mathrm{UB}$, or when
$\mathcal Q$ is empty, in which case $\mathrm{UB}$ is optimal (or the problem
is infeasible if $\mathrm{UB}=+\infty$). Validity requires only that lower
bounds, infeasibility certificates, and tightening steps be valid and that
incumbents be checked for feasibility; the splitting rule and the primal
search affect speed only.\looseness=-1

Many LP-based global solvers
\citep{zhang2025baron,belotti2009branching,vigerske2018scip} branch on variable
bounds, so each node is a box, but they obtain $\mathrm{LB}(B)$ over a polyhedron: an LP relaxation
that may contain auxiliary variables, the model's linear constraints, operator relaxations, and
cutting planes. LBP obtains~\eqref{eq:lbp_corner} without an LP
solve and batches across boxes, but it minimizes over the box alone and by itself does not account
for the constraints of~\eqref{eq:nlp_standard_form}. \Cref{sec:method} shows
how to obtain constraint-aware lower bounds and infeasibility certificates
from LBP alone, and how to pair them with a primal search and a branching
rule, keeping all computations batched.


\vspace{-6pt}
\section{GEAR: Global NLP Solving with GPU-Accelerated Linear Bounds}
\label{sec:method}
\nobreak\vspace{-5pt}

\GEAR{} extends linear bound propagation from verification to constrained global optimization.
The central challenge is to make nonlinear constraints useful throughout the search, not merely for checking feasibility, but also for guiding lower bounding and feasible solution finding, while retaining GPU-batched computation.

\begin{figure}[t]
\centering
\includegraphics[width=\linewidth]{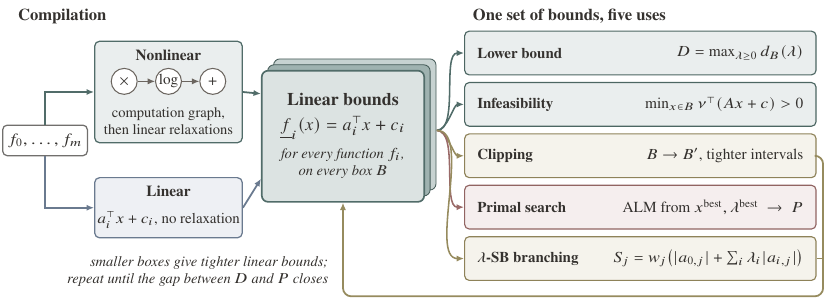}
\caption{\textbf{Overview of \GEAR{}.}  Each of $f_0,\dots,f_m$ is nonlinear or linear.  Nonlinear functions are
represented as computation graphs, and linear bound propagation applies linear relaxations to them; linear
functions need no relaxation.  Either way each function has one linear bound on each box~$B$.  These bounds give the
lower bound $D$, the infeasibility test and clipping (\cref{subsec:dual}), the start of the projected ALM search for
the incumbent $P$ (\cref{subsec:primal}), and the branching score (\cref{subsec:branch}).  Smaller boxes give tighter
linear bounds, and an instance is solved once the relative gap of \cref{eq:exp_solver_criterion} is at most $10^{-4}$.}
\label{fig:walkthrough}
\vspace{-5mm}
\end{figure}

\Cref{fig:walkthrough} connects the parts of this section: the graph representation
(\cref{subsec:representation}), the lower bounds built from linear bounds (\cref{subsec:dual}), the
projected ALM primal search (\cref{subsec:primal}) and the batched branch-and-bound
(\cref{subsec:branch}). Presolve and the remaining
components of the full solver only add valid bounds or verified points (\cref{fig:components},
\cref{appsec:components}).

\vspace{-5pt}
\subsection{Why Verification Alone Is Insufficient for Constrained Global NLP}
\label{subsec:naive}
\label{subsec:gaps}
\nobreak\vspace{-3pt}

Both global NLP and NN verification problems require valid bounds, but differ in what
must be found and certified. A verification problem asks whether a property
holds over an input set; it can finish upon finding a counterexample (SAT) or
verifying its properties (UNSAT). However, constrained optimization must instead find a point satisfying every constraint (feasible solution) and verify that no feasible point improves
its objective by more than the prescribed tolerance (optimal). Its pruning threshold
is therefore not fixed: it is the objective value of an incumbent that
the solver must produce and improve. Constraint-aware verifiers already
incorporate linear inequalities from specifications, activation splits,
and cutting planes \citep{wang2021betaCROWN,zhang2022general,zhou2025clip}.
The challenge is to extend this machinery into an objective-directed
primal-dual search over general nonlinear constraints.

This gap becomes apparent on both sides of the primal-dual search. 
A search that minimizes $f_0$ over the box alone can return points that violate $f_i\le0$.
A lower bound computed from $f_0$ alone minimizes over the whole box, including infeasible points,
so it can stay below $p^\star$ however far we branch.
For example, minimizing $x$ on $[0,1]$ subject to $x\ge1/2$ has optimum
$1/2$, but every box containing $0$ retains an objective-only lower bound
of $0$. More branching cannot remove that box by dominance; a constraint
must rule it out or strengthen its bound. Conversely, a strong lower bound
cannot drive dominance pruning without a feasible incumbent. Global NLP
therefore needs both constraint-aware certificates and a search that
improves feasible objective values, not merely a counterexample finder.

Both remedies must stay batched: the solver has to reuse constraint information across many
subproblems without building a separate optimization model per box, keep linear constraints exact,
and relax products, quotients and powers tightly enough for bounding to be useful.  This motivates
the shared representation of \cref{subsec:representation} and the operator relaxations of
\cref{sec:operators}.


\vspace{-5pt}
\subsection{A GPU-Compatible Representation of Constrained NLPs}
\label{subsec:representation}
\label{subsec:vectorize}
\label{subsec:batched_bounds}
\nobreak\vspace{-3pt}

We begin by organizing the objective and constraints into a common computation graph, allowing bounding and primal search to reuse the same operations across many boxes.
\GEAR{} compiles them as the outputs of $F$ in~\eqref{eq:nlp_graph}, with the first
output designated as the objective. Forward evaluation, automatic
differentiation, and LBP operate on this same graph, and shared
subexpressions can be reused when exposed by the parser. The principal
parallel dimension is the batch of boxes: each graph operation processes
many subproblems with identical structure but different bounds, and all
outputs are bounded by the same linear bound propagation.

Linear constraints need no relaxation: LBP substitutes linear operators exactly, so the linear bound of
$f_i(x)=a_i^\top x+c_i$ is the constraint itself, with the
feasibility tolerance already incorporated as in \cref{subsec:graph_nlp}.
These coefficients are independent of the current box, and when the parser exposes them as an explicit
matrix, the solver reuses that matrix across subproblems whenever it constructs its constraint relaxation,
tightens bounds, or evaluates candidates.


Beyond a shared graph, the components need a common representation of the objective and constraints
on each box, which the linear bounds of propagation provide.  On each box $B$, linear bound
propagation provides
\begin{equation}
\underline f_i(x)=a_i^\top x+c_i\le f_i(x)
\qquad\forall x\in\mathcal F(B),\quad i=0,\ldots,m.
\label{eq:method_linear_bounds}
\end{equation}
Ordinary LBP gives validity throughout $B$; validity on $\mathcal F(B)$
also accommodates the constraint-tightened propagation introduced below.
Across $K$ boxes, the coefficients have logical shapes
$K\times(m+1)\times n$ and $K\times(m+1)$. The same linear bounds support objective bounding, infeasibility
certificates, clipping, primal initialization, and coordinate scoring
without rebuilding a separate optimization model for each operation.
A linear bound that cannot reject a whole box may still exclude part of it,
tighten an intermediate interval, or become useful in combination with
other linear bounds (\cref{app:batched-linear-bound-state}).

\vspace{-5pt}
\subsection{Strengthening Dual Bounds using Linear Bounds}
\label{subsec:dual}
\label{subsec:objective_bound}
\label{subsec:constraint_relaxations}
\label{subsec:lower_bounds}
\label{subsec:dual_ascent}
\label{subsec:valid_dual}
\label{subsec:infeasibility}
\label{subsec:relaxation_clipping}
\label{subsec:complete_clip}
\label{subsec:operators}
\label{subsec:mccormick}
\label{sec:operators}
\nobreak\vspace{-3pt}

\paragraph{Constraint-aware objective bounds.}
Although the shared linear bounds retain constraint information, ordinary box concretization still ignores feasibility and may minimize the objective over infeasible points.
We therefore incorporate the constraints directly into objective bounding through Lagrangian dualization.
Every feasible point satisfies
$\underline f_i(x)\le f_i(x)\le0$ for $i\ge1$. Stacking these linear bounds, together with any linear cuts, as
$A=[a_1,\ldots,a_m]^\top$ and $c=(c_1,\ldots,c_m)^\top$ defines
\begin{equation}
\mathcal P(B)=\{x\in B:Ax+c\le0\}\supseteq\mathcal F(B),
\qquad
p_{\mathrm{LP}}(B)=\min_{x\in\mathcal P(B)}\underline f_0(x)
\le p^\star(B).
\label{eq:node_lp}
\end{equation}
Instead of solving this LP separately at every box, \GEAR{} dualizes
$Ax+c\le0$ while retaining the box as the inner minimization domain.
For $\lambda\ge0$, write $a(\lambda)=a_0+A^\top\lambda$; then
\begin{align}
d_B(\lambda)
&=\min_{x\in B}\bigl[\underline f_0(x)+\lambda^\top(Ax+c)\bigr]=c_0+c^\top\lambda+
\sum_{j=1}^n\min\{a_j(\lambda)\underline x_j,
                  a_j(\lambda)\overline x_j\}.
\label{eq:dual_lb}
\end{align}
We call $d_B(\lambda)$ the \emph{Lagrangian bound}.
At feasible points the added terms are nonpositive, preserving a lower
bound on $f_0$; at infeasible minimizers of the objective's linear bound they can be
positive, allowing the bound to rise. In the example of
\cref{subsec:naive}, $\lambda=1$ makes
$x+\lambda(1/2-x)=1/2$ throughout the box, recovering the constrained
optimum. More generally, every $\lambda\ge0$ gives
$d_B(\lambda)\le p_{\mathrm{LP}}(B)$, and LP strong duality gives
$\max_{\lambda\ge0}d_B(\lambda)=p_{\mathrm{LP}}(B)$ when the relaxation is
feasible. Note that the multipliers $\lambda$ are thus dual variables of this LP, not KKT multipliers of the
NLP: optimizing them closes the gap to the LP over the current linear bounds, not the gap between that
LP and the NLP.  Each evaluation reuses the existing linear bounds and the concretization of
\cref{eq:lbp_corner}, costing $O(mn)$ per box without rerunning bound propagation.
Since $Ax_\lambda+c$ is a supergradient at an inner
minimizer $x_\lambda$, \GEAR{} takes a fixed budget of projected Adam ascent steps and keeps the
best evaluated bound, so an incomplete ascent affects tightness, not validity.
Because a fixed budget of Adam steps can stop far from the LP optimum, \GEAR{} also proposes multipliers with batched PDHG~\cite{applegate2021practical} on the GPU; every $\lambda\ge0$ gives a valid bound, so the best one is kept (\cref{appsec:comp_bab}). 
The constraints also support incumbent-independent infeasibility pruning:
\begin{equation}
\min_{x\in B}\nu^\top(Ax+c)>0,\qquad \nu\ge0,
\quad\Longrightarrow\quad\mathcal F(B)=\emptyset.
\label{eq:infeas_cert}
\end{equation}
\GEAR{} applies this test to every box after propagation. Unit vectors $\nu=e_i$ need no search: a box is
discarded when the linear bound of one constraint is positive on all of $B$. Constraints that are only
jointly inconsistent need a nonnegative combination: for the boxes with the lowest bounds, $\nu$ is taken from
the multipliers of the LP that minimizes the largest violation of $Ax+c\le0$ over $B$, solved on the CPU, and
\eqref{eq:infeas_cert} is then evaluated on the GPU (\cref{appsec:comp_bab}).
Duality, update rules, and certificate derivations
are detailed in \cref{app:dual-proofs,app:sound-invariants}.

\paragraph{Constraint-driven clipping to shrink the feasible region.}

Dual optimization strengthens the use of existing linear bounds, but cannot recover information lost when those bounds were constructed.
To address this limitation, \GEAR{} feeds constraint information back into propagation by tightening the input and intermediate domains via domain clipping~\cite{zhou2025clip}.
Relaxed clipping derives coordinate-wise input bounds, replacing $B$ by
$B'$ with
$\mathcal F(B)\subseteq\mathcal P(B)\subseteq B'\subseteq B$; it need not
compute the smallest enclosing box. Complete clipping instead exploits
$\mathcal P(B)$ when bounding an intermediate variable in the computation graph. For a propagated
linear lower bound $\underline z(x)$, it targets
$\min_{x\in\mathcal P(B)}\underline z(x)$ through a dual of the same form
as~\eqref{eq:dual_lb}, with separate multipliers; upper bounds are handled
symmetrically. While \cite{zhou2025clip} uses coordinate descent to handle multiple constraints, which works well when constraints are few in NN verification settings, it becomes slow as the number of NLP constraints grows. To address this, we use gradient descent to optimize the dual variables over the active constraints.


Narrower intervals give tighter linear relaxations of nonlinear operators and stronger linear bounds for the objective and
constraints, and an empty interval certifies infeasibility.  Clipping uses linear bounds already
available from the parent or an earlier propagation, and the tightened intervals need only contain values
attained at feasible inputs, which is why \eqref{eq:method_linear_bounds} states validity on
$\mathcal F(B)$.  The Lagrangian bound thus extracts more from the existing linear bounds, while
clipping uses them to improve the domains on which later linear bounds are built
(\cref{appendix::clipping,app:node-lifecycle}).

\paragraph{Tighter nonlinear operator relaxations.}

Clipping tightens the domains of nonlinear operators, but the resulting bounds can still be loose when the operator relaxations themselves are weak.
\GEAR{} therefore strengthens these local relaxations while keeping one linear bound per output.
For a product, McCormick's construction gives two lower and two
upper linear planes on a rectangle \citep{McCormick1976}. A lifted LP can
retain all four, while backward linear bound propagation selects a plane according to the
sign of its coefficient. \GEAR{} compares the relevant candidates by
closed-form concretization: for a nonnegative coefficient, it chooses the
lower plane with the largest box minimum; for a negative coefficient,
the upper plane with the smallest box maximum. This maximizes the local contribution among the candidates without
iterative relaxation-parameter optimization. For the product example in
\cref{fig:envelope_binary_ops}, the two lower facets concretize to $-10$
and $2$ on the same rectangle, illustrating that a different valid choice
can matter without adding more stored planes. For quotients,
\GEAR{} constructs planes directly for the two-input operator instead of
composing reciprocal and product relaxations. Candidate planes are used
only on denominator intervals excluding zero and only after their
validity is established over the entire rectangle. This avoids an
unnecessary intermediate relaxation while retaining the same linear
interface. Supported univariate operators use secant-tangent constructions,
with the tangent slope matched to the secant slope on fixed-curvature
intervals. Together, these changes improve the linear bounds that the Lagrangian bound and
clipping use, rather than only re-weighting them. However, the relaxation can still be weaker than an LP
that keeps all McCormick planes at once; \GEAR{} trades that strength for
batched propagation and relies on spatial refinement to reduce the
remaining error.
Operator formulas, domain guards, and comparisons are in
\cref{app:operators,fig:envelope_binary_ops}.\looseness=-1

\vspace{-5pt}
\subsection{Scalable Primal Algorithm on GPU with Projected ALM and Linear Bounds}
\label{subsec:primal}
\label{subsec:residual_primal_branch}
\label{subsec:heuristic_use}
\nobreak\vspace{-3pt}

The mechanisms above strengthen the dual side of global optimization, but even a tight lower bound cannot enable dominance pruning without a feasible incumbent.
To efficiently find such incumbents across many boxes, \GEAR{} introduces a GPU-batched projected augmented Lagrangian method that jointly pursues objective improvement and constraint satisfaction.
The choice is computational as well as algorithmic:
fixed-budget first-order steps need only graph evaluations, automatic
differentiation, and elementwise projection. They batch across searches
without the per-iteration KKT factorization used by conventional
interior-point methods \citep{wachter2006implementation}. This makes
the projected ALM compatible with the many-box workload of the bounding procedure,
rather than requiring a separate nonlinear solver invocation per box.

\paragraph{Adaptive feasibility-objective balancing.}
The augmented Lagrangian keeps objective improvement and feasibility
restoration in the same differentiable search. For $\mu\ge0$ and $\rho>0$,
\GEAR{} uses the inequality-form function
\citep{rockafellar1973augmented}
\begin{equation}
\phi_\rho(x;\mu)=f_0(x)
+\frac{\rho}{2}\sum_{i=1}^m[f_i(x)+\mu_i/\rho]_+^2
-\frac{\|\mu\|_2^2}{2\rho},
\label{eq:hpr_alm}
\end{equation}
where $[t]_+=\max\{t,0\}$. Each outer iteration applies a fixed number
of Adam descent steps projected onto $B$, then updates
$\mu\leftarrow[\mu+\rho f(x)]_+$, with $f=(f_1,\ldots,f_m)$.
The penalty increases when the maximum violation fails to decrease
sufficiently. In the inner gradient, constraint $i$ receives the effective
weight $[\mu_i+\rho f_i(x)]_+$. The multiplier therefore carries constraint
emphasis across outer iterations, while the current violation adjusts
that emphasis within an iteration. This is the mechanism by which ALM
couples objective progress to feasibility restoration rather than relying
on one fixed penalty coefficient. Gradients use the original objective and constraints,
including the exact linear constraints, rather than replacing the NLP by its
relaxation. These updates seek improved incumbents; they are not relied
upon to certify a lower bound or to solve every local subproblem.

\paragraph{How the dual helps the primal.}
Although $\lambda$ weights the linear bounds in the Lagrangian
bound~\eqref{eq:dual_lb} whereas $\mu$ weights the exact nonlinear
constraints in~\eqref{eq:hpr_alm}, the dual evaluation still provides
structural priors that spare the primal search from cold-starting on $B$.
With $\lambda^{\mathrm{best}}$ the best evaluated Lagrangian multipliers
and $x^{\mathrm{best}}$ an inner minimizer of~\eqref{eq:dual_lb},
\GEAR{} starts the search at $x^{\mathrm{best}}$, which already lies in
the low-objective basin where the lower bound is attained, and initializes
$\mu$ with a nonnegative rescaling of $\lambda^{\mathrm{best}}$ so the
first outer iteration immediately emphasizes the constraints that
bottleneck the relaxation rather than accumulating multiplier weights
from zero; neither step assumes that $x^{\mathrm{best}}$ is feasible or
that $\lambda$ equals the NLP multipliers. A squared violation penalty on
the linear constraints can further guide the search by keeping early
gradient steps near the relaxed feasible region. A candidate updates
$\mathrm{UB}$ only after direct evaluation confirms feasibility
(\cref{subsec:graph_nlp}) and a better objective value, so warm starts and
penalties affect progress, not incumbent validity (\cref{app:alm-details}).

\vspace{-5pt}
\subsection{Branching and the Full Algorithm}
\label{subsec:branch}
\label{subsec:driver}
\label{subsec:algorithm_decoupling}
\label{subsec:soundness_summary}
\nobreak\vspace{-3pt}

Even with stronger dual bounds and feasibility-oriented primal search, some boxes remain unresolved because their nonlinear relaxations are still too loose.
\GEAR{} therefore applies spatial refinement, reconstructing tighter relaxations on smaller domains to progressively close the remaining gap.

\paragraph{Multiplier-weighted branching.}
Splitting a box while keeping its linear bounds fixed need not raise the minimum
across both children; the intended gain comes from reconstructing the
nonlinear relaxations on the refined domains. Smart branching (SB)~\cite{bunel2018unified,bunel2020branch}, developed for input
branch-and-bound in NN verification, therefore splits each input coordinate in turn, computes a fast lower bound on both
children, and takes the coordinate with the largest minimum child bound. This needs two bound computations per coordinate and
targets a single output. In constrained optimization a split can help by raising the objective's bound or by ruling out a
constraint, so \GEAR{} by default scores coordinates without bounding any child: it reuses the coefficients of its linear bounds 
and the Lagrangian multipliers $\lambda^{\mathrm{best}}$ to score
coordinates by weighted linear variation:
\begin{equation}
S_j=w_j\left(|a_{0,j}|+\sum_{i=1}^m
\lambda_i^{\mathrm{best}}|a_{i,j}|\right),
\qquad w_j=\overline x_j-\underline x_j.
\label{eq:multiplier_score}
\end{equation}
Width represents the spatial extent available for refinement, while the
weighted coefficients identify directions emphasized by the current bound.
Summing absolute contributions avoids cancellation between objective and
constraint terms. The score is a heuristic, not an exact prediction of
child bounds. Maximization is restricted to coordinates with
$w_j\ge\tau\max_k w_k$, for a fixed $\tau\in(0,1]$, preventing narrow
coordinates from indefinitely displacing wide ones
(\cref{app:branching}).  We call branching on the coordinate that maximizes $S_j$ \SBReweight{}.

\paragraph{The full algorithm.}
Each sBaB subproblem is a box with its own lower bound and linear bounds~\eqref{eq:method_linear_bounds}. Each iteration
selects up to $K$ open boxes with the lowest lower bounds and processes them as one batch. For each box, \GEAR{} computes new
linear bounds by bound propagation with constraint-aware intermediate tightening, tests infeasibility~\eqref{eq:infeas_cert},
and optimizes the Lagrangian bound~\eqref{eq:dual_lb}; the box is discarded if it is infeasible or if its lower bound reaches
$\mathrm{UB}$. Projected ALM on open boxes, started from the Lagrangian bound (\cref{subsec:primal}), can lower $\mathrm{UB}$
and so discard more boxes. Each remaining box is split using \SBReweight{}. Its two children inherit its lower bound and linear
bounds, which remain valid on smaller boxes, so each child can be clipped right away, before its own propagation. New linear
bounds then replace the inherited ones while the linear constraints stay fixed, so the memory per box does not grow with depth.

\paragraph{Soundness decoupled from heuristic progress.}
Correctness follows the sBaB principle of \cref{subsec:sbab}: clipping preserves feasible points,
infeasibility pruning requires a contradiction, dominance pruning requires
$\mathrm{LB}(B)\ge\mathrm{UB}$, and incumbents are checked on the original constraints, so the global
certificate is strictly decoupled from whether any heuristic step (dual ascent, ALM, or branching) succeeds. 
The solver stops at the
relative gap of \cref{eq:exp_solver_criterion}, on an empty queue, or at the time limit.  The full algorithm and invariant proof
appear in \cref{appendix::gear,app:sound-invariants}; finite termination needs the further
assumptions of \cref{app:termination}.
\section{Experiments}
\label{sec:exp}


\paragraph{Benchmarks.}
We evaluate \GEAR{} on 505 bounded continuous NLP instances: 265 from
\texttt{gams\_world} and 240 from \texttt{minlplib}.  Each retained
instance is translated into a nonlinear computation graph with continuous,
box-bounded variables.  Dataset construction, translation and filtering rules,
and validation fields are detailed in \cref{appsec:dataset}.

\paragraph{Baselines and hardware.}
We compare against six global solvers, SCIP~\cite{bestuzheva2023scip}, Gurobi~\cite{gurobi},
BARON~\cite{zhang2025baron}, Couenne~\cite{belotti2009couenne}, LindoGlobal~\cite{lin2009global}
and Xpress~\cite{xpress}, and six local solvers, Knitro~\cite{byrd2006knitro},
IPOPT~\cite{wachter2006implementation}, CONOPT~\cite{drud1994conopt}, LGO~\cite{pinter2003gams},
MINOS~\cite{murtagh2003minos} and SNOPT~\cite{gill2005snopt}, all through AMPL.  All solvers run on
one workstation (8-core 16-thread AMD Ryzen 9700X, 256\,GB RAM, one NVIDIA RTX 5090) with a
180\,s limit per instance.  Every baseline receives all 16 hardware threads in the layout it natively
supports and runs one configuration tuned on a 40-instance subset, with feasibility tolerance
$10^{-6}$ and relative gap $10^{-4}$ as in \cref{eq:exp_solver_criterion} (\cref{app:baseline}).  \GEAR{}
runs one configuration on every instance (\cref{fig:components}).

\paragraph{Solved criterion.}
A run counts as \emph{solved} when (i) its returned point satisfies every constraint and bound of
the original model, with equalities kept exact, to absolute tolerance $10^{-6}$ (baseline points
are re-evaluated in AMPL), and (ii) its primal value $P$ and certified lower bound $D$ have a
relative gap of at most $10^{-4}$:
\begin{equation}
(P - D)\,/\max(1,\lvert P\rvert,\lvert D\rvert) \;\le\; 10^{-4}.
\label{eq:exp_solver_criterion}
\end{equation}
Local solvers are credited only when their feasible $P$ matches a bound
certified by a global solver; flagged or contradicted certificates and infeasibility verdicts are
treated as in \cref{app:baseline}.

\vspace{-5pt}
\subsection{Main Results}
\label{subsec:exp_headline}
\nobreak\vspace{-3pt}

\begin{table}[t]
\centering
\footnotesize
\setlength{\tabcolsep}{2.2pt}
\renewcommand{\arraystretch}{0.92}
\resizebox{\linewidth}{!}{%
\begin{tabular}{@{}lrrrrrrr@{}}
\toprule
 & \multicolumn{2}{c}{\texttt{gams\_world}} & \multicolumn{2}{c}{\texttt{minlplib}} & \multicolumn{3}{c}{Total} \\
\cmidrule(lr){2-3}\cmidrule(lr){4-5}\cmidrule(lr){6-8}
 Global Solver & Feas. & Solved & Feas. & Solved & Feas. & Solved & \% \\
\midrule
GEAR & \textbf{264} & \textbf{251} & \textbf{237} & \underline{170} & \textbf{501} & \textbf{421} & \textbf{83.4} \\
\midrule
Gurobi & 254 & \underline{227} & 218 & \textbf{174} & 472 & \underline{401} & \underline{79.4} \\
SCIP & \underline{262} & 209 & \underline{219} & 155 & \underline{481} & 364 & 72.1 \\
BARON & 217 & 168 & 186 & 129 & 403 & 297 & 58.8 \\
Xpress & 218 & 176 & 173 & 120 & 391 & 296 & 58.6 \\
Couenne & 222 & 182 & 152 & 100 & 374 & 282 & 55.8 \\
LindoGlobal & 180 & 164 & 169 & 110 & 349 & 274 & 54.3 \\
\bottomrule
\end{tabular}%
\hspace{8pt}%
\begin{tabular}{@{}lrrrrrrr@{}}
\toprule
 & \multicolumn{2}{c}{\texttt{gams\_world}} & \multicolumn{2}{c}{\texttt{minlplib}} & \multicolumn{3}{c}{Total} \\
\cmidrule(lr){2-3}\cmidrule(lr){4-5}\cmidrule(lr){6-8}
Local Solver$^{*}$ & Feas. & Solved & Feas. & Solved & Feas. & Solved & \% \\
\midrule
GEAR & \textbf{264} & \textbf{251} & \textbf{237} & \underline{170} & \textbf{501} & \textbf{421} & \textbf{83.4} \\
\midrule
Knitro$^{*}$ & 263 & 231 & 225 & 170 & 488 & 401 & 79.4 \\
IPOPT$^{*}$ & 256 & 203 & 202 & 117 & 458 & 320 & 63.4 \\
CONOPT$^{*}$ & 259 & 205 & 205 & 91 & 464 & 296 & 58.6 \\
LGO$^{*}$ & 249 & 206 & 160 & 78 & 409 & 284 & 56.2 \\
MINOS$^{*}$ & 247 & 193 & 133 & 53 & 380 & 246 & 48.7 \\
SNOPT$^{*}$ & 235 & 182 & 133 & 51 & 368 & 233 & 46.1 \\
\bottomrule
\end{tabular}%
}
\caption{Instances solved in 180\,s (\texttt{gams\_world} 265, \texttt{minlplib} 240, total 505); left: \GEAR{} and the global solvers, right: \GEAR{} and the local solvers, as in \cref{fig:exp_cactus_morandi}. Feas.: a returned point feasible at $10^{-6}$. $^{*}$Local solvers, credited only when they match a certified bound. Bold/underline: best/second among \GEAR{} and the global solvers.}
\label{tab:exp_solver_context}
\vspace{-2mm}
\end{table}

\paragraph{Coverage.}
\Cref{tab:exp_solver_context} reports the comparison.  \GEAR{} certifies 421 of the 505
instances (83.4\%), 20 more than the strongest global baseline, Gurobi (401), and 57 more than SCIP (364).
The lead is concentrated on \texttt{gams\_world}, where \GEAR{} certifies 251/265 instances
(94.7\%) against 227 for Gurobi and 209 for SCIP; on \texttt{minlplib}, \GEAR{} (170) is ahead of SCIP (155) and
behind Gurobi (174).  Among the local solvers, Knitro matches a certified optimum on 401 instances.
\GEAR{} returns a feasible point on 501 instances.

\begin{figure}[t]
\centering
\includegraphics[width=\linewidth]{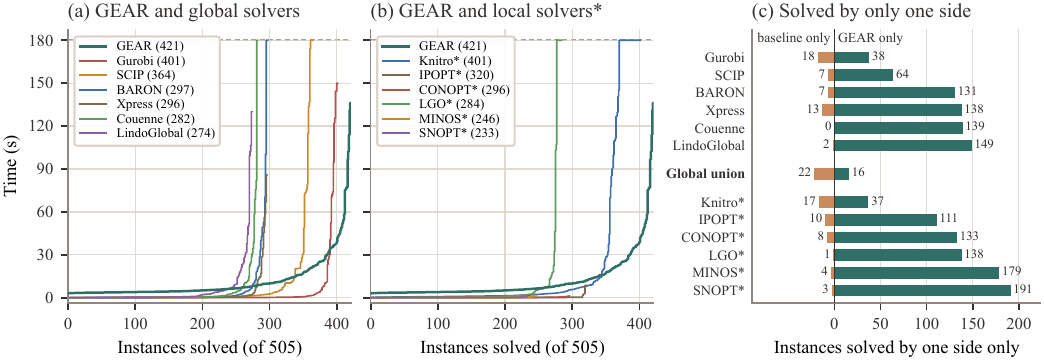}
\caption{\textbf{Solve times and exclusive solves} (180\,s, solved criterion of \cref{eq:exp_solver_criterion}).
(a, b) Cactus plots for \GEAR{} with the global and with the local solvers ($^{*}$, counted only when
matched to a globally certified bound).  (c) Instances solved by \GEAR{} but not by the baseline
(right of zero) and by the baseline but not by \GEAR{} (left of zero), on one scale; Global union: the six global solvers.}
\vspace{-5mm}
\label{fig:exp_cactus_morandi}
\end{figure}

\paragraph{Runtime profile and complementarity.}
The baselines solve easy instances in hundredths of a second, whereas \GEAR{} pays a start-up cost
of a few seconds for graph compilation and GPU initialization (\cref{fig:exp_cactus_morandi}a,b).
It overtakes LindoGlobal, Couenne, BARON and Xpress at 6-9\,s, SCIP at 15\,s and Gurobi at 35\,s,
so its lead comes from the hard tail: on instances both sides certify, the fastest global solver
needs a median of 0.03-0.06\,s against 4-29\,s for \GEAR{}.
\GEAR{} solves 38-149 instances that each global solver misses, against 0-18 in the other
direction (\cref{fig:exp_cactus_morandi}c).  Against the union of the six global solvers (427
instances), \GEAR{} adds 16 instances and misses 22; together they cover 443/505.\looseness=-1

\paragraph{Where the GPU helps.}
By problem class (\cref{appsec:gpu_benefit}), the lead comes from low- to medium-dimensional
instances with few constraints and composite nonlinearities: \GEAR{} certifies 322/330 instances with
at most 20 constraints and 131/147 with transcendental terms, which contain all 16 instances no global
solver certifies.  Factorable relaxations of composed exponential, logarithmic and power terms are
loose.  \GEAR{} closes 293 of its 421 solved instances at the root, including 12 of those 16
(\cref{appsec:win_loss}), and the other 128 in the GPU tree, with a median of 318 boxes and at least $10^5$ boxes on five of them.  On large structured models
(pooling, wastewater, water networks, circle packing), whose certificates rely on RLT and other cutting planes, the global solvers remain stronger, and \GEAR{}
solves no instance with more than 1,000 variables.

\vspace{-5pt}
\subsection{Ablation Studies}
\label{subsec:exp_ablation}
\nobreak\vspace{-3pt}

\Cref{tab:exp_component_ablation} adds one group of components per row, from the verification
pipeline to the full solver.  The first row is the most common setting in NN verifiers (\cref{subsec:naive}): input
branch-and-bound on bound propagation with a PGD attack on the box.  It solves nothing.  A verifier's attack looks for a
counterexample: its points may violate the constraints, so their objective values are not valid pruning thresholds, and it
keeps no feasible point to return.  PGD itself can reach good points: on the 95 instances without constraints, where every point
of the box is feasible, the attack reaches the optimal value on 56, but the pipeline does not keep these points as
incumbents.
Presolve and the components of \cref{appsec:components} other than the lifted relaxation supply feasible points and root
certificates: 72 instances close at the root and 147 in the branch-and-bound that follows.  The projected ALM search
(\cref{subsec:primal}) adds 50.  Constraint-aware dual bounds with clipping (\cref{subsec:dual}) are the largest step of the
GPU method (+97): once incumbents exist, the bounds must use the constraints to prune the boxes that infeasible points keep
open.  The lifted relaxation (\cref{appsec:comp_lifted}), refined operator relaxations (\cref{sec:operators}) and
\SBReweight{} branching (\cref{subsec:branch}) close 55 more.  The dual and branching steps reuse the same propagated linear
bounds, so one set of linear bounds serves every component of the solver.  With the ALM search, \GEAR{} already holds an
optimal point on 471 instances; the later rows raise the solved count from 269 to 421 and add only 7 optimal points, so they
finish the proof rather than find better points.\looseness=-1

\begin{table}[t]
\centering
\small
\setlength{\tabcolsep}{4pt}
\renewcommand{\arraystretch}{0.92}
\begin{tabular}{@{}lrrrr@{}}
\toprule
Configuration & Feasible & Optimal$^{*}$ & Solved & $\Delta$ \\
\midrule
Bound propagation + PGD & 95 & 56 & 0 & - \\
+ presolve and CPU solving components & 480 & 396 & 219 & +219 \\
+ projected ALM primal search (GPU) & 496 & 471 & 269 & +50 \\
+ constraint-aware dual bounds and clipping (GPU) & 498 & 472 & 366 & +97 \\
+ lifted and refined operator relaxation (GPU) & 501 & 478 & 410 & +44 \\
+ \SBReweight{} branching (\GEAR{}) & 501 & 478 & \textbf{421} & +11 \\
\bottomrule
\end{tabular}
\caption{Cumulative ablation of \GEAR{}.  Feasible: the returned point satisfies every constraint to $10^{-6}$.  $^{*}$Optimal: that point is also within the relative gap of the best known value, verified by any global solver.  Solved: \cref{eq:exp_solver_criterion} holds, which requires a certified lower bound.  $\Delta$: change in solved from the previous components.}
\label{tab:exp_component_ablation}
\vspace{-5mm}
\end{table}


\section{Related Works}
\label{sec:limits_related}

\textbf{Linear bound propagation beyond neural network verification.}
Linear bound propagation~\citep{zhang2018crown,singh2018fast,singh2019abstract,wang2018efficient,xu2020automatic}
bounds every output of a computation graph, and constraint-aware verifiers tighten
it with split constraints, cutting planes and domain
clipping~\citep{ijcai2026p59,wang2021betaCROWN,zhang2022general,zhou2024biccos,zhou2025clip}, also for general
nonlinearities and preimages~\citep{shi2025genbab,kotha2023invprop}.  These methods certify a property
against a fixed threshold and never need a feasible point.  In optimization the threshold is the value of
an incumbent that the solver must find and improve, and an objective-only bound can stay below it forever
(\cref{subsec:naive}).  Our view is that the linear bounds verification already computes can drive both
sides of this search: combined with Lagrangian multipliers they give a constraint-aware lower bound and an
infeasibility test, and the same bounds and multipliers start the primal search and score the branching.

\textbf{Per-node strength versus batched throughput.}
Classical global
solvers~\citep{zhang2025baron,belotti2009branching,vigerske2018scip,bestuzheva2023scip,misener2014antigone,lin2009global,gurobi}
relax a factorable program~\citep{tawarmalani2002convexification,liberti2006reformulation,khajavirad2018hybrid}
and solve a relaxation, most often an LP, at every node, strengthened by bound tightening, cutting planes and
RLT~\citep{sherali2013reformulation,gleixner2017obbt}.  An LP-based node gets the strength of an LP that keeps all
McCormick planes~\citep{McCormick1976}, but nodes are processed one after another, each with its own LP.  \GEAR{} makes the opposite trade: it keeps one linear bound per output, computed without an LP solve,
so each box is bounded more loosely, but many boxes are bounded in one batched step.  
The
experiments show the two designs to be complementary (\cref{sec:exp}): \GEAR{} leads on composite
nonlinearities in low and moderate dimension, and per-node LP strength wins on large structured models.\looseness=-1

\textbf{GPU optimization.}
GPUs already accelerate parts of both fields.  In verification, polyhedral and linear bound propagation
run on GPUs~\citep{muller2021scaling,xu2021fast} and branch-and-bound subproblems are
bounded in batches~\citep{xu2021fast,wang2021betaCROWN}.  In mathematical programming, first-order
methods such as PDLP~\citep{applegate2021practical,applegate2025pdlp} replace factorizations with
matrix-vector products to solve one large LP, and commercial solvers parallelize their tree search over
CPU threads~\citep{xpress}.  In deterministic global optimization, GPUs have offloaded McCormick
relaxations~\citep{gottlieb2024automatic} and interval bounds on partitioned subdomains~\citep{zhang2025accelerating}.
\GEAR{} takes the verification route for a global NLP search: the
branch-and-bound frontier, its linear bounds and the primal search are tensors in a deep-learning
framework~\citep{NIPS2019_9015}, so one GPU step can process thousands of boxes, instead of parallelizing
within one LP or across CPU threads.  Detailed discussion in \Cref{appsec:more_related}.\looseness=-1

\section{Conclusion}

\textsc{GEAR} demonstrates that linear bound propagation can serve as a unified GPU-batched interface for bounded continuous global NLP, reusing valid linear relaxations for lower bounding, infeasibility detection, primal recovery, and branching to solve a complementary class of highly nonlinear instances of low and moderate dimension. 
\paragraph{Future works and limitations.}\GEAR{} assumes continuous variables and a finite box; mixed-integer and
unbounded NLPs are left to future work. 
 More limitations are
discussed in Appendix~\ref{subsec:limits}.


\begin{ack}
 Huan Zhang is supported in part by the AI2050 program at Schmidt Sciences (AI2050 Early Career Fellowship) and NSF (2543005, 2331967). Grani A.~Hanasusanto is supported in part by NSF (CCF-2343869 and ECCS-2404413). We thank Yifei Li from UIUC for running part of the AMPL baselines.
 We thank the anonymous reviewers for their constructive feedback. 

\end{ack}

\bibliography{reference}

@article{zhang2025baron,
  title        = {Solving continuous and discrete nonlinear programs with {BARON}},
  author       = {Zhang, Yi and Sahinidis, Nikolaos V.},
  journal      = {Computational Optimization and Applications},
  year         = {2025},
  volume       = {92},
  pages        = {1123--1161},
  doi          = {10.1007/s10589-024-00633-0}
}

@article{zhou2025clip,
  title={Clip-and-verify: Linear constraint-driven domain clipping for accelerating neural network verification},
  author={Zhou, Duo and Chavez, Jorge and Chen, Hesun and Hanasusanto, Grani A and Zhang, Huan},
  journal={Advances in Neural Information Processing Systems},
  volume={38},
  pages={174849--174895},
  year={2026}
}

@article{misener2014antigone,
  title={{ANTIGONE}: Algorithms for continuous/integer global optimization of nonlinear equations},
  author={Misener, Ruth and Floudas, Christodoulos A},
  journal={Journal of Global Optimization},
  volume={59},
  number={2},
  pages={503--526},
  year={2014},
  publisher={Springer}
}

@article{doig1960automaticmethodofsolvingdiscrete,
 ISSN = {00129682, 14680262},
 URL = {http://www.jstor.org/stable/1910129},
 author = {A. H. Land and A. G. Doig},
 journal = {Econometrica},
 number = {3},
 pages = {497--520},
 publisher = {[Wiley, Econometric Society]},
 title = {An Automatic Method of Solving Discrete Programming Problems},
 urldate = {2025-05-13},
 volume = {28},
 year = {1960}
}

@article{lu2019neural,
  title={Neural network branching for neural network verification},
  author={Lu, Jingyue and Kumar, M Pawan},
  journal={arXiv preprint arXiv:1912.01329},
  year={2019}
}

@article{de2021improved,
  title={Improved branch and bound for neural network verification via {Lagrangian} decomposition},
  author={De Palma, Alessandro and Bunel, Rudy and Desmaison, Alban and Dvijotham, Krishnamurthy and Kohli, Pushmeet and Torr, Philip HS and Kumar, M Pawan},
  journal={arXiv preprint arXiv:2104.06718},
  year={2021}
}

@inproceedings{lan2022tight,
  title={Tight neural network verification via semidefinite relaxations and linear reformulations},
  author={Lan, Jianglin and Zheng, Yang and Lomuscio, Alessio},
  booktitle={Proceedings of the AAAI Conference on Artificial Intelligence},
  volume={36},
  pages={7272--7280},
  year={2022}
}

@inproceedings{sun2022romax,
  title={{ROMAX}: Certifiably robust deep multiagent reinforcement learning via convex relaxation},
  author={Sun, Chuangchuang and Kim, Dong-Ki and How, Jonathan P},
  booktitle={2022 International Conference on Robotics and Automation (ICRA)},
  pages={5503--5510},
  year={2022},
  organization={IEEE}
}

@inproceedings{shi2025genbab,
  title={Neural Network Verification with Branch-and-Bound for General Nonlinearities},
  author={Shi, Zhouxing and Jin, Qirui and Kolter, Zico and Jana, Suman and Hsieh, Cho-Jui and Zhang, Huan},
  booktitle={International Conference on Tools and Algorithms for the Construction and Analysis of Systems},
  year={2025}
}

@article{applegate2021practical,
  title={Practical large-scale linear programming using primal-dual hybrid gradient},
  author={Applegate, David and D{\'\i}az, Mateo and Hinder, Oliver and Lu, Haihao and Lubin, Miles and O'Donoghue, Brendan and Schudy, Warren},
  journal={Advances in Neural Information Processing Systems},
  volume={34},
  pages={20243--20257},
  year={2021}
}

@inproceedings{ijcai2026p59,
  title     = {Lookahead Branching for Neural Network Verification},
  author    = {Davis, Liam and Zhou, Duo and Zhang, Huan and Katz, Guy and Barrett, Clark and Wu, Haoze},
  booktitle = {Proceedings of the Thirty-Fifth International Joint Conference on
               Artificial Intelligence, {IJCAI-26}},
  publisher = {International Joint Conferences on Artificial Intelligence Organization},
  editor    = {Diego Calvanese},
  pages     = {518--526},
  year      = {2026},
  month     = {8},
  note      = {Main Track},
  doi       = {10.24963/ijcai.2026/59},
  url       = {https://doi.org/10.24963/ijcai.2026/59},
}

@article{applegate2025pdlp,
  title={{PDLP}: A Practical First-Order Method for Large-Scale Linear Programming},
  author={Applegate, David and D{\'\i}az, Mateo and Hinder, Oliver and Lu, Haihao and Lubin, Miles and O'Donoghue, Brendan and Schudy, Warren},
  journal={arXiv preprint arXiv:2501.07018},
  year={2025}
}

@article{venzke2020verification,
  title={Verification of neural network behaviour: Formal guarantees for power system applications},
  author={Venzke, Andreas and Chatzivasileiadis, Spyros},
  journal={IEEE Transactions on Smart Grid},
  volume={12},
  number={1},
  pages={383--397},
  year={2020},
  publisher={IEEE}
}

@inproceedings{wong2020neural,
  title={Neural network virtual sensors for fuel injection quantities with provable performance specifications},
  author={Wong, Eric and Schneider, Tim and Schmitt, Joerg and Schmidt, Frank R and Kolter, J Zico},
  booktitle={2020 IEEE Intelligent Vehicles Symposium (IV)},
  pages={1753--1758},
  year={2020},
  organization={IEEE}
}

@article{wu2023neural,
  title={Neural {Lyapunov} control for discrete-time systems},
  author={Wu, Junlin and Clark, Andrew and Kantaros, Yiannis and Vorobeychik, Yevgeniy},
  journal={Advances in Neural Information Processing Systems},
  volume={36},
  pages={2939--2955},
  year={2023}
}

@article{chen2024verification,
  title={Verification-Aided Learning of Neural Network Barrier Functions with Termination Guarantees},
  author={Chen, Shaoru and Molu, Lekan and Fazlyab, Mahyar},
  journal={arXiv preprint arXiv:2403.07308},
  year={2024}
}

@inproceedings{yang2024lyapunov,
author = {Yang, Lujie and Dai, Hongkai and Shi, Zhouxing and Hsieh, Cho-Jui and Tedrake, Russ and Zhang, Huan},
title = {Lyapunov-stable neural control for state and output feedback: a novel formulation},
year = {2024},
publisher = {JMLR.org},
booktitle = {Proceedings of the 41st International Conference on Machine Learning},
articleno = {2311},
numpages = {14},
location = {Vienna, Austria},
series = {ICML'24}
}

@article{zhang2022general,
  title={General cutting planes for bound-propagation-based neural network verification},
  author={Zhang, Huan and Wang, Shiqi and Xu, Kaidi and Li, Linyi and Li, Bo and Jana, Suman and Hsieh, Cho-Jui and Kolter, J Zico},
  journal={Advances in Neural Information Processing Systems},
  volume={35},
  pages={1656--1670},
  year={2022}
}

@inproceedings{katz2017reluplex,
  title={Reluplex: An efficient {SMT} solver for verifying deep neural networks},
  author={Katz, Guy and Barrett, Clark and Dill, David L and Julian, Kyle and Kochenderfer, Mykel J},
  booktitle={Computer Aided Verification: 29th International Conference, CAV 2017, Heidelberg, Germany, July 24-28, 2017, Proceedings, Part I 30},
  pages={97--117},
  year={2017},
  organization={Springer}
}

@article{tjeng2017evaluating,
  title={Evaluating robustness of neural networks with mixed integer programming},
  author={Tjeng, Vincent and Xiao, Kai and Tedrake, Russ},
  journal={arXiv preprint arXiv:1711.07356},
  year={2017}
}

@article{bunel2018unified,
  title={A unified view of piecewise linear neural network verification},
  author={Bunel, Rudy R and Turkaslan, Ilker and Torr, Philip and Kohli, Pushmeet and Kumar, M. Pawan},
  journal={Advances in Neural Information Processing Systems},
  volume={31},
  year={2018}
}

@article{bunel2020branch,
  title={Branch and bound for piecewise linear neural network verification},
  author={Bunel, Rudy and Lu, Jingyue and Turkaslan, Ilker and Torr, Philip HS and Kohli, Pushmeet and Kumar, M Pawan},
  journal={Journal of Machine Learning Research},
  volume={21},
  number={42},
  pages={1--39},
  year={2020}
}

@article{anderson2020strong,
  title={Strong mixed-integer programming formulations for trained neural networks},
  author={Anderson, Ross and Huchette, Joey and Ma, Will and Tjandraatmadja, Christian and Vielma, Juan Pablo},
  journal={Mathematical Programming},
  volume={183},
  number={1-2},
  pages={3--39},
  year={2020},
  publisher={Springer}
}

@article{xu2020automatic,
  title={Automatic perturbation analysis for scalable certified robustness and beyond},
  author={Xu, Kaidi and Shi, Zhouxing and Zhang, Huan and Wang, Yihan and Chang, Kai-Wei and Huang, Minlie and Kailkhura, Bhavya and Lin, Xue and Hsieh, Cho-Jui},
  journal={Advances in Neural Information Processing Systems},
  volume={33},
  pages={1129--1141},
  year={2020}
}

@article{morrison2016babsurvery,
title = {Branch-and-bound algorithms: A survey of recent advances in searching, branching, and pruning},
journal = {Discrete Optimization},
volume = {19},
pages = {79-102},
year = {2016},
issn = {1572-5286},
doi = {https://doi.org/10.1016/j.disopt.2016.01.005},
url = {https://www.sciencedirect.com/science/article/pii/S1572528616000062},
author = {David R. Morrison and Sheldon H. Jacobson and Jason J. Sauppe and Edward C. Sewell}
}

@article{ferrari2022complete,
  title={Complete Verification via Multi-Neuron Relaxation Guided Branch-and-Bound},
  author={Ferrari, Claudio and M{\"u}ller, Mark Niklas and Jovanovi{\'c}, Nikola and Vechev, Martin},
  journal={arXiv preprint arXiv:2205.00263},
  year={2022}
}

@book{sherali2013reformulation,
  title={A reformulation-linearization technique for solving discrete and continuous nonconvex problems},
  author={Sherali, Hanif D and Adams, Warren P},
  volume={31},
  year={2013},
  publisher={Springer Science \& Business Media}
}

@inproceedings{bak2021nfm,
  title={nnenum: Verification of {ReLU} Neural Networks with Optimized Abstraction Refinement},
  author={Bak, Stanley},
  booktitle={NASA Formal Methods Symposium},
  pages={19--36},
  year={2021},
  organization={Springer}
}

@article{muller2021scaling,
  title={Scaling polyhedral neural network verification on {GPUs}},
  author={M{\"u}ller, Christoph and Serre, Francois and Singh, Gagandeep and P{\"u}schel, Markus and Vechev, Martin},
  journal={Proceedings of Machine Learning and Systems},
  volume={3},
  pages={733--746},
  year={2021}
}

@article{dathathri2020enabling,
  title={Enabling certification of verification-agnostic networks via memory-efficient semidefinite programming},
  author={Dathathri, Sumanth and Dvijotham, Krishnamurthy and Kurakin, Alexey and Raghunathan, Aditi and Uesato, Jonathan and Bunel, Rudy R and Shankar, Shreya and Steinhardt, Jacob and Goodfellow, Ian and Liang, Percy S and Kohli, Pushmeet},
  journal={Advances in Neural Information Processing Systems (NeurIPS)},
  year={2020}
}

@misc{brown2022unifiedviewsdpbasedneural,
      title={A Unified View of {SDP}-based Neural Network Verification through Completely Positive Programming}, 
      author={Robin Brown and Edward Schmerling and Navid Azizan and Marco Pavone},
      year={2022},
      eprint={2203.03034},
      archivePrefix={arXiv},
      primaryClass={math.OC},
      url={https://arxiv.org/abs/2203.03034}, 
}

@article{fischetti2017deepneuralnetworksasmilp,
  author       = {Matteo Fischetti and
                  Jason Jo},
  title        = {Deep Neural Networks as 0-1 Mixed Integer Linear Programs: {A} Feasibility
                  Study},
  journal      = {CoRR},
  volume       = {abs/1712.06174},
  year         = {2017},
  url          = {http://arxiv.org/abs/1712.06174},
  eprinttype    = {arXiv},
  eprint       = {1712.06174},
  bibsource    = {dblp computer science bibliography, https://dblp.org}
}

@article{fazlyab2020safety,
  title={Safety verification and robustness analysis of neural networks via quadratic constraints and semidefinite programming},
  author={Fazlyab, Mahyar and Morari, Manfred and Pappas, George J},
  journal={IEEE Transactions on Automatic Control},
  year={2020},
  publisher={IEEE}
}

@article{singh2019abstract,
  title={An abstract domain for certifying neural networks},
  author={Singh, Gagandeep and Gehr, Timon and P{\"u}schel, Markus and Vechev, Martin},
  journal={Proceedings of the ACM on Programming Languages},
  volume={3},
  number={POPL},
  pages={1--30},
  year={2019},
  publisher={ACM New York, NY, USA}
}

@inproceedings{singh2018fast,
  title={Fast and Effective Robustness Certification},
  author={Singh, Gagandeep and Gehr, Timon and Mirman, Matthew and P{\"u}schel, Markus and Vechev, Martin},
  booktitle={Advances in Neural Information Processing Systems (NeurIPS)},
  year={2018}
}

@inproceedings{wang2018efficient,
  title={Efficient formal safety analysis of neural networks},
  author={Wang, Shiqi and Pei, Kexin and Whitehouse, Justin and Yang, Junfeng and Jana, Suman},
  booktitle={Advances in Neural Information Processing Systems (NeurIPS)},
  year={2018}
}

@inproceedings{wang2018formal,
  title={Formal security analysis of neural networks using symbolic intervals},
  author={Wang, Shiqi and Pei, Kexin and Whitehouse, Justin and Yang, Junfeng and Jana, Suman},
  booktitle={USENIX Security Symposium},
  year={2018}
}

@inproceedings{raghunathan2018semidefinite,
  title={Semidefinite relaxations for certifying robustness to adversarial examples},
  author={Raghunathan, Aditi and Steinhardt, Jacob and Liang, Percy S},
  booktitle={Advances in Neural Information Processing Systems (NeurIPS)},
  year={2018}
}

@incollection{NIPS2019_9015,
title = {{PyTorch}: An Imperative Style, High-Performance Deep Learning Library},
author = {Paszke, Adam and Gross, Sam and Massa, Francisco and Lerer, Adam and Bradbury, James and Chanan, Gregory and Killeen, Trevor and Lin, Zeming and Gimelshein, Natalia and Antiga, Luca and Desmaison, Alban and K{\"o}pf, Andreas and Yang, Edward and DeVito, Zachary and Raison, Martin and Tejani, Alykhan and Chilamkurthy, Sasank and Steiner, Benoit and Fang, Lu and Bai, Junjie and Chintala, Soumith},
booktitle = {Advances in Neural Information Processing Systems (NeurIPS)},
year = {2019}
}

@article{kingma2014adam,
  title={Adam: A method for stochastic optimization},
  author={Kingma, Diederik P and Ba, Jimmy},
  journal={International Conference on Learning Representations (ICLR)},
  year={2015}
}

@article{huchette2023deep,
  title={When deep learning meets polyhedral theory: A survey},
  author={Huchette, Joey and Mu{\~n}oz, Gonzalo and Serra, Thiago and Tsay, Calvin},
  journal={arXiv preprint arXiv:2305.00241},
  year={2023}
}

@article{zhou2024biccos,
  title={Scalable neural network verification with branch-and-bound inferred cutting planes},
  author={Zhou, Duo and Brix, Christopher and Hanasusanto, Grani A and Zhang, Huan},
  journal={Advances in Neural Information Processing Systems},
  volume={37},
  pages={29324--29353},
  year={2024}
}

@inproceedings{chiu2025sdpcrown,
  title        = {{SDP-CROWN}: Efficient Bound Propagation for Neural Network Verification with Tightness of Semidefinite Programming},
  author       = {Chiu, Hong-Ming and Chen, Hao and Zhang, Huan and Zhang, Richard Y.},
  booktitle    = {International Conference on Machine Learning (ICML)},
  year         = {2025},
  note         = {ICML 2025}
}

@book{tawarmalani2002convexification,
  title = {Convexification and Global Optimization in Continuous and Mixed-Integer Nonlinear Programming: Theory, Algorithms, Software, and Applications},
  author = {Tawarmalani, Mohit and Sahinidis, Nikolaos V.},
  year = {2002},
  publisher = {Springer},
  series = {Nonconvex Optimization and Its Applications}
}

@article{khajavirad2018hybrid,
  title = {A Hybrid {LP/NLP} Paradigm for Global Optimization Relaxations},
  author = {Khajavirad, Aida and Sahinidis, Nikolaos V.},
  journal = {Mathematical Programming Computation},
  volume = {10},
  number = {3},
  pages = {383--421},
  year = {2018}
}

@article{liberti2006reformulation,
  title = {An Exact Reformulation Algorithm for Large Nonconvex {NLPs} Involving Bilinear Terms},
  author = {Liberti, Leo and Pantelides, Constantinos C.},
  journal = {Journal of Global Optimization},
  volume = {36},
  number = {2},
  pages = {161--189},
  year = {2006}
}

@article{belotti2009branching,
  title={Branching and bounds tightening techniques for non-convex {MINLP}},
  author={Belotti, Pietro and Lee, Jon and Liberti, Leo and Margot, Francois and W{\"a}chter, Andreas},
  journal={Optimization Methods \& Software},
  volume={24},
  number={4-5},
  pages={597--634},
  year={2009},
  publisher={Taylor \& Francis}
}

@article{xpress,
  title={Parallelization of the {FICO Xpress-Optimizer}},
  author={Berthold, Timo and Farmer, James and Heinz, Stefan and Perregaard, Michael},
  journal={Optimization Methods and Software},
  volume={33},
  number={3},
  pages={518--529},
  year={2018},
  publisher={Taylor \& Francis}
}

@misc{gurobi,
  author = {{Gurobi Optimization, LLC}},
  title  = {{Gurobi Optimizer Reference Manual}},
  year   = 2026,
  url    = {https://www.gurobi.com}
}

@article{rockafellar1973dual,
  title={A dual approach to solving nonlinear programming problems by unconstrained optimization},
  author={Rockafellar, R Tyrrell},
  journal={Mathematical programming},
  volume={5},
  number={1},
  pages={354--373},
  year={1973},
  publisher={Springer}
}

@article{McCormick1976,
  author  = {G. P. McCormick},
  title   = {Computability of global solutions to factorable nonconvex programs: Part {I} --- Convex underestimating problems},
  journal = {Mathematical Programming},
  volume  = {10},
  number  = {1},
  pages   = {147--175},
  year    = {1976}
}

@misc{belotti2009couenne,
  title={Couenne: a user’s manual},
  author={Belotti, Pietro},
  year={2009},
  publisher={Technical report, Lehigh University}
}

@misc{pinter2003gams,
  title={{GAMS/LGO} nonlinear solver suite: key features, usage, and numerical performance},
  author={Pint{\'e}r, J{\'a}nos D.},
  year={2003}
}

@article{lin2009global,
  title={The global solver in the {LINDO API}},
  author={Lin, Youdong and Schrage, Linus},
  journal={Optimization Methods \& Software},
  volume={24},
  number={4-5},
  pages={657--668},
  year={2009},
  publisher={Taylor \& Francis}
}

@incollection{byrd2006knitro,
  title={KNITRO: An integrated package for nonlinear optimization},
  author={Byrd, Richard H and Nocedal, Jorge and Waltz, Richard A},
  booktitle={Large-scale nonlinear optimization},
  pages={35--59},
  year={2006},
  publisher={Springer}
}

@book{biegler2010nonlinear,
  title={Nonlinear programming: concepts, algorithms, and applications to chemical processes},
  author={Biegler, Lorenz T},
  year={2010},
  publisher={SIAM}
}

@article{rockafellar1973augmented,
  title = {The Multiplier Method of {Hestenes} and {Powell} Applied to Convex Programming},
  author = {Rockafellar, R. Tyrrell},
  journal = {Journal of Optimization Theory and Applications},
  volume = {12},
  number = {6},
  pages = {555--562},
  year = {1973}
}

@techreport{vigerske2018scip,
  title = {{SCIP}: Global Optimization of Mixed-Integer Nonlinear Programs in a Branch-and-Cut Framework},
  author = {Vigerske, Stefan and Gleixner, Ambros},
  institution = {Konrad-Zuse-Zentrum f{\"u}r Informationstechnik Berlin},
  year = {2018}
}

@misc{bestuzheva2023scip,
      title={Global Optimization of Mixed-Integer Nonlinear Programs with {SCIP} 8}, 
      author={Ksenia Bestuzheva and Antonia Chmiela and Benjamin Müller and Felipe Serrano and Stefan Vigerske and Fabian Wegscheider},
      year={2023},
      eprint={2301.00587},
      archivePrefix={arXiv},
      primaryClass={math.OC},
      url={https://arxiv.org/abs/2301.00587}, 
}

@inproceedings{wang2021betaCROWN,
  title = {Beta-{CROWN}: Efficient Bound Propagation with Per-neuron Split Constraints for Neural Network Robustness Verification},
  author = {Wang, Shiqi and Zhang, Huan and Xu, Kaidi and Lin, Xue and Jana, Suman and Hsieh, Cho-Jui and Kolter, J.~Zico},
  booktitle = {Advances in Neural Information Processing Systems (NeurIPS)},
  year = {2021}
}

@inproceedings{xu2021fast,
  title = {Fast and Complete: Enabling Complete Neural Network Verification with Rapid and Massively Parallel Incomplete Verifiers},
  author = {Xu, Kaidi and Zhang, Huan and Wang, Shiqi and Wang, Yihan and Jana, Suman and Lin, Xue and Hsieh, Cho-Jui},
  booktitle = {International Conference on Learning Representations (ICLR)},
  year = {2021}
}

@inproceedings{zhang2018crown,
  title = {Efficient Neural Network Robustness Certification with General Activation Functions},
  author = {Zhang, Huan and Weng, Tsui-Wei and Chen, Pin-Yu and Hsieh, Cho-Jui and Daniel, Luca},
  booktitle = {Advances in Neural Information Processing Systems (NeurIPS)},
  year = {2018}
}

@inproceedings{kotha2023invprop,
  author       = {Suhas Kotha and
                  Christopher Brix and
                  J. Zico Kolter and
                  Krishnamurthy Dvijotham and
                  Huan Zhang},
  title        = {Provably Bounding Neural Network Preimages},
  booktitle    = {Advances in Neural Information Processing Systems (NeurIPS)},
  year         = {2023}
}

@article{gleixner2017obbt,
  author       = {Ambros M. Gleixner and
                  Timo Berthold and
                  Benjamin M{\"{u}}ller and
                  Stefan Weltge},
  title        = {Three enhancements for optimization-based bound tightening},
  journal      = {Journal of Global Optimization},
  volume       = {67},
  number       = {4},
  pages        = {731--757},
  year         = {2017},
  doi          = {10.1007/s10898-016-0450-4}
}

@article{wachter2006implementation,
  title={On the implementation of an interior-point filter line-search algorithm for large-scale nonlinear programming},
  author={W{\"a}chter, Andreas and Biegler, Lorenz T},
  journal={Mathematical programming},
  volume={106},
  number={1},
  pages={25--57},
  year={2006},
  publisher={Springer}
}

@article{gill2005snopt,
  title={{SNOPT}: An {SQP} algorithm for large-scale constrained optimization},
  author={Gill, Philip E and Murray, Walter and Saunders, Michael A},
  journal={SIAM review},
  volume={47},
  number={1},
  pages={99--131},
  year={2005},
  publisher={SIAM}
}

@techreport{murtagh2003minos,
  author      = {Bruce A. Murtagh and Michael A. Saunders},
  title       = {{MINOS} 5.51 {U}ser's {G}uide},
  institution = {Systems Optimization Laboratory, Stanford University},
  number      = {SOL 83-20R},
  year        = {2003},
  note        = {Originally published 1983; revised September 23, 2003}
}

@article{drud1994conopt,
  title   = {{CONOPT}---A Large-Scale {GRG} Code},
  author  = {Drud, Arne Stolbjerg},
  journal = {ORSA Journal on Computing},
  volume  = {6},
  number  = {2},
  pages   = {207--216},
  year    = {1994}
}

@article{adjiman1998global,
  title   = {A global optimization method, {$\alpha$BB}, for general twice-differentiable constrained {NLPs}: {I}. {T}heoretical advances},
  author  = {Adjiman, C. S. and Dallwig, S. and Floudas, C. A. and Neumaier, A.},
  journal = {Computers \& Chemical Engineering},
  volume  = {22},
  number  = {9},
  pages   = {1137--1158},
  year    = {1998}
}

@article{burer2008finite,
  title   = {A finite branch-and-bound algorithm for nonconvex quadratic programming via semidefinite relaxations},
  author  = {Burer, Samuel and Vandenbussche, Dieter},
  journal = {Mathematical Programming},
  volume  = {113},
  number  = {2},
  pages   = {259--282},
  year    = {2008}
}

@article{nohra2021spectral,
  title   = {Spectral relaxations and branching strategies for global optimization of mixed-integer quadratic programs},
  author  = {Nohra, Carlos J. and Raghunathan, Arvind U. and Sahinidis, Nikolaos V.},
  journal = {SIAM Journal on Optimization},
  volume  = {31},
  number  = {1},
  pages   = {142--171},
  year    = {2021}
}

@article{gottlieb2024automatic,
  title   = {Automatic source code generation for deterministic global optimization with parallel architectures},
  author  = {Gottlieb, Robert X. and Xu, Pengfei and Stuber, Matthew D.},
  journal = {Optimization Methods and Software},
  volume  = {41},
  number  = {2},
  pages   = {308--346},
  year    = {2026}
}

@article{zhang2025accelerating,
  title   = {Accelerating deterministic global optimization via {GPU}-parallel interval arithmetic},
  author  = {Zhang, Hongzhen and Kerkenhoff, Tim and Kichler, Neil and Dahmen, Manuel and Mitsos, Alexander and Naumann, Uwe and Bongartz, Dominik},
  journal = {arXiv preprint arXiv:2507.20769},
  year    = {2025}
}
\bibliographystyle{plain}


\raggedbottom
\normalsize
\makeatletter\@secpenalty=-300\makeatother
\let\paragraph\GEARparagraph
\newpage
\appendix
\crefalias{section}{appendix}\crefalias{subsection}{subappendix}
\onecolumn
\makeatletter
\@ifundefined{proposition}{\newtheorem{proposition}[theorem]{Proposition}}{}
\makeatother

\addtocontents{toc}{\protect\setcounter{tocdepth}{2}}
\renewcommand{\contentsname}{Appendix Contents}
\begingroup
\setlength{\parskip}{0pt}
\tableofcontents
\endgroup

\noindent\Cref{app:methodology} gives the derivations and proofs for the GPU method of \cref{sec:method}.
\Cref{appsec:components} describes every component of the solver evaluated in \cref{sec:exp}.
\Cref{appsec:dataset} gives the benchmark, the baselines and how runs are scored.
\Cref{appsec:analysis} analyzes the results and states the limitations.
\Cref{appsec:more_related} covers further related work.
\newpage

\section{Methodology Details and Proofs}
\label{appsec:proofs}
\label{app:methodology}

\subsection{Roadmap and Soundness Requirements}
\label{app:roadmap}

This section supplies the derivations and implementation details omitted
from \cref{sec:method}. It uses the same notation as the main text:
$B=[\underline x,\overline x]$, $\underline f_i(x)=a_i^\top x+c_i$,
$\mathcal P(B)=\{x\in B:Ax+c\le0\}$, and the box-domain LP dual
$d_B(\lambda)$ of \eqref{eq:dual_lb}. The vector $\lambda$ belongs to the
propagated-row LP, whereas $\mu$ belongs to the original-NLP primal
heuristic. A bound on $\mathcal F(B)$ suffices for all pruning arguments;
it need not hold at infeasible inputs excluded during intermediate
clipping.

\begin{table}[h]
\centering
\small
\setlength{\tabcolsep}{4pt}
\begin{tabular}{p{0.29\linewidth}p{0.64\linewidth}}
\toprule
\textbf{Main-text component} & \textbf{Appendix expansion} \\
\midrule
Section~\ref{subsec:representation}: representation &
\Cref{app:batched-linear-bound-state,app:admissibility}: graph, tensor layout, and validity under clipping. \\
Section~\ref{subsec:dual}: objective bound &
\Cref{app:dual-proofs,app:sound-invariants}: LP duality, supergradients, best-iterate updates, and infeasibility certificates. \\
Section~\ref{subsec:dual}: clipping &
\Cref{appendix::clipping,app:concretization}: coordinate updates, intermediate-bound duals, and conservative concretization. \\
Section~\ref{subsec:dual}: operator relaxations &
\Cref{app:operators}: plane selection, relaxation formulas, validity checks, and examples. \\
Section~\ref{subsec:primal}: primal search &
\Cref{app:alm-details}: projected ALM, warm starts, linear-row penalties, and incumbent acceptance. \\
Section~\ref{subsec:branch}: integration &
\Cref{app:node-lifecycle,app:branching,appendix::gear}: data dependencies, branching, and the full algorithm; \cref{app:termination}: conditional termination. \\
\bottomrule
\end{tabular}
\caption{Location of the details supporting each main-text mechanism.}
\label{tab:appendix_to_method}
\end{table}

\paragraph{Certificates versus heuristics.}
The solver may terminate dual and primal optimization after finite
budgets. Pruning nevertheless requires a valid certificate: a lower bound,
a contradiction, or a comparison with a feasible incumbent. The relevant
requirements are summarized in \cref{tab:sound-vs-heuristic}. In this
appendix a \emph{linear feasibility cut} means the same inequality
$\underline f_i(x)\le0$ as a \emph{linear constraint relaxation}; it does
not imply a separately generated cutting-plane family.

\begin{table}[h]
\centering
\small
\setlength{\tabcolsep}{4pt}
\begin{tabular}{p{0.22\linewidth}p{0.35\linewidth}p{0.34\linewidth}}
\toprule
\textbf{Object} & \textbf{Required for validity} & \textbf{Efficiency choice} \\
\midrule
Affine row & $\underline f_i\le f_i$ on the node's feasible set & Relaxation selection and batching \\
Constraint relaxation & $\mathcal F(B)\subseteq\mathcal P(B)$ & Inherited-row replacement schedule \\
Objective bound & Conservative evaluation of $d_B(\lambda)$ with $\lambda\ge0$ & Initialization, optimizer, and iteration budget \\
Infeasibility test & A strictly positive conservative certificate & Individual rows or nonnegative mixtures \\
Clipping & Every feasible input and its intermediate values are retained & Update order and number of sweeps \\
Incumbent & Original-constraint feasibility and direct objective evaluation & ALM schedule, seeds, and warm starts \\
Branching & Children cover the surviving parent domain & Coefficient score and eligible coordinates \\
\bottomrule
\end{tabular}
\caption{Soundness requirements and heuristic choices. The numerical
requirements are conditional on conservative evaluations as detailed in
\cref{app:admissibility}; heuristic choices do not replace them.}
\label{tab:sound-vs-heuristic}
\end{table}

\subsection{Batched Linear Bounds from a Multi-Output Graph}
\label{app:batched-linear-bound-state}
\label{app:gear_representation}

The objective and all constraints are the outputs of the graph $F$ in~\eqref{eq:nlp_graph}.
A constraint $f_i(x)=a_i^\top x+c_i$ that is linear in $x$ is reproduced exactly by propagation, and when
the parser exposes such constraints as an explicit matrix, that matrix is also retained directly. The residual $\max_{i\ge1}f_i(x)$ can be
bounded from constraint rows without appending a hard-max graph operator.

\paragraph{Backward linear bound propagation.}
For a scalar output selector $s$, a direct bound-propagation pass seeks
\begin{equation}
L_s(x)=a_s^\top x+c_s\le s^\top F(x)
\label{eq:selector-linear-bound}
\end{equation}
on its validity set. If $z=h(y)$ has valid linear relaxations on the
intermediate interval, backward linear bound propagation uses
\begin{equation}
\gamma h(y)\ge\gamma_+\underline h(y)+\gamma_-\overline h(y),
\qquad \gamma_+=\max\{\gamma,0\},\quad\gamma_-=\min\{\gamma,0\}.
\label{eq:sign-compatible-composition}
\end{equation}
Affine operators are substituted exactly. Applying this inequality in
reverse topological order proves output-row validity. With ordinary
intervals it holds on $B$; with constraint-tightened intervals it holds
on $\mathcal F(B)$, whose graph trajectories remain within those
intervals. Combining cached lower rows with nonnegative weights is also
valid, but need not give the same row as rerunning bound propagation on the
combined output: relaxation selection and coefficient cancellation can
differ. The dual uses cached rows and does not rerun propagation.

\paragraph{Bound records and tensor layout.}
A mathematical record for an output row is
\begin{equation}
\mathfrak b_i(B)=(i,a_i,c_i,V_i),\qquad
\mathcal F(B)\subseteq V_i,
\quad a_i^\top x+c_i\le f_i(x)\ \text{on }V_i.
\label{eq:bound-record}
\end{equation}
Relaxation parameters and numerical error budgets can accompany this record;
the proof only uses its validity statement. For $K$ boxes, the logical
coefficient arrays are
\begin{equation}
\mathsf A\in\mathbb R^{K\times(m+1)\times n},\qquad
\mathsf c\in\mathbb R^{K\times(m+1)},
\label{eq:batch_tensors}
\end{equation}
where $\mathsf A[k,i,:]=a_i^{(k)}$ and
$\mathsf c[k,i]=c_i^{(k)}$. Exact affine rows are invariant across boxes
and may be broadcast rather than duplicated. Each dual update costs
$O(mn)$ arithmetic per box. Primal steps still require graph evaluations,
and changed intermediate intervals can require new propagation; only
operations on the fixed cached rows avoid a graph traversal.

\paragraph{Residual bounds without a max relaxation.}
Let $r(x)=\max_{1\le i\le m}f_i(x)$, with the unconstrained case handled
without a residual test. Each constraint lower row is a lower bound on
$r$. More generally, for $\omega\ge0$ and $\mathbf1^\top\omega=1$,
\[
\sum_i\omega_i\underline f_i(x)
\le\sum_i\omega_i f_i(x)\le r(x)
\qquad(x\in\mathcal F(B)).
\]
If the box minimum of this row is positive, no feasible point exists.
Normalization gives a residual interpretation; it is not necessary for
the nonnegative-mixture contradiction in \eqref{eq:infeas_cert}.
The selected mixture family and its search budget determine which
infeasible nodes are detected, not whether a successful certificate is
valid.

\subsection{Admissible Linear Constraint Relaxations and Tolerance Convention}
\label{app:admissibility}

Admissibility means preserving the feasible set, including when rows have
been derived under earlier constraints. If a row was obtained at a box
$B_k$, the required statement for use at a descendant $B\subseteq B_k$ is
\begin{equation}
a_k^\top x+c_k\le f_{i_k}(x)
\quad\forall x\in\mathcal F(B_k).
\label{eq:app_linear_validity}
\end{equation}
A record of its constraint index and provenance may be written
\begin{equation}
\mathfrak c_k=(i_k,a_k,c_k,B_k).
\label{eq:cut-record}
\end{equation}
Since $\mathcal F(B)\subseteq\mathcal F(B_k)$, every feasible descendant
point satisfies $a_k^\top x+c_k\le0$. A row need not be valid throughout
$B_k$ if its construction used feasible-only intermediate bounds.
The constraint-tightening argument is inductive: each step starts with
valid rows, retains all feasible trajectories, and derives new valid rows
from the retained intervals.

\paragraph{Storage and forwarding.}
The initial active set contains the original affine constraints. At a
child it additionally contains the nonlinear rows inherited from its
parent. After fresh propagation these nonlinear rows are replaced, not
accumulated with all ancestors. Thus the refreshed set does not grow with depth. Replacing rows can lose relaxation strength, but not
feasibility preservation. Scalar lower bounds are separately retained
by taking the maximum with the inherited bound.

\paragraph{Feasibility tolerances.}
The functions $f_i$ throughout the method already include the tolerances
of \cref{subsec:graph_nlp}. For example, an equality becomes
$g(x)-\tau_{\mathrm{feas}}\le0$ and
$-g(x)-\tau_{\mathrm{feas}}\le0$. Feasibility is then tested against zero,
not against the same positive tolerance a second time. Separate inequality
and equality tolerances can be absorbed into their respective rows under
the same convention. These tolerances define the problem being solved;
they are distinct from conservative floating-point error allowances.

\paragraph{Conservative numerical evaluation.}
The proofs below are exact-arithmetic statements. To use a computed row
$(\widehat a,\widehat c)$, suppose an error allowance satisfies
\[
| (\widehat a-a)^\top x+(\widehat c-c)|\le\epsilon_{\mathrm{aff}}(B)
\quad\text{on the required validity set}.
\]
Then $\widehat a^\top x+\widehat c-\epsilon_{\mathrm{aff}}(B)$ remains a
valid lower row. A row used for clipping must be weakened before it can
remove points; subtracting an error only from a final reported scalar
does not protect an earlier clipping step. Scalar minima must be rounded
or bounded downward, retained intervals outward, and incumbent values
conservatively relative to the original model. The allowances must cover
operator relaxations, graph translation where relevant, propagation,
coefficient mixing, and concretization. The existence and implementation
of such allowances are requirements, not consequences of writing a final
subtraction.

\subsection{Why Linear Mixing Can Improve a Bound}
\label{app:dual-proofs}
\label{app:gear_dual}

The relaxation targeted by dual ascent is the fixed-row LP in
\eqref{eq:node_lp}, not the original NLP or a lifted multi-facet LP.
For a nonempty finite box and fixed affine rows, set
\begin{equation}
a(\lambda)=a_0+A^\top\lambda,\qquad
c(\lambda)=c_0+c^\top\lambda.
\label{eq:mixed_row}
\end{equation}
The dual function is exactly \eqref{eq:dual_lb}, with inner domain $B$.

\begin{proposition}[Dual of the propagated-row LP]
\label{prop:dual}
For every $\lambda\ge0$,
$d_B(\lambda)\le p_{\mathrm{LP}}(B)\le p^\star(B)$.
If $\mathcal P(B)\ne\emptyset$, then
$\max_{\lambda\ge0}d_B(\lambda)=p_{\mathrm{LP}}(B)$.
If $\mathcal P(B)=\emptyset$, then $d_B$ is unbounded above.
\end{proposition}
\begin{proof}
For $x\in\mathcal P(B)$, the mixed objective is at most
$\underline f_0(x)$, proving weak duality. A nonempty $\mathcal P(B)$
is compact, so the finite LP optimum exists; LP strong duality applied
with the box bounds retained in the inner minimization gives equality.
For the empty case, the compact convex set
$\{Ax+c:x\in B\}$ is disjoint from the closed negative orthant.
A strict separating vector can be chosen as $\nu\ge0$ and satisfies
$\delta=\min_{x\in B}\nu^\top(Ax+c)>0$. Hence
$d_B(t\nu)\ge d_B(0)+t\delta\to+\infty$.
\end{proof}

\paragraph{Inner minimizer and supergradient.}
The inner minimization is coordinate-separable. One minimizing corner is
\[
x_{\lambda,j}=
\begin{cases}
\underline x_j,&a_j(\lambda)\ge0,\\
\overline x_j,&a_j(\lambda)<0.
\end{cases}
\]
If a coefficient is zero, any value in its interval is minimizing.
The function $d_B$ is concave and piecewise linear. Evaluating the mixed
objective for $\lambda'$ at $x_\lambda$ proves
\begin{equation}
d_B(\lambda')\le d_B(\lambda)
+(Ax_\lambda+c)^\top(\lambda'-\lambda).
\label{eq:dual-supergradient}
\end{equation}
Thus the row residual vector is a supergradient. It identifies constraints
violated by the inner minimizer, which need not be feasible for the LP.

\paragraph{Fixed-budget updates and inherited bounds.}
From a nonnegative initialization, take $T$ projected Adam ascent steps
using the supergradient, retaining the best conservative evaluation:
\begin{equation}
\mathrm{LB}_{\mathrm{new}}(B)=
\max_{0\le t\le T}\widehat d_B(\lambda^{(t)}),
\qquad\widehat d_B(\lambda)\le d_B(\lambda).
\label{eq:adam_step}
\end{equation}
The stored node bound is
$\max\{\mathrm{LB}_{\mathrm{parent}},\mathrm{LB}_{\mathrm{new}}(B)\}$.
Let $\lambda^{\mathrm{best}}$ maximize the newly evaluated values, and
$x^{\mathrm{best}}=x_{\lambda^{\mathrm{best}}}$; this is the warm-start
pair, even if the inherited scalar happens to be stronger. Evaluating
$\lambda=0$ provides an objective-only backstop. The theory of
\cref{app:termination} explicitly assumes a sufficiently accurate
backstop; its use is not inferred merely from a finite Adam budget.

\paragraph{An illustrative improvement.}
For $\min_{x\in[0,1]}x$ subject to $1/2-x\le0$,
\[
d_B(\lambda)=\lambda/2+\min\{0,1-\lambda\}.
\]
The objective-only value $d_B(0)=0$ ignores feasibility; $d_B(1)=1/2$
recovers the exact LP and NLP value. In general an improving direction
requires ruling out the entire minimizing face, not only one arbitrary
minimizing corner.

\paragraph{Active-face improvement.}
Let $X_0=\arg\min_{x\in B}\underline f_0(x)$ and
$h_i(x)=a_i^\top x+c_i$. For this finite polyhedral problem,
\begin{equation}
d_B'(0;e_i)=\min_{x\in X_0}h_i(x)=:\underline v_i.
\label{eq:app_directional_deriv}
\end{equation}
Let $\mathcal V(B)$ denote the vertices of $B$, and define
\begin{equation}
\rho_i^B=
\min_{\substack{v\in\mathcal V(B):\ \underline f_0(v)>d_B(0)\\
 h_i(v)<\underline v_i}}
\frac{\underline f_0(v)-d_B(0)}{\underline v_i-h_i(v)},
\label{eq:app_stability_radius}
\end{equation}
where the minimum over the empty set is $+\infty$.
\begin{theorem}[Strict linear-mixture improvement]
\label{thm:strict_improve}
For $0\le\eta<\rho_i^B$,
$d_B(\eta e_i)-d_B(0)=\eta\underline v_i$.
The minimizing set is contained in $X_0$; if $\underline v_i>0$,
sufficiently small positive steps strictly raise the bound.
For a cap $\Lambda=[0,\lambda_{\max}]^m$ with $\lambda_{\max}>0$,
\[
\max_{\lambda\in\Lambda}d_B(\lambda)-d_B(0)
\ge [\underline v_i]_+\min\{\rho_i^B/2,\lambda_{\max}\}.
\]
\end{theorem}
\begin{proof}
On $X_0$, minimizing the perturbed objective gives
$d_B(0)+\eta\underline v_i$. A vertex outside $X_0$ with
$h_i(v)\ge\underline v_i$ cannot improve it for $\eta\ge0$.
Any other vertex can first match it at the ratio in
\eqref{eq:app_stability_radius}. Before the minimum such ratio, an
optimizer therefore remains in $X_0$. Evaluating an admissible step
proves the final bound when $\underline v_i>0$; otherwise $\lambda=0$
proves its zero right-hand side. Monotonicity is asserted only on the
stated interval, not along the whole ray.
\end{proof}

\paragraph{Best-iterate rate for a reference optimizer.}
The following bound concerns projected supergradient ascent, not Adam.
\begin{theorem}[Best-iterate rate for projected supergradient ascent]
\label{thm:dual_ascent_rate}
Let $g$ be concave on a compact convex set $\Lambda$, let $\lambda^*$
maximize it, and let the chosen supergradients satisfy $\|s^t\|\le G$.
For $\lambda^{t+1}=\Pi_\Lambda(\lambda^t+\eta_t s^t)$ with
$\eta_t=\eta/\sqrt T$,
\[
g(\lambda^*)-\max_{1\le t\le T}g(\lambda^t)
\le\frac{\|\lambda^1-\lambda^*\|^2+\eta^2G^2}{2\eta\sqrt T}.
\]
With $\eta_t=\eta/\sqrt t$, the corresponding weighted-iterate argument
gives an $O(\log T/\sqrt T)$ best-iterate bound.
\end{theorem}
\begin{proof}
Projection is nonexpansive, and the supergradient inequality gives
\[
2\eta_t\bigl[g(\lambda^*)-g(\lambda^t)\bigr]
\le \|\lambda^t-\lambda^*\|^2
-\|\lambda^{t+1}-\lambda^*\|^2+\eta_t^2G^2.
\]
Sum over $t$ and divide by $2\sum_t\eta_t$; the weighted mean gap bounds
the smallest iterate gap. The stated schedules give the two bounds.
\end{proof}
A cap can lower the best attainable dual value relative to the unrestricted
LP dual. Projected Adam does not inherit this theorem; finite-iterate
validity follows independently from Proposition~\ref{prop:dual} and conservative
evaluation. Rescaling an LP multiplier to initialize ALM remains a
warm-start heuristic, not a claim of original-NLP stationarity.

\subsection{How Clipping Works in \GEAR{}}
\label{appendix::clipping}
\label{app:gear_clipping}

Both clipping procedures use currently available valid rows. Write
$\mathcal P=\{x\in B:Ax+c\le0\}$ for this temporary relaxation, which
may use inherited rather than freshly propagated rows. Its only required
relation to the original problem is $\mathcal F(B)\subseteq\mathcal P$.
The procedures are the relaxed and complete clipping mechanisms of
Zhou et al.~\cite{zhou2025clip}, applied to NLP constraint rows.

\paragraph{Relaxed clipping.}
Let $L_i=\min_{x\in B}(a_i^\top x+c_i)$. For $a_{i,j}>0$, the minimum
contribution of coordinate $j$ to $L_i$ is $a_{i,j}\underline x_j$.
Holding that coordinate at a feasible $x_j$ and minimizing all other
coordinates implies $L_i+a_{i,j}(x_j-\underline x_j)\le0$.
The analogous argument for $a_{i,j}<0$ gives
\begin{equation}
x_j\le\underline x_j-L_i/a_{i,j}\quad(a_{i,j}>0),\qquad
x_j\ge\overline x_j-L_i/a_{i,j}\quad(a_{i,j}<0).
\label{eq:relaxed_clipping}
\end{equation}
Intersect these bounds with the original box and across all rows. A
single sweep costs $O(mn)$ per box and yields
$\mathcal F(B)\subseteq\mathcal P\subseteq B'\subseteq B$.
Zero coefficients yield no coordinate update. A positive $L_i$ or crossed
coordinate bounds certifies infeasibility. The procedure computes a
containing box, not necessarily the tightest axis-aligned hull of the
joint polyhedron. Parallel sweeps and repeated sequential updates are
both sound when each uses the appropriate current intervals.

\paragraph{Complete clipping.}
For an intermediate lower row $\underline z(x)=a_z^\top x+c_z$, the
constraint-aware target is
\[
\min_{x\in\mathcal P}\underline z(x).
\]
Its dual is
\[
\max_{\beta\ge0}\min_{x\in B}
\bigl[a_z^\top x+c_z+\beta^\top(Ax+c)\bigr],
\]
with multipliers distinct from the objective multipliers $\lambda$.
Every evaluated dual value is a valid intermediate lower bound. For an
upper row $\overline z$, applying the lower-bound construction to
$-\overline z$ and negating the result gives a valid upper bound.
The single-constraint procedure of Zhou et al.~\cite{zhou2025clip} evaluates the
piecewise-linear dual by sorting breakpoints in $O(n\log n)$ time;
multiple constraints are treated through coordinate ascent. Finite
updates need not reach the joint LP optimum. Retaining the previous bound
alongside new candidates prevents an unfinished update from weakening
an intermediate interval.

\paragraph{Why the procedures differ.}
For $x_1,x_2\in[0,1]$ and $x_1+x_2\le1$, each coordinate can still
attain $1$, so no valid input hull can decrease either upper endpoint.
For $z=x_1+x_2$, however, the row yields $z\le1$ rather than the box bound
$2$. This tighter interval can strengthen every subsequent relaxation
that uses $z$. The resulting new rows need only hold for feasible graph
trajectories, as formalized in \cref{app:admissibility}.

\paragraph{Scheduling and infeasibility.}
Before a child's first propagation, its parent's rows and the original
affine rows are available for clipping. Bounds produced later in the pass
may be used only in subsequent computations or a repeated pass, not
retroactively in their own derivation. Empty input or intermediate
intervals certify absence of feasible trajectories. Empty intermediate
intervals need not establish emptiness of a separately defined
input-only polyhedron that omitted those intermediate conditions; the
required conclusion is $\mathcal F(B)=\emptyset$.

\subsection{Safe Linear Concretization over Clipped Relaxation Domains}
\label{app:concretization}

Box concretization and polyhedral concretization have different roles.
For an affine function $q(x)=a^\top x+c_0$,
\[
b_B(q)=c_0+\sum_j\min\{a_j\underline x_j,a_j\overline x_j\}
\]
is the exact box minimum in real arithmetic. The stronger polyhedral
target is the LP
\begin{equation}
\min_x\quad a^\top x+c_0
\qquad\text{s.t.}\quad
\underline x\le x\le\overline x,\quad Ax+c\le0.
\label{eq:exact-cut-lp}
\end{equation}
This LP explains the target of complete clipping and the objective dual.
The validity of the batched dual construction requires no LP solve; LP solvers only propose
multipliers (\cref{appsec:comp_bab}).

\paragraph{Domain containment.}
If clipping replaces $B$ by $B'$, the invariant is
\begin{equation}
\mathcal F(B)\subseteq\mathcal P\subseteq B'\subseteq B,
\qquad\mathcal F(B')=\mathcal F(B).
\label{eq:clip_containment}
\end{equation}
Later concretizations use the current box $B'$; inherited rows remain
valid on its feasible set. The polyhedron is a relaxation target, whereas
the inner minimization of $d_{B'}$ is still over the box. Minimizing
$\underline f_0+\lambda^\top(Ax+c)$ over $\mathcal P$ would impose the
same rows twice: their weighted sum is already nonpositive everywhere
in that domain, so its maximum over $\lambda\ge0$ cannot improve on
$\lambda=0$. This is not the dual used in \eqref{eq:dual_lb}.

\paragraph{Equivalent coordinate formulas.}
For a row $\alpha^\top x\le b$, relaxed clipping can also be written
\begin{equation}
\overline x_j\leftarrow\min\!\left\{\overline x_j,
\frac{b-\sum_{k\ne j}\min\{\alpha_k\underline x_k,
\alpha_k\overline x_k\}}{\alpha_j}\right\}
\qquad(\alpha_j>0),
\label{eq:lin-cut-upper-prop}
\end{equation}
\begin{equation}
\underline x_j\leftarrow\max\!\left\{\underline x_j,
\frac{b-\sum_{k\ne j}\min\{\alpha_k\underline x_k,
\alpha_k\overline x_k\}}{\alpha_j}\right\}
\qquad(\alpha_j<0).
\label{eq:lin-cut-lower-prop}
\end{equation}
Both use the minimum activity of all other coordinates. The sign of the
denominator determines the inequality direction. Early stopping leaves
a larger containing box and is therefore conservative.

\begin{table}[h]
\centering
\small
\setlength{\tabcolsep}{4pt}
\begin{tabular}{p{0.24\linewidth}p{0.33\linewidth}p{0.33\linewidth}}
\toprule
\textbf{Quantity} & \textbf{Target} & \textbf{Safe approximation} \\
\midrule
Objective bound & $d_B(\lambda)$, evaluated over the current box & Return no larger value \\
Infeasibility certificate & $\min_{x\in B}\nu^\top(Ax+c)$ & Return no larger value; prune only when positive \\
Intermediate lower bound & $\min_{x\in\mathcal P}\underline z(x)$ & Any feasible dual value, evaluated downward \\
Retained intervals & Contain every feasible trajectory & Round outward \\
Incumbent & Feasible value of the original objective & Check feasibility and objective directly \\
\bottomrule
\end{tabular}
\caption{Directions of conservative evaluation. Feasibility shifts have
already been absorbed into $f_i$; numerical margins are separate.}
\label{tab:oracle-safety-directions}
\end{table}

\subsection{Operator Relaxations and Validity Domains}
\label{app:operators}
\label{appsec:envelopes}
\label{app:gear_operators}

\begin{table}[h]
\centering
\footnotesize
\setlength{\tabcolsep}{4pt}
\begin{tabular}{@{}p{0.27\linewidth}p{0.31\linewidth}p{0.34\linewidth}@{}}
\toprule
Operator & Change relative to the baseline & Mechanism (this work) \\
\midrule
$x\!\cdot\!y$ (\textsc{BoundMul})        & revised default (upstream fixed one McCormick vertex)        & coefficient-aware $r\!\in\!\{0,1\}$ per element \\
$x/y$ on $y_l{>}0$ (\textsc{BoundDiv})   & new opt-in fused mode (upstream decomposes $x\!\cdot\!(1/y)$) & 8-plane candidate set with per-element coverage \\
$\log,\sqrt{\cdot},1/x,\exp$             & revised tangent side (upstream midpoint-style tangent)       & parallel-tangent (tangent slope $=$ secant slope) \\
$x^\alpha,\,\alpha\!\notin\!\mathbb{Z}_{\ge2}$ (\textsc{BoundPow}) & new branch & non-integer / $x\!\ge\!0$ branch with secant-tangent \\
\textsc{BoundScatterElements} / \textsc{BoundATenIndex\_add} & new operators & tensor-mode bound\_backward / bound\_forward / interval\_propagate \\
\textsc{ClampedMultiplication}           & revised                                                       & \texttt{where(A==0)} guards remove $0\!\cdot\!\infty$ NaN paths \\
\bottomrule
\end{tabular}
\caption{Operator changes relative to the baseline \textsc{auto\_LiRPA}
implementation. New and revised primitives retain a tensorized affine-bound
interface.}
\label{tab:operator_delta}
\end{table}

This section spells out the per-operator relaxations summarized by
\cref{tab:operator_delta} in \cref{subsec:operators}. Linear bound
propagation composes these local relaxations backward through the
computation graph; the rule is local but the resulting linear bound is
valid on the retained feasible trajectories, because every local
replacement is valid on its recorded interval domain and selected
according to the backward-coefficient sign.

\begin{lemma}[Composition of local relaxations]
\label{lem:operator-composition}
Let $S$ be a set of inputs whose graph trajectories lie within all current intermediate intervals. Let every nonlinear graph node $v=h(u)$ be equipped, on its current interval domain $U$, with linear functions $\underline h$ and $\overline h$ satisfying $\underline h(u)\le h(u)\le\overline h(u)$ for all $u\in U$.
If backward linear bound propagation for a lower bound replaces each occurrence of $\gamma h(u)$ by $\gamma\underline h(u)$ for $\gamma\ge0$ and by $\gamma\overline h(u)$ for $\gamma<0$, then the final expression at the input is an affine lower bound for the selected graph output on $S$. Ordinary propagation takes $S=B$; constraint-tightened propagation may take $S=\mathcal F(B)$.
\end{lemma}

The proof is a reverse topological induction using \cref{eq:sign-compatible-composition}.
The same argument applies componentwise to batched coefficients and to multi-input operators after flattening their local linear relaxations into the parent variables.

\paragraph{Linear and separable arithmetic.}
Affine nodes are exact.
For $z=\alpha x+\beta y+c$, the backward coefficient on $z$ is distributed as $\alpha\gamma$ to $x$, $\beta\gamma$ to $y$, and $\gamma c$ to the bias.
No relaxation gap is introduced.
Elementwise negation, scalar multiplication, addition, subtraction, indexing, reshaping, and broadcasting are handled in this exact category after canonicalization.
This exact propagation matters because it keeps large linear subgraphs from paying relaxation cost before nonlinear operators are reached.

\paragraph{Univariate secant-tangent template.}
Let $z=h(t)$ on $t\in[\alpha,\beta]$ with $\alpha<\beta$, and define the secant
\begin{equation}
S_{\alpha,\beta}^h(t)
=h(\alpha)+\frac{h(\beta)-h(\alpha)}{\beta-\alpha}(t-\alpha).
\label{eq:univariate-secant}
\end{equation}
For any $t_0\in[\alpha,\beta]$, define the tangent
\begin{equation}
T_{t_0}^h(t)=h(t_0)+h'(t_0)(t-t_0).
\label{eq:univariate-tangent}
\end{equation}
If $h$ is convex on $[\alpha,\beta]$, then $T_{t_0}^h(t)\le h(t)\le S_{\alpha,\beta}^h(t)$.
If $h$ is concave, the inequalities reverse:
$S_{\alpha,\beta}^h(t)\le h(t)\le T_{t_0}^h(t)$.
The implementation chooses $t_0$ by the secant-slope rule
\begin{equation}
h'(t_0)=\frac{h(\beta)-h(\alpha)}{\beta-\alpha}
\label{eq:secant-slope-rule}
\end{equation}
when the solution lies in the interval; otherwise it uses the endpoint tangent that gives the smaller valid bound gap.
For $h\in C^2$ with a bounded second derivative on the interval, the one-dimensional gap is $O((\beta-\alpha)^2)$. This local estimate is distinct from the graph-wide accuracy assumption in Assumption~\ref{ass:relaxation_modulus}.

\paragraph{Logarithm.}
For $h(t)=\log t$, the domain guard is $0<\alpha\le t\le\beta$.
The function is concave, so the lower side is the secant and the upper side is a tangent.
The secant-slope tangent point is the logarithmic mean
\begin{equation}
t_0=\frac{\beta-\alpha}{\log\beta-\log\alpha}.
\label{eq:log-tangent-point}
\end{equation}
If the lower endpoint is at or below zero after outward interval tolerances, the node must be split, rewritten with a stricter domain guard, or rejected; clipping the value of $\log t$ without a corresponding domain proof would invalidate the linear bound.

\paragraph{Square root.}
For $h(t)=\sqrt t$, the guard is $0\le\alpha\le t\le\beta$.
The function is concave on the nonnegative axis, so the lower side is the secant and the upper side is a tangent.
For $\alpha<\beta$, the secant-slope tangent point is
\begin{equation}
t_0=\left(\frac{\sqrt{\alpha}+\sqrt{\beta}}{2}\right)^2 .
\label{eq:sqrt-tangent-point}
\end{equation}
When $\alpha=0$, the tangent slope at the endpoint is singular; the implementation uses the interior point from \cref{eq:sqrt-tangent-point} for the upper side and keeps the secant for the lower side.
If numerical interval propagation produces $\alpha<0$ for a source expression that is mathematically nonnegative, the linear-bound engine cannot use that fact unless it is represented as a verified bound; otherwise the expression is split or rejected.

\paragraph{Reciprocal.}
For $h(t)=1/t$, the interval must be sign-fixed and bounded away from zero.
On $\alpha>0$, $h$ is convex; tangents are lower bounds and the secant is an upper bound.
On $\beta<0$, $h$ is concave; the same secant-tangent construction is used with the lower/upper sides reversed.
In both cases the slope-matching condition $h'(t_0)=(h(\beta)-h(\alpha))/(\beta-\alpha)$ has the
sign-aware solution
\begin{equation}
t_0=\operatorname{sign}(\alpha)\sqrt{\alpha\beta}
\quad\Longleftrightarrow\quad
t_0=
\begin{cases}
+\sqrt{\alpha\beta}, & \alpha>0,\\
-\sqrt{\alpha\beta}, & \beta<0,
\end{cases}
\label{eq:reciprocal-tangent-point}
\end{equation}
since $\alpha,\beta$ have the same sign under the sign-fixed guard.
The naive form $t_0=\sqrt{\alpha\beta}$ is positive on negative intervals and would lie outside $[\alpha,\beta]$, producing an invalid relaxation; the sign correction restores validity.
Intervals crossing zero are never relaxed by a finite linear bound.
The solver must split the denominator interval, rewrite the model, or mark the instance unsupported.

\paragraph{Exponential.}
For $h(t)=\exp t$, the function is convex on every finite interval.
The upper side is the secant and the lower side is a tangent.
The secant-slope tangent point is
\begin{equation}
t_0=\log\frac{e^\beta-e^\alpha}{\beta-\alpha}.
\label{eq:exp-tangent-point}
\end{equation}
When the interval is extremely narrow, the implementation evaluates this expression with a stable
expm1-style form and otherwise falls back to the midpoint tangent only if the slope-matching
calculation is numerically indistinguishable from the midpoint at floating-point precision.

\begin{figure}[t]
\centering
\includegraphics[width=\linewidth]{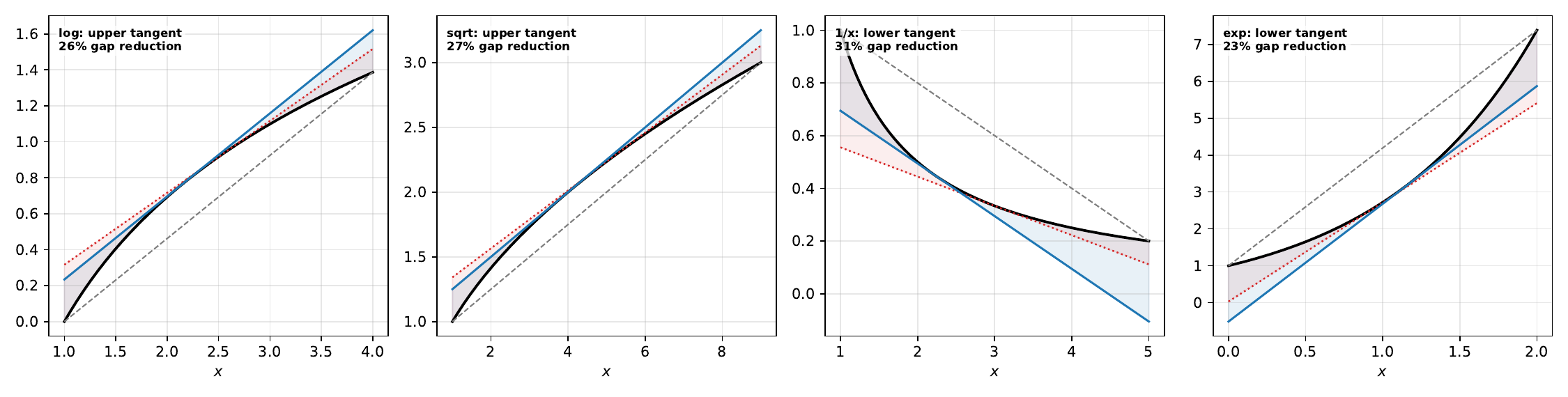}
\caption{Tangent-side comparison for four convex/concave operators. Black is the operator, gray is
the unchanged secant, red is the upstream midpoint-style tangent, and blue is the \GEAR{}
parallel-tangent. For concave $\log$ and $\sqrt{\cdot}$ the changed side is the upper tangent; for
convex $1/x$ and $\exp$ it is the lower tangent.}
\label{appfig:envelope_parallel_tangent}
\end{figure}

\paragraph{Power functions.}
For $h(t)=t^p$ with non-integer $p$, the base interval must be nonnegative and, when derivatives are singular at zero, guarded away from zero according to the implementation tolerance.
For $0<p<1$, $h$ is concave on $[0,\infty)$ and uses the same lower-secant / upper-tangent pattern as $\sqrt t$.
For $p>1$, $h$ is convex on $[0,\infty)$ and uses lower tangents / upper secants.
For negative powers, the implementation treats $t^p$ as a reciprocal-power primitive on a positive sign-fixed interval and applies the corresponding convexity rule; if the base can approach zero, no finite linear bound is emitted.
Integer powers with odd exponents can be represented on sign-changing intervals, but their curvature changes at zero for most exponents; the interval is split at curvature changes before applying the monotone-curvature template.

\paragraph{Selection among valid affine planes.}
For $z=h(y)$ on a box $Y$ and backward coefficient $\gamma$, let
$\mathcal L$ and $\mathcal U$ contain valid lower and upper planes.
The rule used in the main text is
\begin{equation}
q^\star\in
\begin{cases}
\displaystyle\arg\max_{q\in\mathcal L}\min_{y\in Y}q(y),&\gamma\ge0,\\[2pt]
\displaystyle\arg\min_{q\in\mathcal U}\max_{y\in Y}q(y),&\gamma<0.
\end{cases}
\label{eq:coeff_aware_lower}
\end{equation}
This maximizes $\min_{y\in Y}\gamma q(y)$ among the eligible candidates.
For $\gamma=0$ all choices have zero contribution. Each score is a box
concretization, so selection can be tensorized. The guarantee is local:
other terms, shared intermediate variables, and coefficient cancellation
can change which selection gives the best final composed bound.

\paragraph{Bilinear terms.}
For $z=xy$ on $X\times Y=[x_L,x_U]\times[y_L,y_U]$, the McCormick facets are
\begin{align}
z&\ge x_Ly+y_Lx-x_Ly_L, &
z&\ge x_Uy+y_Ux-x_Uy_U, \label{eq:mccormick-lower}\\
z&\le x_Uy+y_Lx-x_Uy_L, &
z&\le x_Ly+y_Ux-x_Ly_U. \label{eq:mccormick-upper}
\end{align}
An LP relaxation would keep all four facets simultaneously.
\GEAR{} keeps one linear bound per output direction, so it must project these facets to one lower or one upper linear plane for each backward coefficient.
For a positive parent coefficient $\gamma$, lower candidates in \cref{eq:mccormick-lower} are scored by the scalar linear minimum of $\gamma p(x,y)$ on the current local rectangle or cut-tightened box.
For a negative parent coefficient, upper candidates in \cref{eq:mccormick-upper} are scored after multiplication by $\gamma$.
The selected plane is sound because it is one of the McCormick facets, but the projection loses the simultaneous-facet closure; this is the precise reason LP-native solvers can dominate on bilinear-heavy root relaxations.
The maximum McCormick envelope width over the rectangle is proportional to $(x_U-x_L)(y_U-y_L)$, so spatial refinement closes this local gap quadratically in the box diameter when both factors shrink.

\begin{figure}[t]
\centering
\includegraphics[width=0.96\linewidth]{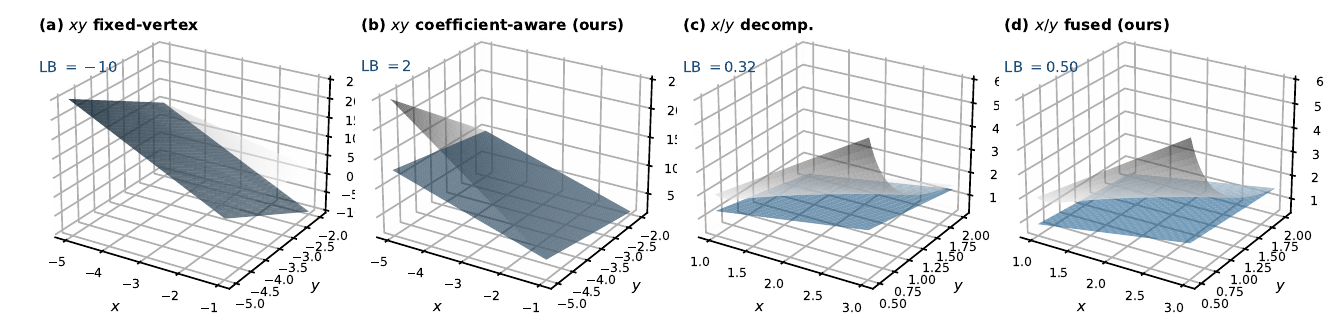}
\caption{Linear lower bounds for two bivariate operators. On
$[-5,-1]\times[-5,-2]$, the fixed product plane gives $-10$ and the
selected McCormick plane gives $2$. On $[1,3]\times[0.5,2]$, the
reported reciprocal-product decomposition gives $0.32$ and the fused
quotient plane gives $0.5$. The selected bounds equal the true minima
on these illustrated rectangles.}
\label{fig:envelope_binary_ops}
\end{figure}

\paragraph{Fused division.}
For $z=x/y$, the denominator interval must satisfy either $y_L\ge\delta_{\mathrm{div}}>0$ or $y_U\le-\delta_{\mathrm{div}}<0$.
The fused implementation constructs candidate planes $p(x,y)=\alpha x+\beta y+\gamma$ directly for the quotient on the guarded rectangle.
On the positive-denominator branch used in \cref{fig:envelope_binary_ops}, the candidate pool has
eight planes: four corner tangents
\begin{equation}
P_T^{(x^\star,y^\star)}(x,y)
=\frac{x}{y^\star}-\frac{x^\star y}{(y^\star)^2}+\frac{x^\star}{y^\star},
\qquad (x^\star,y^\star)\in\{x_L,x_U\}\times\{y_L,y_U\},
\label{eq:division-corner-tangent}
\end{equation}
and four three-corner secant planes obtained by interpolating $x/y$ at any three corners of
$[x_L,x_U]\times[y_L,y_U]$.
The negative-denominator branch is handled by the same guarded-quotient construction with the
inequality orientation reversed by the denominator sign.
Because $x/y$ is not globally convex or concave on an arbitrary guarded rectangle, candidate
generation alone is not trusted: every candidate is accepted only after a residual validity check.
For $y>0$, a lower plane must satisfy
\begin{equation}
p(x,y)\le x/y
\quad\Longleftrightarrow\quad
R_p(x,y):=y\,p(x,y)-x\le0
\quad \forall (x,y)\in X\times Y,
\label{eq:division-lower-check-positive}
\end{equation}
and an upper plane must satisfy $R_p(x,y)\ge0$.
For $y<0$, the inequalities reverse because multiplication by $y$ changes sign.
The residual $R_p$ is a quadratic-bilinear expression; the implementation checks the required sign with interval bounds, McCormick bounds on the $xy$ term, and exact endpoint checks for the $y^2$ term.
If the residual check is inconclusive, the plane is rejected and the engine falls back to a reciprocal-plus-bilinear relaxation, interval splitting, or unsupported-instance reporting.
This accept-after-check rule is the soundness guard for fused division.
After invalid candidates are filtered out, the surviving lower and upper planes are scored by the
same box-extremum rule as \cref{eq:coeff_aware_lower}.

\paragraph{Trigonometric functions.}
For $\sin t$ and $\cos t$, global secant-tangent relaxations over long intervals are too loose and may cross curvature changes.
The interval is first localized by splitting at points where the second derivative changes sign:
\[
\sin t:\ t=k\pi,\qquad
\cos t:\ t=\pi/2+k\pi,\qquad k\in\mathbb{Z}.
\]
On each resulting monotone-curvature segment, \cref{eq:univariate-secant,eq:univariate-tangent} applies.
If an interval spans multiple periods and splitting would exceed the configured segment budget, a valid coarse interval bound or a documented fallback is required; an unsupported primitive cannot generate a pruning certificate.
A local segment bound may be forwarded only when the retained feasible trajectories lie within its validity domain.

\paragraph{Max, min, and residual outputs.}
Nonlinear NLP constraints are represented as separate graph outputs and affine constraints as exact rows, not through a hard max in the objective.
For residual infeasibility, \GEAR{} avoids needing a tight linear bound for the nonsmooth max by reusing the fact that every constraint lower bound is also a lower bound for $r(x)=\max_i f_i(x)$.
If a source model contains explicit $\max$ or $\min$ primitives, they must be rewritten into inequalities, split into supported piecewise-smooth regions, or handled by a separate valid relaxation rule before they can participate in pruning.

\paragraph{Validity-domain bookkeeping.}
Every local relaxation records the interval on which it is valid.
The output row is valid on the set of inputs whose trajectories remain within these intervals; after constraint-driven clipping, validity on $\mathcal F(B)$ is sufficient.
When a row is forwarded, each child's feasible set is a subset of the parent's, so the validity implication is preserved. Fresh rows replace the inherited nonlinear rows after the child's pass.
No clipping or operator-domain guard may exclude an original feasible trajectory without a valid certificate; otherwise the affected relaxation must be recomputed or not used for pruning.

\begin{table}[h]
\centering
\small
\setlength{\tabcolsep}{3pt}
\begin{tabular}{p{0.15\linewidth}p{0.24\linewidth}p{0.27\linewidth}p{0.22\linewidth}}
\toprule
\textbf{Primitive} & \textbf{Required domain} & \textbf{Relaxation rule} & \textbf{Failure action} \\
\midrule
linear arithmetic & finite input intervals & exact coefficient propagation & none \\
$xy$ & finite rectangle & one coefficient-aware McCormick facet & keep loose facet or split \\
$x/y$ & denominator sign-fixed and guarded & accept-after-check fused quotient plane & reciprocal rewrite, split, or reject \\
$1/t$ & sign-fixed and guarded & secant-tangent with side by curvature & split or reject near zero \\
$\log t$ & $t>0$ with guard & concave lower secant / upper tangent & split, strengthen guard, or reject \\
$\sqrt t$ & $t\ge0$ with derivative guard if needed & concave lower secant / upper tangent & split or reject if nonnegativity is not established \\
$t^p$ & domain compatible with exponent & curvature-local secant-tangent & split at sign/curvature changes or reject \\
$\sin t,\cos t$ & localized monotone-curvature segment & segmentwise secant-tangent & split by period/curvature or reject \\
\bottomrule
\end{tabular}
\caption{Operator validity conditions for the linear-bound engine. Outside these domains, \GEAR{} must split, rewrite, fall back to a looser valid rule, or reject the primitive; it must not emit an unchecked linear bound.}
\label{tab:operator-validity}
\end{table}

\subsection{Primal Recovery Details}
\label{app:alm-details}
\label{appsec:adam_details}
\label{app:gear_primal}

The primal heuristic evaluates the original graph and the original affine
constraints, not only their relaxation. It applies projected first-order
steps to the Hestenes-Powell-Rockafellar function in
\eqref{eq:hpr_alm}, followed by multiplier and penalty updates. These
operations seek a useful incumbent; successful solution of the local NLP
is not a prerequisite for the lower-bound certificates.

\begin{algorithm}[tbp]
\caption{\textsc{ALM-PGD} on the current box $B$}
\label{alg:alm_pgd}
\begin{algorithmic}[1]
\REQUIRE Current box $B=[\underline x,\overline x]$;
original objective $f_0$ and tolerance-shifted constraints $f=(f_1,\ldots,f_m)$;
seed $x\in B$; $\mu\ge0$, $\rho>0$; outer and inner budgets
$K_{\mathrm{out}},K_{\mathrm{in}}$; step schedule $\{\eta_t\}$.
\STATE $\bar x\leftarrow\varnothing$
\FOR{$k=1,\ldots,K_{\mathrm{out}}$}
    \FOR{$t=1,\ldots,K_{\mathrm{in}}$}
        \STATE $g\leftarrow\nabla_x\phi_\rho(x;\mu)$, with
        $\phi_\rho$ defined in \eqref{eq:hpr_alm}
        \STATE $x\leftarrow\Pi_B(\operatorname{AdamStep}(x,g,\eta_t))$
    \ENDFOR
    \STATE Evaluate $f_0(x)$ and all original shifted constraints at $x$
    \IF{$f_i(x)\le0$ for every $i$ and
    ($\bar x=\varnothing$ or $f_0(x)<f_0(\bar x)$)}
        \STATE $\bar x\leftarrow x$
    \ENDIF
    \STATE $\mu\leftarrow[\mu+\rho f(x)]_+$
    \STATE Increase $\rho$ if the maximum violation failed to decrease
    by the prescribed ratio
\ENDFOR
\ENSURE Best accepted candidate $\bar x$ among the tested iterates, or $\varnothing$.
\end{algorithmic}
\end{algorithm}

\paragraph{Multiplier update and gradient.}
The outer update is
\begin{equation}
\mu\leftarrow[\mu+\rho f(x)]_+.
\label{eq:app_alm_update}
\end{equation}
For differentiable original functions, the inner gradient is
\begin{equation}
\nabla_x\phi_\rho(x;\mu)=\nabla f_0(x)
+\rho\sum_{i=1}^m[f_i(x)+\mu_i/\rho]_+\nabla f_i(x).
\label{eq:alm-gradient}
\end{equation}
The squared positive-part penalty is continuously differentiable in its
scalar argument; its derivative is zero at the hinge. Any nonsmoothness
of the original graph still requires the implementation's derivative
convention. Affine constraints contribute their exact gradients $a_i$.
The objective remains in every inner step, so the method does not switch
to a feasibility-only objective upon reaching the constraint set.

\paragraph{Projected Adam inner loop.}
For an inner step $t\ge1$ with $g_t=\nabla_x\phi_\rho(x_t;\mu)$,
a standard projected-Adam update is
\begin{align}
m_t&=\beta_1m_{t-1}+(1-\beta_1)g_t,\\
v_t&=\beta_2v_{t-1}+(1-\beta_2)g_t\odot g_t,\\
\widehat m_t&=m_t/(1-\beta_1^t),\qquad
\widehat v_t=v_t/(1-\beta_2^t),\\
\widetilde x_{t+1}&=x_t-\eta_t
\frac{\widehat m_t}{\sqrt{\widehat v_t}+\epsilon_{\mathrm{adam}}},\\
x_{t+1}&=\Pi_B(\widetilde x_{t+1}).
\label{eq:projected-adam}
\end{align}
All vector operations and box projections batch across starts and boxes.
The inner solver uses \texttt{AdamClipping\_NLP};
its learning-rate normalization and handling of clipped coordinates are
implementation details beyond the certificate proof. Adam's second-moment
scaling is applied before projection. A smooth box parameterization is an
alternative implementation path, not assumed by the projected update
above.

\paragraph{Dual-guided seeds.}
The distinguished seed is
\[
x^{\mathrm{best}}\in\arg\min_{x\in B}
\bigl[\underline f_0(x)+(\lambda^{\mathrm{best}})^\top(Ax+c)\bigr].
\]
This is the box inner minimizer of the evaluated dual, not an optimizer
over $\mathcal P(B)$ and not necessarily an LP-feasible point. Additional
starts may use the midpoint, a projected incumbent,
Sobol points, and minimizers of individual constraint rows. Their use and
budget affect search success, not whether an accepted point is feasible.

\paragraph{Multiplier initialization and penalty schedule.}
A nonnegative rescaling can transfer constraint emphasis from the LP dual
to the ALM multipliers, for example
\[
\mu_i^{(0)}=\kappa\lambda_i^{\mathrm{best}},\qquad
\kappa=\frac{1+\|\nabla f_0(x^{\mathrm{best}})\|_\infty}
{1+\|\sum_i\lambda_i^{\mathrm{best}}\nabla f_i(x^{\mathrm{best}})\|_\infty}.
\]
Zero initialization is an alternative. The
rescaling is a warm start, not a statement that LP and original-NLP
multipliers coincide. The penalty is increased when violation reduction
is insufficient. Neither this schedule nor finite projected-Adam steps
establish a guarantee of feasible recovery on every nonconvex instance.

\paragraph{Use of linear feasibility cuts.}
Retained rows can additionally guide the primal search through
\begin{equation}
\phi_{\rho,\sigma}^{\mathrm{cut}}(x;\mu)
=\phi_\rho(x;\mu)+\frac{\sigma}{2}\|[Ax+c]_+\|_2^2,
\qquad\sigma\ge0.
\label{eq:cut-soft-penalty}
\end{equation}
Every original feasible point satisfies the rows, so the extra term
vanishes on the feasible set. It can nevertheless change trajectories
away from feasibility. This is the ``ALM with linear constraints''
variant; no extra penalty is used as a substitute for checking the
original nonlinear constraints.

\paragraph{Acceptance gate.}
A candidate improves the incumbent only when
\begin{equation}
\bar x\in B_0,\qquad f_i(\bar x)\le0\ (i=1,\ldots,m),\qquad
f_0(\bar x)<\mathrm{UB}.
\label{eq:feasibility_gate}
\end{equation}
All functions here use the same tolerance shifts as in the lower-bound
problem. Objective and feasibility are evaluated afresh; neither a low
ALM value nor satisfaction of the linear rows is sufficient. For numerical
certification the original-model checks and conservative evaluation
requirements are those of \cref{app:admissibility,appsec:translation_soundness}.

\subsection{Branching Calibration}
\label{app:branching}

\SBReweight{} branches on the coordinate that maximizes the $\lambda$-weighted score $S_j$ of
\cref{eq:multiplier_score}, a $\lambda$-weighted form of smart branching. Its fixed-row
interpretation explains the coefficient magnitudes without asserting
that the score predicts the nonlinear bounds reconstructed on children.

\paragraph{Fixed-row midpoint calculation.}
For $q(x)=a^\top x+c_0$ on $B$, a minimizing coordinate is
$e_j=\underline x_j$ when $a_j\ge0$ and $e_j=\overline x_j$ otherwise.
A midpoint split leaves this coordinate's contribution unchanged in
one child and increases its minimum by $|a_j|w_j/2$ in the other, where
$w_j=\overline x_j-\underline x_j$. For a uniformly chosen child of
this midpoint split,
\begin{equation}
\mathbb E[\Delta_j]=\frac{|a_j|w_j}{4},\qquad
\Delta_j\in\{0,|a_j|w_j/2\}.
\label{eq:expected_reduction}
\end{equation}
The minimum across both children is unchanged for this fixed affine row;
refinement is useful because later propagation rebuilds the nonlinear
relaxations, not because splitting a fixed affine function necessarily
raises the global lower bound.

\begin{lemma}[One-step linear branching proxy]
\label{lem:branching_proxy}
For the mixed row with coefficient
$a_j(\lambda)=a_{0,j}+\sum_i\lambda_i a_{i,j}$, the average fixed-row
midpoint improvement is bounded by $S_j/4$, since
\[
|a_j(\lambda)|\le |a_{0,j}|+\sum_i\lambda_i|a_{i,j}|.
\]
Here $S_j$ is evaluated using the same nonnegative multipliers as the
mixed row.
\end{lemma}
\begin{proof}
Apply the triangle inequality and then \eqref{eq:expected_reduction}.
\end{proof}
The implemented score avoids cancellation between outputs. This is a
proxy for where refinement may matter; it is not strong branching or
an exact estimate of the children's recomputed LP values.

\paragraph{Eligible coordinates and shrinking.}
For $\tau\in(0,1]$, choose
\begin{equation}
j^\star\in\arg\max_{j:\,w_j\ge\tau\max_k w_k}S_j.
\label{eq:hybrid_branching}
\end{equation}
If a positive maximum width remains, midpoint-split the selected
coordinate. Along any infinite branch, suppose the maximum width did
not tend to zero. Then infinitely many selected widths would exceed
a fixed positive number. At least one of finitely many coordinates
would be selected infinitely often at such widths, contradicting its
repeated halving. Clipping only reduces widths further. Thus every
infinite branch shrinks to a point. This argument concerns shrinking,
not feasible recovery or finite termination on its own.

\subsection{Overall Algorithm of \GEAR{}}
\label{appendix::gear}

\Cref{alg:driver} expands the main-text description with explicit initialization
and global-gap bookkeeping. Node operations are batched even though the
notation describes their per-box effect. The queue policy is an
implementation choice; the termination theorem additionally requires
that each batch include a minimum-bound box. No queue policy can invalidate
a certificate if all unresolved boxes remain represented.

\begin{algorithm}[H]
\caption{\GEAR{}: batched spatial branch-and-bound}
\label{alg:driver}
\label{alg:gear_main}
\begin{algorithmic}[1]
\REQUIRE Computation graph $F$,
initial finite box $B_0$, gap tolerance $\varepsilon>0$, batch size $K$,
and time limit.
\STATE Initialize $\mathrm{UB}\leftarrow+\infty$, $\bar x\leftarrow\varnothing$,
and $\mathcal Q\leftarrow\{B_0\}$ with $\mathrm{LB}(B_0)=-\infty$
and the original affine rows available for clipping.
\WHILE{$\mathcal Q\ne\emptyset$ and time remains}
    \STATE $L_{\mathrm{glob}}\leftarrow
    \min\{\mathrm{UB},\min_{B\in\mathcal Q}\mathrm{LB}(B)\}$
    \IF{$\mathrm{UB}<+\infty$ and $\mathrm{UB}-L_{\mathrm{glob}}\le\varepsilon$}
        \STATE \textbf{break}
    \ENDIF
    \STATE Pop up to $K$ boxes; they stay unresolved
    until their processing finishes.
    \STATE Apply relaxed clipping with available rows; discard certified
    infeasible boxes.
    \STATE Propagate nonlinear output bounds with complete clipping using
    already available rows; refresh the nonlinear constraint rows.
    \STATE Prune boxes with empty intermediate intervals or a positive
    conservative certificate from \eqref{eq:infeas_cert}.
    \STATE Evaluate fixed-budget dual updates and set each surviving bound
    to the maximum of its inherited value and \eqref{eq:adam_step}.
    \STATE Prune boxes with $\mathrm{LB}(B)\ge\mathrm{UB}$ when an incumbent exists.
    \STATE Run scheduled ALM-PGD searches using the dual warm starts;
    update $(\mathrm{UB},\bar x)$ only by \eqref{eq:feasibility_gate}.
    \STATE Recheck dominance under the updated incumbent.
    \STATE Midpoint-split survivors using \eqref{eq:hybrid_branching};
    push children with inherited lower bounds and fresh nonlinear rows.
\ENDWHILE
\STATE $L_{\mathrm{glob}}\leftarrow
\min\{\mathrm{UB},\min_{B\in\mathcal Q}\mathrm{LB}(B)\}$,
with $\min\emptyset=+\infty$.
\ENSURE Incumbent $(\mathrm{UB},\bar x)$ if found, global lower bound
$L_{\mathrm{glob}}$, and status: gap certified, infeasible, or limit reached.
\end{algorithmic}
\end{algorithm}

At a time-limit interruption within a batch, unfinished boxes retain
valid bounds and are included in the final active set. An unsupported
operator, an inconclusive numerical test, or an unsplittable unresolved
box is not an infeasibility certificate: it must remain unresolved or
cause an explicitly uncertified termination. In particular, no loop
condition evaluates $+\infty-(+\infty)$ before checking that an incumbent
exists.

\subsection{Node Lifecycle and Data Dependencies}
\label{app:node-lifecycle}
\label{app:gear_bab}

A node stores its current box, a scalar inherited lower bound, and the
currently available nonlinear constraint rows. Original affine rows
are shared. The following order makes the data dependencies explicit;
iteration counts and the primal-call schedule remain implementation
parameters.

\paragraph{Step 1: assemble available rows.}
At the root use the original affine constraints; at descendants combine
them with the parent's nonlinear rows. The latter need only be valid
on the parent's feasible set, which contains the child's feasible set.

\paragraph{Step 2: clip the input box.}
Apply relaxed clipping to the available rows. Retain the resulting box
and discard the node only if a conservative contradiction is obtained.
An inability to certify emptiness is not evidence of feasibility.

\paragraph{Step 3: propagate and tighten intermediate bounds.}
Process the nonlinear graph over the clipped box, using available rows
for complete clipping. Replace inherited nonlinear rows only after the
new pass has produced valid rows. A further propagation sweep may use
those refreshed rows; a row is never used before it exists.

\paragraph{Step 4: test infeasibility.}
Test individual rows or selected nonnegative combinations using
\eqref{eq:infeas_cert}, and check empty intermediate intervals. These
certificates do not require a finite incumbent.

\paragraph{Step 5: update the objective lower bound.}
Optimize $\lambda$ over the cached coefficients. Set the node bound to
the maximum of its inherited bound and the best conservatively evaluated
new dual value. Prune if it is at least the finite incumbent $\mathrm{UB}$.
The exact-dominance convention is used throughout this algorithm; stopping
at a finite global gap is a separate test.

\paragraph{Step 6: search for a primal incumbent.}
Run the prescribed primal searches on selected surviving boxes. Accept
only candidates satisfying \eqref{eq:feasibility_gate}. The linear rows
and dual outputs are heuristic guidance, not replacements for the
original constraints.

\paragraph{Step 7: recheck dominance and prepare the children.}
A newly improved incumbent can prune other surviving boxes. For those
remaining, retain the fresh nonlinear rows, original affine rows, and
best scalar lower bound. Row storage does not grow with tree depth.

\paragraph{Step 8: split and return children.}
Use \eqref{eq:hybrid_branching} and midpoint splitting. Each child inherits
the parent's lower bound and current nonlinear rows. The children cover
the clipped parent, which still contains all of its feasible points.

\begin{table}[h]
\centering
\small
\setlength{\tabcolsep}{4pt}
\begin{tabular}{p{0.21\linewidth}p{0.22\linewidth}p{0.46\linewidth}}
\toprule
\textbf{Action} & \textbf{Effect} & \textbf{Validity dependency} \\
\midrule
Clipping contradiction & Infeasibility pruning & Feasible trajectories are preserved \\
Positive row mixture & Infeasibility pruning & Nonnegative weights and valid constraint rows \\
Objective comparison & Dominance pruning & Valid lower bound and feasible incumbent \\
Primal recovery & Incumbent update & Original-model feasibility gate \\
Row inheritance & Reuse before propagation & Feasible-set containment \\
Spatial split & Refine active coverage & Child boxes cover the surviving parent \\
\bottomrule
\end{tabular}
\caption{Node operations and their certificate dependencies.}
\label{tab:node-dependencies}
\end{table}

\subsection{Proofs of the Sound Node Invariants}
\label{app:sound-invariants}

All local pruning rules follow from feasible-set preservation and valid
one-sided bounds. The claims hold even if primal search fails or dual
optimization terminates early.

\paragraph{Constraint validity and clipping.}
For $x\in\mathcal F(B)$,
$\underline f_i(x)\le f_i(x)\le0$, so
$\mathcal F(B)\subseteq\mathcal P(B)$. Coordinate clipping follows from
necessary conditions of these inequalities. An intermediate dual lower
bound is at most $\underline z(x)\le z(x)$ at every feasible point, and
an upper bound follows symmetrically. Thus interval tightening retains
all feasible trajectories. Reapplying backward linear bound propagation to those
intervals preserves output-row validity on $\mathcal F(B)$.

\paragraph{Objective lower bound.}
For every feasible $x$ and every $\lambda\ge0$,
\begin{equation}
d_B(\lambda)\le
\underline f_0(x)+\lambda^\top(Ax+c)
\le f_0(x)+\sum_i\lambda_i f_i(x)\le f_0(x).
\label{eq:app_combined}
\end{equation}
Taking an infimum over feasible points gives $d_B(\lambda)\le p^\star(B)$.
The maximum of any collection of valid lower bounds is valid, including
bounds inherited from a parent. A multiplier cap changes the values
reachable by the optimizer but does not invalidate any nonnegative
iterate.

\paragraph{Infeasibility.}
If $\nu\ge0$, every feasible point has
$\nu^\top(Ax+c)\le0$. A strictly positive box minimum therefore
contradicts feasibility. The normalized residual interpretation in
\cref{app:batched-linear-bound-state} gives the same conclusion. Failure
to find such a vector in a finite search establishes neither feasibility
nor the absence of a stronger certificate.

\paragraph{Incumbent validity and active coverage.}
The acceptance gate ensures that a finite incumbent is feasible for the
same tolerance-shifted problem used by the bounds. At an iteration
boundary, let $\mathcal Q$ contain every unresolved box. Then
\begin{equation}
\begin{gathered}
\mathrm{LB}(B)\le p^\star(B)\quad(B\in\mathcal Q),\\
\mathrm{UB}=f_0(\bar x)\text{ for a checked feasible }\bar x,
\quad\text{or}\quad\mathrm{UB}=+\infty,\\
\{x\in\mathcal F(B_0):f_0(x)<\mathrm{UB}\}
\subseteq\bigcup_{B\in\mathcal Q}B.
\end{gathered}
\label{eq:bab_invariants}
\end{equation}
Clipping preserves the third relation; infeasibility pruning removes no
feasible point; exact dominance removes only points unable to improve
the incumbent; and children cover the surviving parent. Decreasing
$\mathrm{UB}$ only shrinks the set on the left. Consequently
\[
L_{\mathrm{glob}}=
\min\{\mathrm{UB},\min_{B\in\mathcal Q}\mathrm{LB}(B)\}\le p^\star.
\]
If $\mathcal Q$ is empty, a finite incumbent is optimal, and no incumbent
means the original problem is infeasible. A finite gap at most
$\varepsilon$ certifies $\varepsilon$-optimality. Mid-batch reporting must
include boxes being processed as well as boxes still queued.

\paragraph{Tolerance-based pruning is a different convention.}
Pruning at $\mathrm{LB}(B)\ge\mathrm{UB}-\varepsilon$ also supports an
$\varepsilon$-optimal result, but does not preserve the exact coverage
invariant in \eqref{eq:bab_invariants}. An algorithm using that convention
must retain lower-bound certificates for discarded boxes, or report the
appropriate tolerance certificate rather than the minimum over the
remaining queue alone. \Cref{alg:driver} uses exact dominance and a
separate global-gap stopping rule.

\subsection{Conditional Finite Termination}
\label{app:termination}

Finite termination needs assumptions beyond the sound invariants in
\cref{app:sound-invariants}. The statement below is conditional on local
accuracy, successful primal recovery, detected infeasibility, and an
explicit queue policy. It is not a completeness proof for a fixed-budget
ALM implementation. We use $\operatorname{diam}(B)=\max_j w_j$, the
infinity-norm diameter, throughout this subsection.

\begin{assumption}[Relaxation modulus and value oscillation]
\label{ass:relaxation_modulus}
On the finite, guarded domains under consideration, suppose the original
functions are continuous and there are
$C_{\mathrm{rel}}>0$, $q>0$, and $L_f>0$ such that the current affine
rows obey
\begin{align}
\sup_{x\in B}\{f_i(x)-\underline f_i(x)\}
&\le C_{\mathrm{rel}}\operatorname{diam}(B)^q,
&&i=0,\ldots,m,
\label{eq:relaxation_modulus}\\
\sup_{x,y\in B}|f_0(x)-f_0(y)|
&\le L_f\operatorname{diam}(B).
\label{eq:value_oscillation}
\end{align}
Row validity on $\mathcal F(B)$ is separately maintained as in the main
text. The error inequality above is an additional whole-box accuracy
assumption, not a consequence of feasible-only validity.
\end{assumption}

\paragraph{Local relaxations versus a graph-wide modulus.}
For a $C^2$ univariate operator with bounded second derivative,
secant-tangent error is quadratic in the local interval width. Bilinear
error scales with the product of the two factor widths. For $\sqrt t$
on $[0,\beta]$, the parallel secant-tangent gap is $\sqrt\beta/4$.
These are local statements. A uniform graph-wide exponent and constant
also depend on composition, intermediate interval behavior, and domain
guards; the local formulas alone do not prove a universal $q=2$ for
every supported graph. Assumption~\ref{ass:relaxation_modulus} retains this as an
explicit assumption.

\begin{assumption}[Sound oracle and primal recovery]
\label{ass:oracle_and_primal}
Suppose the following properties hold for the evaluated node procedure.
\begin{enumerate}[leftmargin=1.5em]
\item \emph{Objective backstop.} For a function $e_{\mathrm{clip}}(s)\to0$
as $s\to0$, with $e_{\mathrm{clip}}(0)=0$, its stored bound satisfies
$\mathrm{LB}(B)\ge d_B(0)-e_{\mathrm{clip}}(\operatorname{diam}(B))$.
\item \emph{Local primal recovery.} For every $\varepsilon_P>0$, there is
$\delta_P(\varepsilon_P)>0$ such that processing a box containing a
feasible point with diameter at most $\delta_P$ includes a primal search
returning a feasible $\bar x_B$ with
$f_0(\bar x_B)\le p^\star(B)+\varepsilon_P$.
\item \emph{Small-scale infeasibility detection.} There is $\delta_I>0$
such that every infeasible box with diameter at most $\delta_I$ is
certified by the performed clipping or mixture tests. For a mixture
certificate with exact value $r_B$, it suffices that
$r_B>\tau_{\mathrm{num}}+e_r(\operatorname{diam}(B))$.
\item \emph{Residual accuracy.} The corresponding computed certificate
satisfies $r_B-e_r(\operatorname{diam}(B))\le\widehat r_B\le r_B$, with
$e_r(s)\to0$ and $e_r(0)=0$, and pruning uses $\widehat r_B>\tau_{\mathrm{num}}$.
Here $\tau_{\mathrm{num}}\ge0$ is a numerical safety margin, not the
feasibility tolerance already absorbed into the functions.
\end{enumerate}
\end{assumption}

These assumptions concern actual evaluated bounds and successful tests,
not merely the existence of an LP optimum or an infeasibility vector.
In particular, local primal recovery is not established for ALM-PGD by
this paper. Small-scale infeasibility detection is likewise stronger
than feasible-set preservation. The assumptions make the distinction
between soundness and conditional completeness explicit, in the spirit
of spatial global-optimization analysis
\citep{tawarmalani2002convexification}.

\begin{theorem}[Conditional $\varepsilon$-global convergence]
\label{thm:eps_global}
Assume Assumptions~\ref{ass:relaxation_modulus} and~\ref{ass:oracle_and_primal}, midpoint
splitting with \eqref{eq:hybrid_branching}, and preservation of inherited
lower bounds. Suppose each batch includes a box with the smallest stored
lower bound in the queue. In the absence of a time limit,
\cref{alg:driver} terminates with an $\varepsilon$-optimal incumbent when
the problem is feasible, or certifies infeasibility otherwise.
Choose $0<\varepsilon_P<\varepsilon$ and $\delta>0$ such that
$\delta\le\min\{\delta_P(\varepsilon_P),\delta_I\}$ and
\[
\sup_{0\le s\le\delta}
\{C_{\mathrm{rel}}s^q+L_f s+e_{\mathrm{clip}}(s)\}
+\varepsilon_P\le\varepsilon.
\]
For initial widths $w_j^0$, define
$M_j=\max\{0,\lceil\log_2(w_j^0/(\tau\delta))\rceil\}$ when $w_j^0>0$,
and $M_j=0$ otherwise, with $M=\sum_j M_j$.
The number of batch iterations before certification is bounded by
\begin{equation}
N_{\mathrm{iter}}\le2^{M+1}-1,
\label{eq:eps_global_rate}
\end{equation}
so a batch size at most $K$ processes at most $K(2^{M+1}-1)$ boxes.
\end{theorem}

\begin{proof}
For a feasible box, let $x_0$ minimize its objective row over $B$ and
let $x_F$ be a feasible minimizer. The two accuracy assumptions give
\[
d_B(0)=\underline f_0(x_0)
\ge f_0(x_0)-C_{\mathrm{rel}}\operatorname{diam}(B)^q
\ge p^\star(B)-L_f\operatorname{diam}(B)
-C_{\mathrm{rel}}\operatorname{diam}(B)^q.
\]
Hence processing a feasible box of diameter at most $\delta$ supplies
both a lower bound and an incumbent with
$\mathrm{UB}-\mathrm{LB}(B)\le\varepsilon$.
An infeasible box at that scale is detected by assumption.

Call a queued box \emph{unresolved at tolerance} when no finite incumbent
exists or $\mathrm{LB}(B)<\mathrm{UB}-\varepsilon$.
Once a box is resolved at tolerance, all its descendants remain so:
they inherit its lower bound, and the incumbent can only decrease.
Before the global stopping test succeeds, the best-bound queue policy
therefore processes at least one previously unprocessed unresolved box
in each batch. An unresolved path cannot continue through a processed
box of diameter at most $\delta$.

Every split on such a path has maximum width greater than $\delta$,
so its selected width exceeds $\tau\delta$. After $t$ previous splits
of coordinate $j$, its width is at most $w_j^0/2^t$; clipping can only
reduce it. Thus coordinate $j$ is split at most $M_j$ times on an
unresolved path. The unresolved binary tree has depth at most $M$ and
at most $2^{M+1}-1$ nodes. Each batch processes at least one of these
nodes, proving the iteration bound. Additional resolved boxes processed
in the same batch cannot delay termination beyond this bound. When no
unresolved box remains, the global gap test succeeds or the queue is
empty. The reported certificate follows from the sound invariants.
\end{proof}

\clearpage
\section{Complete Description of the Solver}
\label{appsec:components}

\begin{figure}[h!]
\centering
\includegraphics[width=\linewidth]{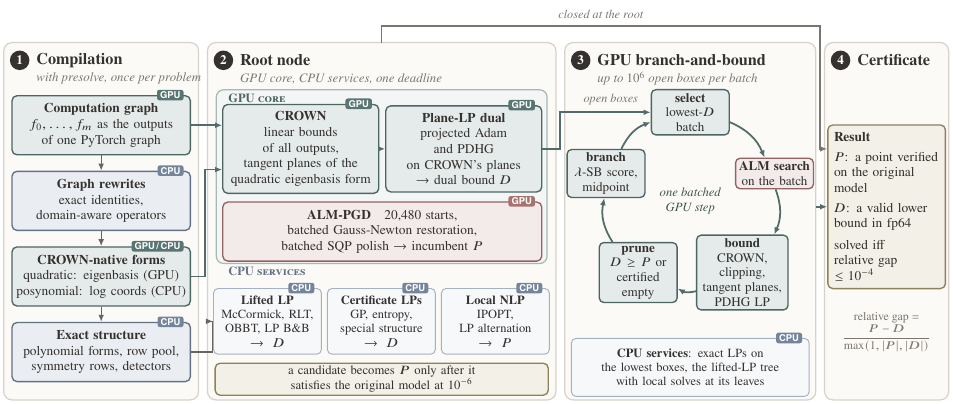}
\caption{\textbf{All components of the solver evaluated in \cref{sec:exp}.}  (1) Compilation and presolve:
the model becomes one computation graph whose outputs are the objective and all constraints, and
recognized structure is rewritten into forms whose tangent planes CROWN computes itself.  (2) At the
root, the GPU core computes CROWN bounds, the plane-LP dual bound $D$ and the batched ALM search for an
incumbent $P$; CPU services add bounds and points.  A candidate becomes $P$ only after it satisfies the original model at $10^{-6}$, and
the run stops once the relative gap of \cref{eq:exp_solver_criterion} is at most $10^{-4}$.  (3) Open boxes enter the batched GPU branch-and-bound loop.
(4) The result is certified by the gap rule of \cref{sec:exp}.  The tag on each box shows whether it
runs on the GPU or the CPU.}
\label{fig:components}
\end{figure}

This section lists every technique used by the solver evaluated in \cref{sec:exp}, in the order in which
a run meets them (\cref{fig:components}).  The GPU core of \cref{sec:method} is summarized only where it
interacts with the other components.  Every technique either adds a valid bound, removes a region that
provably contains no better feasible point, or proposes a point that is accepted only after the
feasibility check of \cref{eq:feasibility_gate}; the soundness argument of \cref{app:sound-invariants} therefore
covers all of them, for valid relaxations and one consistent tolerance-relaxed feasible set.  All statements below
refer to that set, to the stated box and cutoff scope, and to the stated numerical enclosures: optimization routines
only propose points or multipliers, and algebraic identities with independently evaluated enclosures decide validity.

\subsection{Compilation and Presolve}
\label{appsec:comp_compile}

\paragraph{Computation graph.}  The AMPL model is translated into one PyTorch graph whose outputs are the
objective and all constraints (\cref{appsec:dataset}); linear rows keep their exact coefficients and are also
appended as outputs, so that every bound below sees every row.  Tracing and the rewrites of this subsection are
Python passes on the CPU, run once per graph; the graph is then evaluated and bounded on the GPU, batched over
boxes.

\begin{figure}[h]
\centering
\includegraphics[width=\linewidth]{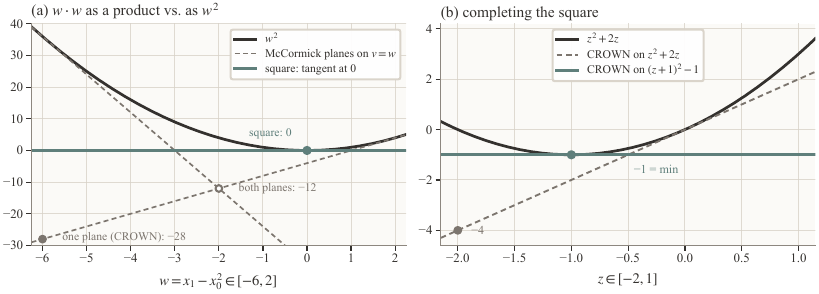}
\caption{Exact rewrites that tighten CROWN.  (a) $w=x_1-x_0^2\in[-6,2]$ on $[-2,2]^2$: restricted to $v=w$,
the McCormick lower planes of $v\cdot w$ are the endpoint tangents of $w^2$; CROWN keeps the one with the larger
box minimum and certifies $-28$ ($-2800$ for the Rosenbrock term $100w^2$), both planes would give $-12$, and the
square's tangent at $\mathrm{proj}(0)$ gives the exact minimum $0$.  (b) On $[-2,1]$, CROWN bounds $z^2+2z$ by
$2z$, i.e.\ $-4$ (as interval arithmetic); after completing the square, the tangent of $(z+1)^2$ at $0$ gives the
exact minimum $-1$.}
\label{fig:appb_exact_rewrites}
\end{figure}

\paragraph{Exact graph rewrites.}  Before bounding, the graph is normalized by identities that hold exactly, so
that CROWN's local relaxations see the structure (\cref{fig:appb_exact_rewrites}).  Identical gathers and common
subexpressions are merged.  A product $u\cdot u$, $u\in[\ell,\upsilon]$, becomes $\mathrm{Sqr}(u)$, whose
tangent at $\mathrm{proj}_{[\ell,\upsilon]}(0)$ attains the exact minimum; the two McCormick lower planes of $u\cdot v$
restricted to $v=u$ are the tangents at $\ell$ and $\upsilon$ (the upper planes restrict to the secant), and CROWN
keeps only one lower plane.  A product $u\cdot v$ with $v=cu$ proved
elementwise in rational polynomial arithmetic becomes $c\,\mathrm{Sqr}(u)$, and a nearly proportional pair is
split as $uv=cv^2+v(u-cv)$ with the residual kept exactly.  For $c\ne0$, a quadratic $cz^2+kz$ in an input or shared node $z$
becomes $c(z+s)^2-cs^2+rz$ with $s=\mathrm{fl}(k/2c)$ and the exact rational remainder $r=k-2cs$, which centres the
square near the quadratic's vertex (exactly at the vertex when $r=0$).  A variable exponent becomes $x^y=\exp(y\log x)$, exact for $x>0$.

\paragraph{Domains and poles.}  These guards act elementwise in every bounding pass on the GPU.  An operator with
a restricted domain ($\log$, $\sqrt{\cdot}$, non-integer powers) intersects its input interval with the domain,
$[\ell,\upsilon]\leftarrow[\max(\ell,0),\upsilon]$, since no point outside it belongs to the model; an endpoint at
$0$ stands for an excluded limiting boundary (for $\log$ and negative powers the domain is $x>0$), poles and undefined
sides are handled separately, and the resulting bounds apply only where the original graph is defined.  A power $x^p$
with constant $p<0$ is convex and decreasing on $(0,\infty)$: it is bounded by a tangent below and the secant
above, and on $[0,\upsilon]$ by the tangent alone (at a point of $[\upsilon/2,\upsilon]$), because the secant
would pass through the pole.  Where an element $e$ has no relaxation on one side (a pole or an infinite endpoint in
its interval), each backward coefficient $\Lambda_{re}$ that selects that side is replaced by the element's
interval contribution, $\Lambda_{re}z_e\ge\min(\Lambda_{re}\underline z_e,\Lambda_{re}\overline z_e)$; only rows
whose contribution is not finite get the always valid bound $-\infty$, and all other rows propagate unchanged.  A
reciprocal across a pole, $0\in(\ell,\upsilon)$, keeps the excluded gap $1/x\notin(1/\ell,1/\upsilon)$ through
shifts and scalings, so that, e.g., $(1/x)^2\ge\min(\ell^{-2},\upsilon^{-2})>0$.  Every intermediate interval is the
intersection of its CROWN and interval bounds.

\paragraph{Exact structure and the row pool.}  A symbolic walk on the CPU attempts to compute, for every supported node, an exact sparse
polynomial with rational coefficients in the inputs and in transcendental atoms ($\exp$, $\log$, powers, $\sin$,
$\cos$), with a budget counted in term operations rather than time; nodes with unsupported operators, exhausted
budget or errors get an unknown form and are excluded from the exact structural proofs.  It feeds the
certificates, the lifted LP, the exact sparse derivatives used by IPOPT and the pattern propagation of
\cref{appsec:comp_bab}.  All rows enter one pool as $\ell_k\le q_k(x)\le h_k$, with $q_k$ the non-constant part
and the feasibility tolerance $\tau$ added once; rows $q_j(x)\le\beta_j$ with $q_j=s_jq$ for one direction $q$ merge
exactly into $\max_{s_j<0}\beta_j/s_j\le q\le\min_{s_j>0}\beta_j/s_j$ (rounded outward once), and rows implied by
the root box are masked.  Two global facts are proved on the exact forms.  A constraint row that equals
$\alpha f_0+\beta$ with $\alpha<0$ and threshold $t$ gives the floor $f_0\ge(t+\tau-\beta)/\alpha$, applied to every
box.  A row $a\,h(x)^2+\beta\le t$ with $a>0$ gives $|h|\le r=\sqrt{(t+\tau-\beta)/a}$ (rounded up), and the CROWN
lower planes $p_\pm(x)\le\pm h(x)$ on the root box turn it into two linear rows $p_\pm(x)\le r$, global on the root box; for nonlinear $h$ each batch recomputes
stronger planes on its own boxes, and these rows are local to those boxes and their descendants.  A negative radicand
certifies that the square row has no feasible point.

\paragraph{Symmetry.}  On the exact forms (CPU), a chain of disjoint variable blocks $B_1,\dots,B_m$ with equal root
boxes is accepted when every exchange of two consecutive blocks maps the objective to itself and the multiset of rows
to itself, compared exactly.  Adjacent exchanges generate all permutations of the blocks, so every feasible point can
be permuted, at the same objective value, into $x_{B_1[1]}\le\dots\le x_{B_m[1]}$ (the first coordinates of the
blocks).  These $m-1$ ordering rows are added as global linear rows, which constrained concretization uses at every
node on the GPU, and incumbents are sorted into the ordered region.  Chains used together are disjoint.  This removes symmetric copies of the optimum
without excluding all of them; the transformation preserves the real objective and feasible set exactly, although a
floating-point evaluation order may differ for the permuted point, which is not re-gated when an incumbent is sorted.

\paragraph{Rewrites for native CROWN bounds.}  Two proved structures are rewritten into small graphs whose
convexity CROWN's own relaxations exploit.  A quadratic objective (\cref{appsec:comp_dual}) is written at an anchor
$\hat x$ in the eigenbasis of its Hessian on the free coordinates, $H=V\operatorname{diag}(\lambda)V^\top+E$ with
$E$ the residual of the decomposition:
\begin{equation}
f_0(x)=f_0(\hat x)+g^\top d+\tfrac12\textstyle\sum_k\lambda_kz_k^2+\tfrac12d^\top Ed,
\qquad d=x-\hat x,\quad z=V^\top d,
\label{eq:appb_eigen_rewrite}
\end{equation}
an exact identity whose last term is at least $-\tfrac12\|E\|_2\max_{x\in B}\|d\|^2$ on a box $B$.  The graph
$\mathrm{Linear}\to\mathrm{Sqr}\to\mathrm{Linear}$ is built on the GPU (Hessian by automatic differentiation,
symmetric eigendecomposition) and rebuilt when the anchor moves.  Each posynomial $p(x)=\sum_kc_kx^{a_k}$
($c_k>0$, real exponent vectors $a_k$) is written in $y=\log x$ as $\sum_k\exp(a_k^\top y+\log c_k)$, a sum of
convex terms (\cref{fig:appb_gp_logcoords}); this graph is built and bounded on the CPU.

\subsection{Root Schedule}
\label{appsec:comp_root}
The root steps below run in a fixed order under the single deadline of the run.  After every step the gap
rule \eqref{eq:exp_solver_criterion} is checked against the best verified point and the best bound, and the
run ends as soon as it holds; otherwise the remaining open boxes start the GPU branch-and-bound of
\cref{sec:method}.  In the run of \cref{sec:exp}, 293 of the 421 solved instances close at the root and 128 in
the tree.  By the first bound that reaches $P-\varepsilon$ on each instance, the lifted LP and its
branch-and-bound close 214, CROWN 35, the separable floors 20, the posynomial and relative-entropy
certificates 16, the special-structure certificates 7 and another root bound 1; in the GPU tree, 94 close from
the root box and 33 from the lifted-LP leaves, and the infeasibility of \texttt{mosel\_start\_folioqpgraph} is
proved in the tree.  The first bound to close an instance need not be the only one that would.  The primal searches of
\cref{appsec:comp_primal} overlap these steps in parallel root mode.

\subsection{Root Dual Bounds and Certificates}
\label{appsec:comp_dual}

\paragraph{CROWN and separable floors.}  CROWN, without $\alpha$ optimization, gives the planes
\eqref{eq:method_linear_bounds} of all outputs over the root box on the GPU; these planes are also the rows of the
dual LP.  When the objective decomposes along the graph as $f_0=\kappa+\sum_jc_jt_j(x)$ with at most 64 terms, each
term is bounded separately,
\[
f_0\ \ge\ \kappa+\sum_{c_j>0}c_j\max\bigl(\underline t_j^{\mathrm{IBP}},\underline t_j^{\mathrm{CROWN}}\bigr)
+\sum_{c_j<0}c_j\min\bigl(\overline t_j^{\mathrm{IBP}},\overline t_j^{\mathrm{CROWN}}\bigr)
\quad\text{on the root box,}
\]
and this value is a floor for every box.  The terms are bounded on the GPU one after another and summed on the host.

\begin{figure}[h]
\centering
\includegraphics[width=\linewidth]{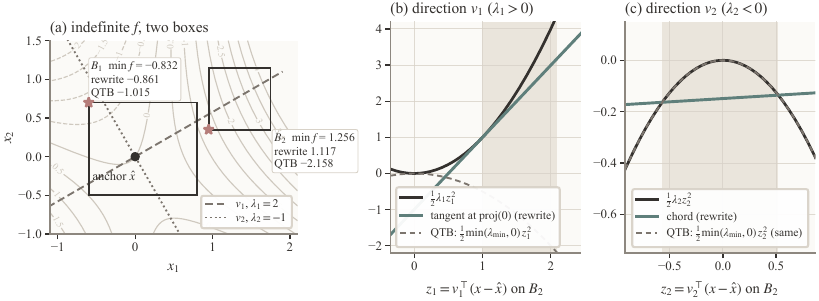}
\caption{QTB vs.\ CROWN on the eigenbasis rewrite, for $f(x)=g^\top x+\frac12x^\top Hx$ with eigenvalues $(2,-1)$,
anchor $\hat x=0$ and $g=(0.4,-0.3)$.  (a) Contours, eigen-directions and two boxes ($\star$: box minimizers).
On $B_1\ni\hat x$ the rewrite certifies $-0.861$ and QTB $-1.015$ against the minimum $-0.832$; on $B_2$, away from
the anchor, $1.117$ and $-2.158$ against $1.256$.  (b,c) On $B_2$, the rewrite bounds $\frac12\lambda_kz_k^2$ on
the interval of $z_k$ (shaded) by the tangent at $\mathrm{proj}(0)$ ($\lambda_1>0$) or by the chord ($\lambda_2<0$);
QTB replaces every $\lambda_k$ by $\min(\lambda_{\min},0)$ and then $\|x-\hat x\|^2$ by its maximum over the box.}
\label{fig:appb_qtb_eigenbasis}
\end{figure}

\paragraph{Quadratic tangent bound.}  A degree propagation on the graph proves on the CPU a sufficient bound
$\widehat{\deg}f_0\le2$ ($\deg x_i=1$, $\deg(uv)=\deg u+\deg v$, and $\infty$ for non-polynomial operators; failure
does not imply that $f_0$ is not quadratic).  If it succeeds and
$n\le2{,}000$, then $f_0(x)=f_0(\hat x)+g^\top d+\tfrac12d^\top Hd$ exactly for any anchor $\hat x$, with $g=\nabla f_0(\hat x)$ and a
constant Hessian $H$.  $H$ and an estimate $\lambda_{\mathrm{lb}}=\lambda_{\min}(H)-10^{-9}(1+\|H\|_F)$ of a lower bound on its
smallest eigenvalue are computed on the GPU in floating point; the cushion is a numerical safety allowance, not a
verified enclosure of the true Hessian, so the bounds below are rigorous only if it holds.  The anchor $\hat x$ is an
approximate KKT point of the root problem with every row relaxed by the tolerance (accepted from IPOPT's exit status and
a finite feasibility tolerance, without a KKT proof), found by IPOPT on the CPU (at most 10\,s), whose
solve also proposes incumbents.  The anchor then moves to each improved incumbent, except that an anchor of a convex
objective with affine rows, whose plane then bounds the feasible set by $f_0(\hat x)$ when it is an exact minimizer,
is kept unless an incumbent improves on $f_0(\hat x)$.  Any anchor gives a valid plane when the curvature and data
errors are enclosed, but tightness at the anchor is not guaranteed.  The tangent plane alone gives the quadratic tangent bound (QTB)
\[
f_0(x)\ \ge\ f_0(\hat x)+g^\top d+\tfrac12\min(\lambda_{\mathrm{lb}},0)
\textstyle\sum_j\max\bigl\{(\hat x_j-\underline x_j)^2,(\overline x_j-\hat x_j)^2\bigr\}\qquad(x\in B),
\]
which lowers every direction by the most negative curvature.  The solver instead runs CROWN on the rewrite
\eqref{eq:appb_eigen_rewrite}, per box on the GPU: with $[\underline z_k,\overline z_k]$ the interval of $z_k$ on
$B$, the term $\tfrac12\lambda_kz_k^2$ is bounded by the tangent at
$\mathrm{proj}_{[\underline z_k,\overline z_k]}(0)$ if $\lambda_k>0$ and by the chord if $\lambda_k<0$
(\cref{fig:appb_qtb_eigenbasis}).  On a box that holds $\hat x$ the positive directions reproduce QTB's plane, and
elsewhere they can be tighter, when the interval of $z_k$ excludes $0$; a negative direction is bounded by its chord, exact at both ends of its interval, instead
of by the most negative curvature at the farthest corner of the box.  The
resulting plane, lowered by $\tfrac12\|E\|_2\max_{x\in B}\|d\|^2$ (with $\|E\|_2$ computed in floating point) and an
empirical rounding cushion, replaces the CROWN objective
plane in a second dual LP over the same rows (\cref{appsec:comp_bab}).  When $\lambda_{\mathrm{lb}}>0$ the curvature
is also kept: for the multipliers $y\ge0$ of that LP (and for $y=0$),
$f_0(x)+y^\top(Gx-h)\ge f_0(\hat x)+g^\top d+\tfrac{\lambda_{\mathrm{lb}}}2\|d\|^2+y^\top(Gx-h)$ is separable and
minimized at $x_j^\star=\mathrm{proj}_{[\underline x_j,\overline x_j]}(\hat x_j-q_j/\lambda_{\mathrm{lb}})$ with
$q=g+G^\top y$; the bound is evaluated rigorously after linearizing
$\tfrac\lambda2t^2\ge st-\tfrac{s^2}{2\lambda}$ at $s=\lambda(x_j^\star-\hat x_j)$.

\paragraph{Convex objectives and rows.}  A non-quadratic objective proved convex on the root box (composition rules
on the CPU, each univariate second derivative enclosed in interval arithmetic on 64 uniform slices of its range, augmented with
dyadic slices at both ends) contributes
its tangent plane $f_0(\hat x)+\nabla f_0(\hat x)^\top(x-\hat x)$, lowered by the enclosure errors, to the same second
LP.  Constraint rows $g_i(x)\le t_i$ that are quadratic and admitted as numerically near-convex ($\lambda_{\min}$ estimate
$\ge-10^{-8}$, the negative-curvature term below being kept), or convex by the same rules, give at the incumbent $\hat x$ the global linear rows
$g_i(\hat x)+\nabla g_i(\hat x)^\top(x-\hat x)+\tfrac12\min(\lambda_{\mathrm{lb},i},0)\max_{x\in B_0}\|x-\hat x\|^2\le
t_i+\tau$ (the curvature term only for quadratic rows, $B_0$ the root box, minus a rounding cushion).  They are
refreshed when the incumbent improves, for at most four rounds, and only when there are at most 50 candidate quadratic
rows; the Hessians are computed on the GPU, one row at a time.

\begin{figure}[h]
\centering
\includegraphics[width=\linewidth]{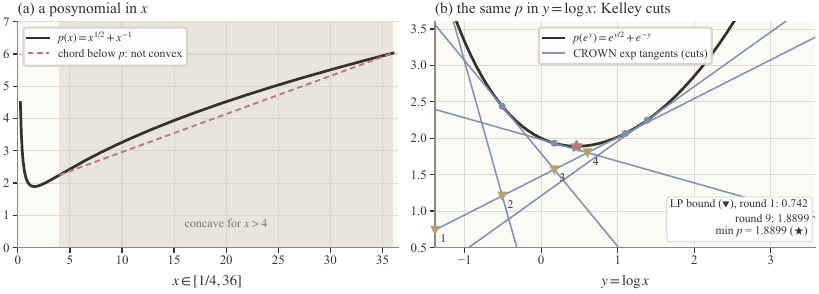}
\caption{Posynomial certificate for $p(x)=x^{1/2}+x^{-1}$ on $[1/4,36]$ (CPU).  (a) $p$ is concave for $x>4$, so
it is not convex in $x$.  (b) In $y=\log x$, $p(e^y)=e^{y/2}+e^{-y}$ is convex and CROWN's tangents of the
exponential terms are global lower lines.  The Kelley LP over these cuts (points: box centre, incumbent $x=4$, then
each LP point; $\blacktriangledown$: LP minimizers of rounds 1-4) certifies $0.742$ in round 1 and the minimum
$1.8899$ to four digits in round 9.}
\label{fig:appb_gp_logcoords}
\end{figure}

\paragraph{Posynomial (geometric-program) certificate.}  When the root box is strictly positive and the exact forms
show the objective to be a posynomial plus a constant, the substitution $x=e^y$ makes the objective and every
posynomial row convex in $y$; other rows are dropped, which only weakens the bound.  A positive monomial equality with tolerance is equivalent to
two log-affine inequalities when both endpoints are positive; in the enabled CROWN certificate path both sides are
instead written as positive monomial or reciprocal rows and relaxed by exponential tangents.  A Kelley cutting-plane LP
in the epigraph of the objective, whose constant offset $\kappa$ is added back,
\[
\min\ \kappa+t\quad\text{s.t.}\quad t\ge\pi_0^{(j)}(y),\quad \pi_i^{(j)}(y)\le t_i+\tau,\quad
\log\underline x\le y\le\log\overline x,
\]
is solved by an LP solver on the CPU, where each cut $\pi^{(j)}(y)=\sum_ke^{m_k}(1+a_k^\top y+\log c_k-m_k)\le p(e^y)$ is CROWN's
tangent of the exponential terms of the log-coordinate rewrite at $m_k=a_k^\top y^{(j)}+\log c_k$ (clamped to the
term's range), taken at the box centre, the incumbent and then each LP point (\cref{fig:appb_gp_logcoords}).  CROWN
only proposes the slopes: the intercepts and the LP's dual value are re-evaluated in 50-digit interval arithmetic,
so the reported bound is rigorous for the cut data.  The loop stops when a round after the fifth no longer improves the
bound, or after 20\,s.  The certificate, CROWN included, runs on the CPU.

\paragraph{Relative-entropy certificate.}  When the variables enter through $y_j=e^{x_j}$, the objective is a sum
of relative-entropy terms $w\,u\log(u/v)$ plus an affine part, with weights $w\ge0$ and $u,v$ proved strictly
positive on the transformed box (then jointly convex), and the rows are affine in $y$.  The implementation recognizes
numerators $y_j$ or $1$ and affine denominators, retains usable affine rows and drops unsupported ones.  Each term gets its own epigraph variable with outward-rounded tangent cuts, and the LP dual value
of the LP solver's multipliers, evaluated in interval arithmetic, is the bound (CPU).

\paragraph{Special-structure certificates.}  Three detectors match exact algebraic patterns of the forms and rows
and prove a floor with a dedicated argument; numerical solvers only propose.  (i) A polygon area in polar
coordinates whose law-of-cosines rows bound every pairwise distance by $d$: the convex hull of the $N$ vertices has
diameter at most $d$, hence perimeter at most $2Nd\sin\frac{\pi}{2N}$ and area at most
$Nd^2\sin^2\!\frac{\pi}{2N}/\tan\frac{\pi}{N}$, evaluated in interval arithmetic on the CPU; for the
tolerance-relaxed polar model the floor additionally carries the allowance for the order rows' backward angles
(with $d$ the relaxed distance threshold).  (ii) A normalized symmetric eigenvalue block with an affine objective
in the eigenvalue: for any computed factor $L$ and $\rho\ge\|A-sI-LL^\top\|_\infty$ (induced infinity norm, the maximum absolute row
sum; rounded outward), $A-(s-\rho)I\succeq LL^\top\succeq0$, so $\lambda_{\min}(A)\ge s-\rho$.  Because the
eigenpair rows hold only within a tolerance, with $\|q\|^2\ge L_q>0$, $|(Aq-\lambda q)_i|\le\epsilon_i$ and
$f_0=a\lambda+b$, $a>0$, the floor is $b+a\bigl(s-\rho-\sqrt{\sum_i\epsilon_i^2/L_q}\bigr)$; the eigendecomposition, the Cholesky factor
and the residual run on the GPU.  (iii) Near-binary assignment rows: for $\tau<\frac14$, $|x(x-1)|\le\tau$ puts $x$ within
$\delta=2\tau/(1+\sqrt{1-4\tau})$ of $\{0,1\}$; after proving $0\le x\le1$, the assignment structure and
$m\delta+\epsilon<\frac12$ for every cardinality row, every relaxed-feasible point rounds to a partition.  Each
enumerated partition receives an outward-rounded lower bound over its allowed near-binary deviations, and a cover
search minimizes the sum of these bounds in GPU kernels (prices from an LP solver on the CPU).  Each
pattern matches a single family of the benchmark (\texttt{polygon}, \texttt{eigmin}, \texttt{grouping}).

\subsection{Lifted Linear Relaxation}
\label{appsec:comp_lifted}

The plane LP~\eqref{eq:node_lp} relaxes every row by its own plane; the lifted LP keeps each nonlinear term as
one column shared by all rows.  It is built once at the root.  Its LP solves (bound tightening, the lifted root
bound and the LP branch-and-bound) are floating-point LP solves on the CPU, serial and synchronous between GPU
iterations, whose multipliers are then evaluated by the certified weak-duality bound below; the GPU batches bound the
boxes admitted by the dense budget over the same rows with PDHG, and the other boxes keep their other bounds.

\paragraph{Lifting.}  Problems with at most 2{,}000 variables and a finite root box are lifted.  A row $f_i$
($i=0$ is the objective) whose exact form (\cref{appsec:comp_compile}) is a polynomial in the inputs and in
transcendental atoms (at most 64 in the problem) becomes linear in the columns:
\begin{equation}
f_i(x)=s_i\bigl(z(x)\bigr),\qquad s_i(z)=k_i+g_i^\top z,\qquad z=(x,\,t,\,w).
\label{eq:appb_lifted_row}
\end{equation}
There is one column $t_a=\phi_a(\alpha_a^\top z+\alpha_{a0})$ per atom of an affine argument
($\phi\in\{\exp,\log,x^p,\sin,\cos\}$, $1/x$ included) and one column $w_k=z_{a_k}z_{b_k}$ per distinct product; a monomial of higher degree is split into halves ($x^4=(x^2)(x^2)$), so both
factors of a product are earlier columns.  Quadratic rows take their coefficients from automatic differentiation at
$x=0$ on the GPU, with a numerical sanity check at four random points (not a coefficient-error enclosure); the other rows take them from the exact forms on the CPU.  On a
box $B$, interval arithmetic over the columns gives a column box $Z(B)=[l^z,u^z]$.  With the relaxation and RLT rows
$E_Bz\le e_B$ below, the base linear rows included, and an incumbent value $P$, define
\begin{equation}
\mathcal L(B,P)=\bigl\{z\in Z(B):\ s_i(z)\le0\ (i\ge1),\ E_Bz\le e_B,\ s_0(z)\le P\bigr\}.
\label{eq:appb_lifted_lp}
\end{equation}
In exact arithmetic every $x\in B$ that satisfies the rows with $f_0(x)\le P$ lifts to $z(x)\in\mathcal L(B,P)$, so
$\min\bigl(P,\min_{\mathcal L(B,P)}s_0\bigr)$ bounds the problem on $B$.  The first-level product endpoints of $Z(B)$
are computed in floating point without outward rounding (deeper levels are widened), so this containment holds up to
rounding, which the rows' cushion absorbs only numerically: each right-hand side is relaxed by
$10^{-12}\bigl(1+|k_i|+\sum_j|g_{ij}|\max(|l^z_j|,|u^z_j|)\bigr)$.  The implementation invokes the standalone lifted
root and tree bound only when the objective and every row with a finite threshold are lifted (a partial relaxation
would also be valid, if weaker); the other cases use it for bound tightening and for the per-box bound (with CROWN
planes for the other rows).  A convex quadratic objective over linear rows is
left to the quadratic tangent bound (\cref{appsec:comp_dual}), which is exact there.

\begin{figure}[h]
\centering
\includegraphics[width=\linewidth]{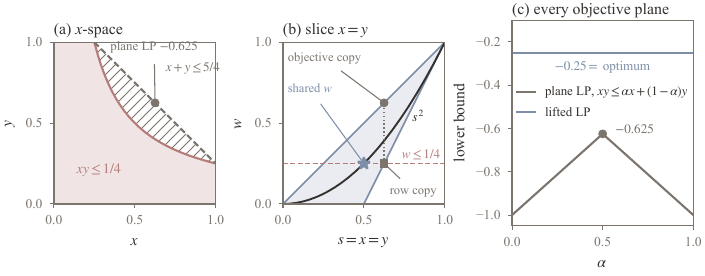}
\caption{Per-row planes versus one shared column on $\min -xy$ subject to $xy\le\frac14$ over $[0,1]^2$ (optimum
$-\frac14$).  (a)~The plane LP relaxes the objective to $-(\alpha x+(1-\alpha)y)$ (from $xy\le x$, $xy\le y$) and
the row to $x+y-1\le\frac14$; it admits the hatched points.  (b)~On the diagonal $x=y=s$, its optimum gives the same
product two values: $\frac58$ in the objective's copy and $\frac14$ in the row's.  Full McCormick envelopes with one
copy per row still give $-0.625$; one shared column $w$ gives $-0.25$.  (c)~Over all objective planes the plane LP
reaches at most $-0.625$ (at $\alpha=\frac12$); the lifted LP reaches $-0.25$, the optimum.  All values are LP
optima.}
\label{fig:appb_lifting_vs_planes}
\end{figure}

\paragraph{Why lifting helps.}  Let $R_i$ hold the points $(x,v,w)$ that satisfy the relaxation of row $i$ (for the
objective, $v\ge s_0$) and the relaxations of the columns it uses, and let $\pi$ project onto $(x,v)$.  A row relaxed
by itself in $x$-space keeps a private copy of $w$, so the best per-row relaxation minimizes $v$ over
$\bigcap_i\pi(R_i)$; the lifted LP shares $w$ and minimizes over
$\pi\bigl(\bigcap_iR_i\bigr)\subseteq\bigcap_i\pi(R_i)$.  The inclusion can be strict
(\cref{fig:appb_lifting_vs_planes}).

\begin{figure}[h]
\centering
\includegraphics[width=\linewidth]{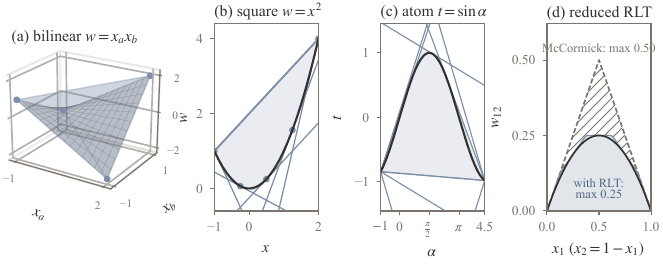}
\caption{Relaxation rows on one box.  (a)~The four McCormick faces of $w=x_ax_b$, the convex hull of the saddle over
the box.  (b)~$w=x^2$: the secant above, and tangents at $l+\{0,\frac14,\frac12,\frac34,1\}(u-l)$ below.
(c)~$t=\sin\alpha$ on $[-1,4.5]$, across two inflection points: support lines
$\sigma(t-m\alpha)\le H_\sigma(m)$ with tangent slopes at $\frac14,\frac12,\frac34$ of the range and the tightest
slopes through both ends.  (d)~Reduced RLT for $x_1+x_2=1$ on $[0,1]^2$ with the product $x_1x_2$.  Multiplying by
$x_1$ and by $x_2$ adds $w_{11}+w_{12}=x_1$ and $w_{12}+w_{22}=x_2$.  The projection on $(x_1,w_{12})$ shrinks from
the McCormick triangle (hatched, $\max x_1x_2=0.50$) to a band over $x_1(1-x_1)$ ($\max 0.25$, exact).}
\label{fig:appb_lifted_envelopes}
\end{figure}

\paragraph{Operator relaxations.}  A product $w=x_ax_b$ gets the four McCormick faces of
$(x_a-l_a)(x_b-l_b)\ge0$, $(u_a-x_a)(u_b-x_b)\ge0$, $(u_a-x_a)(x_b-l_b)\ge0$ and $(x_a-l_a)(u_b-x_b)\ge0$; for
example, $w\ge l_bx_a+l_ax_b-l_al_b$.  A square gets its secant $w\le(l+u)x-lu$ and the tangents
$w\ge2\theta x-\theta^2$ at $\theta\in l+\{0,\frac14,\frac12,\frac34,1\}(u-l)$.  An atom of fixed curvature on its
argument range $[\underline\alpha,\overline\alpha]$ ($\exp$; $\log$ on a positive range, $x^p$ on a one-signed range of its domain) gets five
tangents on its convex side and the secant on the other.  $\sin$ (and $\cos\alpha=\sin(\alpha+\frac\pi2)$) gets
support lines on both sides $\sigma=\pm1$:
\begin{equation}
\sigma\,(t-m\alpha)\le H_\sigma(m)=\max_{\alpha\in[\underline\alpha,\overline\alpha]}\sigma(\sin\alpha-m\alpha).
\end{equation}
These hold on a range of any width.  For $|m|\le1$, $H_\sigma(m)$ is the largest of $\sigma(\sin\alpha-m\alpha)$ at the
two ends and at the first and last $\alpha=\sigma\arccos m+2\pi k$ in range; for $|m|>1$ the ends suffice.  The slopes $m$ are $\cos$ at $\frac14,\frac12,\frac34$ of the range and
the tightest slopes through either end (by bisection).  At the root, two LPs per trigonometric atom certify its
argument range, and 16 cutting-plane rounds add the slope at the LP point.  Floating-point atom values are widened
by $10^{-14}$ relative, a numerical allowance that assumes libm accuracy, including for large-argument phase reduction.  The relaxations are rebuilt on every box: in dense batches on the GPU for PDHG, and sparse on
the CPU for the LP solver (\cref{fig:appb_lifted_envelopes}).

\paragraph{Cleared denominators and convex objectives.}  A row $p(x)/m(x)\le\tau$ whose denominator $m$ is a monomial
of strictly positive root variables adds the cut $m(x)\bigl(p(x)/m(x)-\tau\bigr)\le0$ when this polynomial has degree
at most two.  The cut uses exact rational coefficients with an outward conversion error, and it is dropped when a
later threshold is looser.  A non-quadratic objective proved to be $q_0+\sum_kc_k\varphi_k(q_k)$ ($c_k>0$,
$\varphi_k$ convex, $q_k$ quadratic) gets one epigraph column $e_k\ge\varphi_k(q_k)$ per term, with linear cuts
$e_k\ge\varphi_k(\theta)+\varphi_k'(\theta)\bigl(s_{q_k}(z)-\theta\bigr)$, $s_{q_k}$ being $q_k$ as an affine function of
the columns (convexity of $\varphi_k\circ q_k$ in $x$ is not needed).  Its tangents are
taken at 12 grid points and up to 32 retained points, and on the CPU each LP adds tangents at its own solution (at
most 8 rounds).  Both structures are recognized once on the CPU.

\paragraph{Reduced RLT.}  A pair of opposite linear rows $g^\top x-b\in[-\varepsilon^-,\varepsilon^+]$ (an equality
with its tolerance) is multiplied by $x_i$ when $x_ix_j$ is a product column for some $j$ in the support of $g$.
Since $x_i(g^\top x-b)=\sum_jg_jw_{ij}-bx_i$, this gives $-\varepsilon^-x_i\le\sum_jg_jw_{ij}-bx_i\le\varepsilon^+x_i$
when $l_i\ge0$, the reversed pair $\varepsilon^+x_i\le\sum_jg_jw_{ij}-bx_i\le-\varepsilon^-x_i$ when $u_i\le0$, and
the interval product $[l_i,u_i]\cdot[-\varepsilon^-,\varepsilon^+]$ otherwise (when $x_i$ crosses zero), with
$\varepsilon^\pm$ taken from the outward-relaxed source rows.  Missing
products $x_ix_j$ become new columns with their own relaxations.  The rows are added after the McCormick warm-up below,
all of them or none: at most 1{,}200 RLT rows and 7{,}000 columns.  The selection runs once on the CPU, and the RLT
block of each box is assembled on the host, also for the GPU batches.

\begin{figure}[h]
\centering
\includegraphics[width=\linewidth]{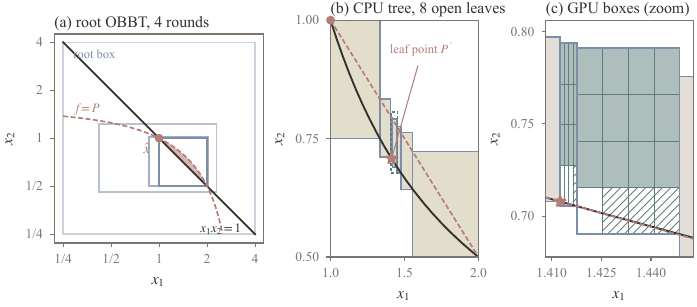}
\caption{From the root box to the GPU tree on $\min x_1+2x_2$ subject to $x_1x_2\ge1$ over $[\frac14,4]^2$, with
incumbent $\hat x=(1,1)$, $P=3$ (optimum $2\sqrt2\approx2.828427$).  (a)~Root OBBT with the cutoff row
$x_1+2x_2\le P$ (log axes).  Four rounds shrink the box to approximately $[1,2]\times[\frac12,1]$ at the displayed precision, the bounding box of
$\{x_1x_2\ge1,\,f\le P\}$ (shaded), and raise the node bound (with propagation) from 2.588 to 2.667.  The toy LP rows
carry no modelling slack; candidate acceptance uses $10^{-6}$.  Without the cutoff row the
box does not move.  (b)~The lifted-LP tree on this box, stopped at 8 open leaves for display (the solver stops at
$\max(256,32n)$); node propagation removed the white parts.  A leaf's LP point verifies as $P'=2.828429$ (star).
(c)~Against $P'$, 6 leaves fall to their own bounds in the GPU tree.  The other two (dashed window in (b)) are cut
into 16 boxes each.  Of these 32 boxes, the leaf's certificate $\kappa+r^\top z$ alone prunes 24 (filled), fresh
lifted LPs prune 4 more (hatched), and 4 stay open.  Computed with a small re-implementation of the same rows.}
\label{fig:appb_lifted_root_to_gpu}
\end{figure}

\paragraph{Bound tightening against the cutoff (OBBT).}  Two minimization LPs per column give certified lower bounds $D_j^-\le\min_{\mathcal L(B,P)}z_j$ and
$D_j^+\le\min_{\mathcal L(B,P)}(-z_j)$, and
\begin{equation}
l_j\leftarrow\max\bigl(l_j,\,D_j^-\bigr),\qquad u_j\leftarrow\min\bigl(u_j,\,-D_j^+\bigr).
\end{equation}
The new box keeps every point that satisfies the rows with $f_0\le P$, the incumbent included (an audit checks it).
The targets are the variables of the products and, for the auxiliary bounds below, the atom and product columns,
admitted by their measured shrink per LP.  The relaxations are rebuilt before each target.  At the root, a McCormick
warm-up (at most 5\,s), RLT and a second pass (10\,s in all) run at most five rounds each, until the mean relative
shrink falls below 1\%.  Later passes run after an incumbent improvement of more than $10^{-4}$ (relative), or while
a pass cut short by its budget was still shrinking (5\,s each, 20\,s in all, until a pass shrinks the mean width by
less than 2\%).  The LPs are serial solves on the CPU; every GPU box is intersected with the latest box,
elementwise on the GPU (\cref{fig:appb_lifted_root_to_gpu}a).

\paragraph{Lifted root bound.}  Once an incumbent exists, the unified root (\cref{appsec:comp_root}) runs the node
solve below on the current box, after the LP is built and again after the warm-up and RLT.  The run ends if this
bound closes the gap.

\paragraph{LP branch-and-bound.}  The tree runs on the CPU, best first and one node at a time, and carries the cutoff
row $s_0\le P$.  At a node $B'$, up to three passes of outward-rounded interval propagation over the rows of
$\mathcal L(B',P)$ tighten the box, and the relaxations are rebuilt after each pass.  Dual simplex then solves
the LP from the parent's basis (at most $5(m+N)$ iterations and 10\,s).  OBBT probes run on the 32 variables with the
largest violation times width when the LP has at most 1{,}200 rows; probing stops when the first ten probed nodes
show that fewer than 20\% of the probes pay.  A node keeps $\beta(B')=\max(\beta(\text{parent}),D(y))$, and the tree
bound is $D=\min(P,\min_{\text{open}}\beta)$.  The branching variable maximizes
\begin{equation}
\min\Bigl(1,\frac{\rho_j}{0.05\max_i\rho_i}\Bigr)\Bigl(\sum_{k\ni j}|\hat w_k-\hat z_{a_k}\hat z_{b_k}|
+\sum_{a\ni j}|\hat t_a-\phi_a(\hat\alpha_a)|\Bigr)
\end{equation}
at the LP point, where $\rho$ is the width relative to the tree's root box.  The split is at the LP point clamped to
the middle 80\%.  The tree stops when the gap closes, after 30\,s, when its root bound lies further below the GPU
root bound than that bound lies below $P$, on a plateau, or, for a hand-off, at $\max(256,32n)$ open leaves.  For
incumbents, IPOPT starts from the best node's LP point at nodes $128,256,\dots,2048$ until one start improves the
incumbent.  At the end, the leaf LP points, leaf LPs re-solved with their rows at half the tolerance, and blends
towards the incumbent are checked at $10^{-6}$, and IPOPT starts from the three lowest leaves.  A verified
improvement lowers the cutoff.

\paragraph{Weak-duality and Farkas certificates.}  For an LP $\min\{c^\top z+c_0:Gz\le h,\ z\in[l^z,u^z]\}$ and any
$y\ge0$,
\begin{equation}
D(y)=c_0-h^\top y+\sum_j\min\bigl(r_jl^z_j,\,r_ju^z_j\bigr)\ \le\ \text{LP optimum},\qquad r=c+G^\top y.
\label{eq:appb_weak_dual}
\end{equation}
The CPU LP solver and PDHG only propose $y$; $D(y)$ is evaluated with outward rounding (\cref{appsec:comp_numerics}) by one
routine, on CPU tensors for the CPU solver and on the GPU for PDHG.  A box is declared empty only when some $y\ge0$ gives
$\min_{z\in[l^z,u^z]}y^\top(Gz-h)>0$, that is,~\eqref{eq:appb_weak_dual} with $c=0$, $c_0=0$ is positive.  If the rows
used include the incumbent cutoff $s_0\le P$, this proves only that no point with $f_0\le P$ exists in the box (the box
is retired by dominance, not shown infeasible for the original model).  Such a $y$ is a
Farkas ray, the multipliers of the phase-one LP $\min\{\eta:Gz-\eta\mathbf1\le h,\ z\in[l^z,u^z],\ \eta\ge0\}$, a single row of the propagation, or a
probe interval that misses the box.  Selected lifted-LP solves, notably those with RLT, are repeated without presolve when presolve calls the LP infeasible
and budget remains; an infeasibility status alone never suffices without a separately evaluated certificate.

\paragraph{Hand-off of open leaves.}  The open leaves $(B_i,\beta_i)$ cover every point of the box that satisfies the
rows with $f_0<P$, so the root bound rises to $\min(\min_i\beta_i,P)$.  At most 8{,}192 leaves are copied to the GPU
once, and they become the first boxes of the GPU tree with their bounds (\cref{fig:appb_lifted_root_to_gpu}b).

\paragraph{Certificate transfer.}  For the last LP of each open leaf (box $B$, cutoff $P_0$), $r=c+G^\top y$ and
$\kappa=c_0-h^\top y$ are computed once in exact rational arithmetic and stored as floats with error bounds.  On
$\mathcal L(B,P_0)$, $c^\top z+c_0=\kappa+r^\top z+y^\top(h-Gz)\ge\kappa+r^\top z$, and the leaf's rows, built on $B$,
hold on every $B'\subseteq B$.  So, while the incumbent is at most $P_0$, every GPU box $B'\subseteq B$ has
\begin{equation}
\beta(B')\ \ge\ \min\Bigl\{P_0,\ \kappa+\sum_j\min\bigl(r_jl'_j,\,r_ju'_j\bigr)\Bigr\},\qquad
[l',u']=Z(B')\cap[l^z,u^z],
\end{equation}
and it keeps the larger of this and its own bound.  The value grows as $B'$ shrinks, through the columns with
$r_j\ne0$; at an exact optimal primal-dual pair these sit at a bound in the leaf LP, but this is only intuition, since
validity does not require complementarity.  In exact arithmetic the value is nondecreasing for nested intervals; the
floating-point values are kept monotone by the explicit maximum with the box's own bound.  Each batch matches its boxes to the scopes on the host (a box
hierarchy) and evaluates the bounds on the GPU, without an LP (\cref{fig:appb_lifted_root_to_gpu}c).

\paragraph{Auxiliary bounds for CROWN.}  Inside its box, an OBBT interval of an atom or product column holds for
every point that satisfies the rows with $f_0\le P_0$; so does the column box of a leaf certificate.  Every graph node
whose exact form is affine in the columns, $v=a_0+\sum_ja_jz_j$, gets the interval
$a_0+\sum_j[\min,\max](a_j\underline z_j,a_j\overline z_j)$, summed in exact rationals and rounded outward, and this
interval is intersected with CROWN's intermediate bounds of $v$.  These bounds hold only below $P_0$, so the
objective bound they produce is capped at $P_0$.  Matching and the exact sums run on the host, once per new scope;
the intersection runs on the GPU.

\paragraph{Per-box lifted bound in the GPU tree.}  Every GPU batch also bounds its boxes over $\mathcal L(B,\infty)$,
with no cutoff row and no propagation.  Rows that are not lifted, and the objective when it is not lifted, enter as
their CROWN planes on $B$.  The LPs of the lowest-bound boxes that fit in $2.5\cdot10^7$ dense coefficients run 500
PDHG iterations on the GPU, with $D(y)$ certified every 50.  A CPU LP solver handles the 16 lowest in floating point (0.5\,s
per batch).  The RLT and epigraph rows are assembled per box on the host.  Where this bound sets a box's bound, the
box branches on the largest violation (\cref{appsec:comp_bab}).

\paragraph{Continuation.}  The open leaves still cover every point below a lower incumbent, so the lifted tree can
resume from them; if the root had no incumbent, the tree starts here instead.  It runs only while the GPU tree's
bound has not passed the lifted tree's $D$, and the GPU tree counts as ahead until it has visited twice as many boxes
as it received.  The tree is first continued once at the same incumbent, to measure its progress.  After that, a
continuation needs an incumbent improvement of more than $10^{-4}$ (relative) or a previous continuation that halved
its gap.  It also needs a measured gain per second (bound gain beyond the GPU bound, plus incumbent gain) that is
positive and at least the GPU tree's.  A call takes at most $\min(30\,\mathrm{s},\,0.3\times$ the remaining time$)$,
and all calls together at most 60\,s unless the last one halved its gap.  Each call runs on the CPU while the GPU
loop waits.

\subsection{One Step of the GPU Branch-and-Bound}
\label{appsec:comp_bab}

\begin{figure}[h]
\centering
\includegraphics[width=\linewidth]{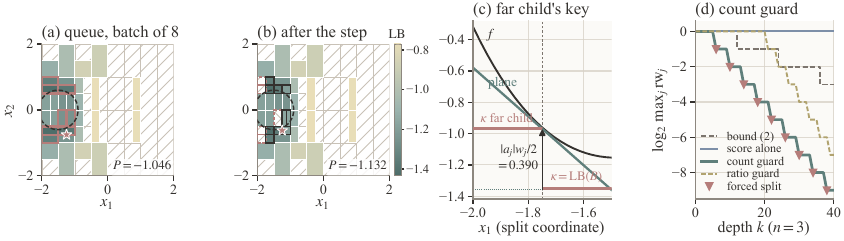}
\caption{\textbf{One GPU branch-and-bound step on a toy}: $\min x_1^4/4-x_1^2+0.3x_1+x_2^2/2$ subject to
$(x_1+1.5)^2+x_2^2\ge0.36$ on $[-2,2]^2$, with lowered tangent planes in place of CROWN.  (a)~Queue at
iteration~9, coloured by $\mathrm{LB}$; the 8 lowest boxes form the batch (red); hatching marks regions removed
by $\mathrm{LB}\ge P$ (/) or certified empty ($\times$).  (b)~After the step: 7 boxes split (black), 4 removed by
$\mathrm{LB}\ge P$ and 1 certified empty (dashed); $D=-1.434\le p^\star\approx-1.255\le P=-1.132$.  (c)~The far
child's objective-plane minimum rises by $|a_j|w_j/2=0.390$ (\cref{eq:appb_far_child}), and its key by 0.382 above
the parent's bound.  (d)~A path on which the score never picks
$x_3$ ($n=3$, coefficients $(10^3,10^3,0)$): the count guard forces a widest split every fourth step and stays
below \cref{eq:appb_count_guard}.}
\label{fig:appb_bab_step}
\end{figure}

An iteration selects a batch, runs the scheduled primal search on it (\cref{appsec:comp_primal}), bounds and
prunes it, and splits the survivors (\cref{fig:appb_bab_step}); per-box work is batched on the GPU, while the
host keeps $P$ (the incumbent value, $\mathrm{UB}$ in \cref{subsec:sbab}) and $D$.

\paragraph{Selection (GPU).}  The batch $\mathcal B_t$ holds the $k_t$ boxes of $\mathcal Q$ with the lowest
keys, $\max_{\mathcal B_t}\mathrm{LB}\le\min_{\mathcal Q\setminus\mathcal B_t}\mathrm{LB}$ (a top-$k$ in every
iteration), so it contains a box of minimum bound, as \cref{app:termination} requires.  $k_0=2$; $k$ doubles
while the last batch was full and less than $45\%$ of the device memory is reserved, and is halved for good
after an out-of-memory error, the batch returning to $\mathcal Q$ unchanged.

\begin{figure}[h]
\centering
\includegraphics[width=\linewidth]{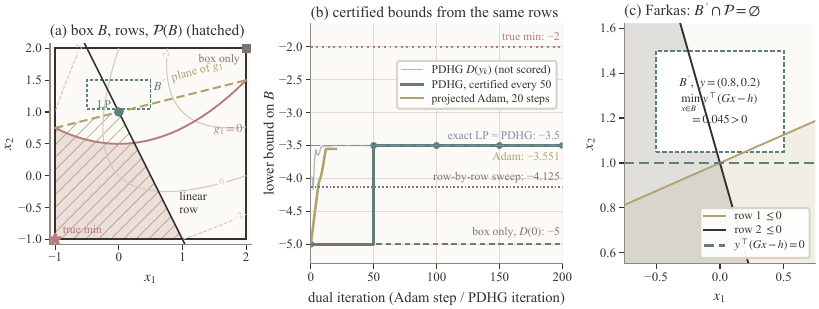}
\caption{Bounding from CROWN planes on a toy: $\min\,x_1x_2-\frac12x_1^2-\frac12x_2^2+x_1+x_2$ s.t.\
$g_1=x_2-\frac14x_1^2-\frac12\le0$ and $x_1+\frac12x_2\le\frac12$ on $B=[-1,2]^2$.  (a) Feasible set (shaded),
$\mathcal P(B)$ from the CROWN plane of $g_1$ and the linear row (hatched), and the minimizers of the objective plane
over $B$ and over $\mathcal P(B)$ and of the problem.  (b) Bounds from the same rows: box only $-5$ (interval
arithmetic: $-8$); the row-by-row single-row duals of constrained concretization, $-4.125$; 20 projected Adam steps,
$-3.551$; PDHG, certified at iteration 50, and the exact LP, $-3.5=p_{\mathrm{LP}}(B)$; the minimum is $-2$.  (c) On a child
box $B'$ each row alone is satisfiable, but the phase-one multipliers $y=(0.8,0.2)$ give
$\min_{x\in B'}y^\top(Gx-h)=0.045>0$, so $B'$ is removed.}
\label{fig:appb_plane_lp}
\end{figure}

\paragraph{Bounding.}  The tensor bounding kernels and the admitted dense PDHG batches run on the GPU, batched over the boxes.  Host work
includes the symbolic and scope processing, the RLT and epigraph row assembly and the scheduled LP solves;
budget-excluded boxes keep their other bounds.
Interval bounds are computed first; the boxes they do not settle are bounded by CROWN, which yields the planes
$\underline f_i(x)=a_i^\top x+c_i$ of \eqref{eq:method_linear_bounds}.  Every concretization (objective and
intermediate nodes) is constrained by the rows available before the pass: the parent's planes (the objective plane
cut at the incumbent included), the linear rows and the presolve cuts, at most the 16 non-redundant rows with the
largest maximum over the box.  The rows act one after another, each through its single-row dual
\[
\beta_k\in\arg\max_{\beta\ge0}\ \min_{x\in B}\bigl[(a+\beta a_k)^\top x+b+\beta c_k\bigr],\qquad
(a,b)\leftarrow(a+\beta_ka_k,\,b+\beta_kc_k),
\]
a concave piecewise-linear function of $\beta$ maximized over $\beta=0$ and the nonnegative breakpoints
$-a_j/a_{kj}$, $a_{kj}\ne0$ (a fused kernel in double precision); a row whose slope stays positive as $\beta\to\infty$
proves the box empty.  After branching, each
child is clipped with the fresh planes by one sweep of \eqref{eq:relaxed_clipping}.  The objective bound is raised by
the dual $d_B(\lambda)$ of \eqref{eq:dual_lb} over the fresh rows $Ax+c\le0$ (constraint planes, linear rows and
presolve cuts; we write $G=A$, $h=-c$).  Every $\lambda\ge0$ gives $d_B(\lambda)\le p_{\mathrm{LP}}(B)$
(\cref{prop:dual}), so three methods only propose multipliers (\cref{fig:appb_plane_lp}).  (i) Twenty projected
Adam steps \eqref{eq:adam_step} from $\lambda=0$ with learning rate $0.05$, keeping the best iterate (not
outward-rounded, \cref{appsec:comp_numerics}).  (ii) PDHG on the same LP: with $x=\underline x+w\odot z$, $z=D_c\tilde z$,
$w=\overline x-\underline x$, ten rounds of Ruiz equilibration $\tilde G=D_rG\operatorname{diag}(w)D_c$ and diagonal
steps $\tau_j=1/\max(\sum_i|\tilde G_{ij}|,10^{-12})$, $\sigma_i=1/\max(\sum_j|\tilde G_{ij}|,10^{-12})$,
\[
\tilde z^+=\Pi_{[0,\bar z]}\bigl(\tilde z-\tau\odot(\tilde a_0+\tilde G^\top\tilde y)\bigr),\qquad
\tilde y^+=\bigl[\tilde y+\sigma\odot(\tilde G(2\tilde z^+-\tilde z)-\tilde h)\bigr]_+,
\]
with $\tilde a_0=D_c\operatorname{diag}(w)a_0$, $\tilde h=D_r(h-G\underline x)$ and $\bar z=D_c^{-1}\mathbf 1$,
for 500 iterations (which only improve the candidate multipliers; certification needs no convergence), warm-started from the batch mean of the previous multipliers; every 50 iterations, and at
$\lambda=0$, $d_B$ is evaluated at $\lambda=D_r\tilde y$ with outward rounding and the best value is kept.  (iii) For
the 16 boxes with the lowest bound (when $mn\le2\cdot10^5$), the LP solved in floating point on the CPU in a pool of at most
eight threads; its multipliers $\max(-\text{marginals},0)$ are scored by the same rigorous $d_B$ on the GPU.  When
the solver reports infeasibility, the multipliers $\nu$ of the phase-one LP $\min\{t:Ax+c\le t\mathbf 1,\ t\ge0,\ x\in
B\}$ are checked as a Farkas certificate \eqref{eq:infeas_cert}, $\min_{x\in B}\nu^\top(Ax+c)>0$, under the same
rounding, and the box is removed.  The quadratic-tangent plane (a second LP over the same rows, again PDHG plus exact
LP solves on its 16 lowest boxes), the lifted LP rows, the transferred certificates and the objective floors of
\cref{appsec:comp_compile,appsec:comp_dual} each contribute a valid bound, and the box keeps the largest.

\paragraph{Pattern propagation.}  When the objective is affine and every row is affine, a single product, or a
signed sum of squared differences (in this benchmark, only the \texttt{kall} circle-packing models), three
rounds of outward interval propagation per box replace the bounding step, and boxes are split along the
coordinate with the largest uncertainty of the squared terms.

\paragraph{Pruning (GPU).}  A box is removed if $\mathrm{LB}(B)\ge P$ (equality included: no strictly better
point remains), if a row's lower bound exceeds its threshold, $\underline f_i(B)>\tau_i$, or if
$\mathrm{LB}(B)=+\infty$, an emptiness certificate of the constrained concretization or of an LP
(\cref{eq:infeas_cert}; with a cutoff row among the rows used, a dominance certificate).  Children are tested again on entry, together with $l_j>u_j$ after clipping and a
positive minimum of the linear rows.

\paragraph{\SBReweight{} branching score (GPU).}  With $a_r^\top x+b_r$ ($r=0,\dots,m$; $r=0$ the objective) the linear bounds
of the box (\cref{eq:method_linear_bounds}), $w=u-l$ and the Lagrangian multipliers $\lambda^{\mathrm{best}}\ge0$, the default choice is
\[
j^\star=\operatorname*{arg\,max}_j\ S_j,\qquad S_j=w_j\Bigl(|a_{0j}|+\sum_{i=1}^{m}\lambda_i^{\mathrm{best}}|a_{ij}|\Bigr),
\]
the \SBReweight{} score of \cref{eq:multiplier_score} (the relative-width fallback and the count guard below are
part of the evaluated solver).  Each term $\frac12|a_{rj}|w_j$ is, for its row separately, the
gain of that row's box minimum: for the half $B^{+}$ of a midpoint split along $j$ that excludes the minimizing vertex of row
$r$, only the $j$-th term of \cref{eq:lbp_corner} changes,
\begin{equation}
\min_{x\in B^{+}}\bigl(a_r^\top x+b_r\bigr)=\min_{x\in B}\bigl(a_r^\top x+b_r\bigr)+\tfrac12|a_{rj}|\,w_j ,
\label{eq:appb_far_child}
\end{equation}
which for $r=0$ raises the objective plane's minimum on the far child and for a constraint row can prune it.  The
$\lambda$-weighted sum over rows makes $S_j$ a branching heuristic only: it does not guarantee an equal increase in either
child's recomputed bound or in the global lower bound.

\paragraph{Overrides (GPU; lifted LP points partly from the CPU).}  With $\mathrm{rw}_j=w_j/w_j^0$ the width
relative to the root box, four rules replace $j^\star$ in turn.  (i)~Where the lifted LP
(\cref{appsec:comp_lifted}) sets the bound, $j^\star$ maximizes $\sum_{k\ni j}|z^\star_k-\phi_k(z^\star)|$, the
violation at the LP point $z^\star$ of the lifted products and atoms containing $x_j$, over
$\mathrm{rw}_j\ge0.05\max_i\mathrm{rw}_i$ (if positive); halving $x_j$ halves the McCormick gap bound $w_iw_j/4$
of each product $x_ix_j$.  (ii)~If $\mathrm{rw}_{j^\star}<10^{-3}\max_i\mathrm{rw}_i$, the relatively widest
coordinate is taken.  (iii)~The relatively widest variable of the rotation pairs $a(x_i^2+x_j^2)=k$ is split
while its relative width exceeds $0.15$.  (iv)~Without constraint rows and with $n\le6$, both children of every
coordinate are bounded for the $10^4$ lowest boxes, and $j^\star$ maximizes
$\max(\Delta_j^-,10^{-12})\max(\Delta_j^+,10^{-12})$, $\Delta_j^\pm$ the children's gains over $\mathrm{LB}(B)$.

\paragraph{Count guard (GPU).}  With $c(B)$ the number of splits since the box was last split along its
relatively widest coordinate, the last rule is $j^\star\leftarrow\arg\max_j\mathrm{rw}_j$ if
$c(B)\ge n_{\mathrm{open}}(B)$, the number of coordinates of positive width.  A widest split then comes at least
once in every $n+1$ splits, and in exact arithmetic (or while every committed midpoint split genuinely halves its
interval) $n$ of them halve $\max_j\mathrm{rw}_j$, since widths never grow; a box at depth
$k$ thus satisfies
\begin{equation}
\max_j\mathrm{rw}_j(B)\le2^{-\lfloor k/(n(n+1))\rfloor},
\label{eq:appb_count_guard}
\end{equation}
the exhaustive subdivision of \cref{app:termination}, whatever rules (i) to (iv) chose (\cref{fig:appb_bab_step}d),
with $\mathrm{rw}_j=0$ when $w_j^0=0$.  In binary64 an interval of adjacent floats has its midpoint at an endpoint, so
the rate cannot be claimed below representable resolution; finite termination also needs the relaxation, primal-recovery
and infeasibility-detection assumptions of \cref{app:termination}.

\paragraph{Split, clipping and child keys (GPU).}  The box is cut at $m=(l_{j^\star}+u_{j^\star})/2$ into
$\{x_{j^\star}\le m\}$ and $\{x_{j^\star}\ge m\}$, which share the floating-point $m$ and cover $B$ exactly.  Each
child $C$ is clipped by the parent's rows, the objective row at right-hand side $P$ included
(\cref{subsec:dual}), and queued with the key
\begin{equation}
\kappa(C)=\max\Bigl(\mathrm{LB}(B),\ \operatorname{rd}_{\downarrow}\min_{x\in C}\bigl(a_0^\top x+b_0\bigr)\Bigr),
\label{eq:appb_child_key}
\end{equation}
where $a_0^\top x+b_0$ is the parent's objective plane and $\operatorname{rd}_\downarrow$ subtracts the
forward-error bound $(n+4)\mathbf u\bigl(\sum_i|a_{0i}|(|c_i|+r_i)+|b_0|\bigr)$ of the vertex formula ($c$, $r$ the
centre and half-width of $C$, $\mathbf u$ the unit roundoff) and rounds toward $-\infty$; the allowance assumes finite
intermediates, covered underflow and $2(n+4)^2\mathbf u\le1$.  With $q_B=\min_{x\in B}(a_0^\top x+b_0)$ and
$\Delta=\tfrac12|a_{0j^\star}|w_{j^\star}$, by \cref{eq:appb_far_child} the far child's key is, up to rounding,
$\max(\mathrm{LB}(B),q_B+\Delta)$: it rises above $\mathrm{LB}(B)$ by $\max(0,q_B+\Delta-\mathrm{LB}(B))$, which is
less than $\Delta$ when the parent's bound is stronger than its objective plane.

\paragraph{Reported bound and gap rule (GPU and host).}  A removed box holds no feasible point or has
$\mathrm{LB}(B)\ge\theta$, where $\theta$ is the floor stored at its removal: the incumbent value $P'\ge P_t$ for
dominance, and for the epsilon pruning of the evaluated configuration (the $\mathrm{M3}$ reduction) the rounded-up
threshold $P-\max(\epsilon_{\mathrm{abs}},\epsilon_{\mathrm{rel}}\max(1,|P|))$.  With $\theta_{\min}$ the smallest such
$\theta$, the run reports
\begin{equation}
D_t=\max\Bigl(D_{t-1},\ \min\bigl(\theta_{\min},\ \min_{B\in\mathcal Q_t}\mathrm{LB}(B)\bigr)\Bigr)\ \le\ p^\star,\qquad
p^\star\le P_t,
\label{eq:appb_reported_bound}
\end{equation}
($D_0$ the root bound, $p^\star$ with the tolerance absorbed), which is a valid lower bound when every stored floor
$\theta$ is.  Under epsilon pruning the floor is below the incumbent, so this bound may be weaker than $P_t$ even
when the queue is empty.  The cap matters: the box holding an optimum is
removed once its bound reaches an exact incumbent, and the survivors' minimum may then exceed $p^\star$.  The run
stops when \cref{eq:exp_solver_criterion} holds with $P-D\ge-10^{-5}$, when $\mathcal Q=\emptyset$ (then
$P=p^\star$ if only exact dominance retirement was used, or infeasibility if $P=+\infty$; with tolerance-based
retirement, queue exhaustion certifies the stated gap, not exact optimality), or at the deadline; a bound above $P$ by more than $10^{-5}$ is an
error and never certifies.

\subsection{Primal Search}
\label{appsec:comp_primal}

\paragraph{Acceptance gate (GPU, one point at a time).}  Every proposed point is evaluated on the compiled
original model and accepted only if
\begin{equation}
\bar x\in B_0,\qquad\max_i\,[f_i(\bar x)]_+\le10^{-6},\qquad\max_i\,[a_i^\top\bar x-b_i]_+\le10^{-6}
\label{eq:appb_gate}
\end{equation}
over the nonlinear and linear rows (equalities as two rows), i.e.\ \cref{eq:feasibility_gate} with
$\tau_{\mathrm{feas}}=10^{-6}$; then $P\leftarrow\min(P,f_0(\bar x))$.  The in-tree ALM applies the same test to
its batch residuals.  Candidate generators enforce $\bar x\in B_0$ by construction (clamping or projection), and the routine checks the
unscaled nonlinear and linear residuals in floating point.  No merit or rescaled residual decides, so $P\ge p^\star$
up to the floating-point evaluation of $f_0$ and the rows (no outward-rounded objective evaluation).

\paragraph{Root ALM-PGD (GPU).}  At the root, ALM-PGD runs on the GPU, IPOPT in a CPU thread beside it, and LP
alternation in a CPU process when the exact forms are a QCQP with a product term.  The ALM minimizes
\cref{eq:hpr_alm} over $R=20{,}480$ starts on rescaled rows,
\[
\phi_\rho(x;\mu)=f_0(x)+\frac\rho2\sum_i\Bigl[\frac{f_i(x)}{s_i}+\frac{\mu_i}{\rho}\Bigr]_+^2-\frac{\|\mu\|^2}{2\rho},
\qquad s_i=\max\bigl(1,\mathrm{RMS}_r[f_i(x^{(r)})]_+\bigr),
\]
$s_i$ halved toward 1 after each outer iteration, with $\mu\leftarrow\mu\odot s_{\mathrm{new}}/s_{\mathrm{old}}$
to keep raw multiplier units; the scales are computed separately for projected and unprojected root starts when the
affine projection below is active.  Inner step $t\le T=200$ of outer iteration $k=0,1,\dots,19$ is
\[
x\leftarrow\Pi_{B_0}\Bigl(x-\operatorname{clip}\bigl(\eta_k(1-0.9\tfrac tT)\max_jw_j\,\hat d_t,\ \pm\eta_kw\bigr)\Bigr),
\qquad\eta_k=0.01\cdot0.6^{k},
\]
with $\hat d_t=m_t/\bigl(\sqrt{v_t/(1-\beta_2^t)}+\epsilon\bigr)$ the Adam direction (the first moment is not
bias-corrected, unlike \cref{eq:projected-adam}); then $\mu\leftarrow[\mu+\rho f/s]_+$, and
$\rho\leftarrow\min(2\rho,10^6)$ ($\rho_0=1$) where the largest violation fell by less than $25\%$.  Start 0 is
the box centre and the rest are uniform, except a quarter scaled log-uniformly around $\Pi_{B_0}(0)$ and a
quarter projected onto the affine equalities $Ex=d$ detected as opposite rows with constant Jacobians,
$x\leftarrow\Pi_{B_0}(x-E^\top(EE^\top)^{+}(Ex-d))$, at the start and after every Adam step.  The early-exit test
$\|\nabla\phi_\rho\|\le\max(\rho^{-0.8},10^{-5})$ synchronizes with the host at every inner step.

\paragraph{Restoration first (GPU, host control).}  Once IPOPT has a verified value $P$, progress is measured by
\[
r(x)=\max\Bigl(\frac{v(x)}{10^{-6}},\ 1+\frac{[f_0(x)-(P-\tau)]_+}{\tau}\Bigr)\ge1,\qquad
\tau=\max\bigl(10^{-5},10^{-4}\max(1,|P|)\bigr),
\]
$v$ the largest violation, so $r=1$ at a feasible point better than $P$ by the gap tolerance.  The ALM runs its
endgame (below) once and stops if the restored basin matches $P$ within $\tau$; otherwise it continues only while
the measured decrease of $\log\min r$ per inner step can reach 0 within the remaining inner steps.

\paragraph{Gauss-Newton/Levenberg-Marquardt restoration (GPU).}  Domains without a point that passes the gate
are restored from their archived point and least violated end points (256 seeds and 50 iterations at the root;
8 seeds per box, 30 iterations and 2\,s in the tree).  In unit-box coordinates $z$, with rows normalized as
$t_i=(f_i+\tfrac12\cdot10^{-6})/\|\nabla_zf_i\|$, an approximate restoration step is computed (not the exact solution of a constrained trust-region subproblem): the
damped Gauss-Newton direction
\begin{equation}
(J_{\mathcal A}^\top J_{\mathcal A}+\lambda I)\,d_{\mathrm{GN}}=-J_{\mathcal A}^\top t_{\mathcal A},\qquad
d=d_{\mathrm{GN}}\min\bigl(1,\Delta/\|d_{\mathrm{GN}}\|_\infty\bigr),
\label{eq:appb_lm}
\end{equation}
over the coordinates not held at a bound, with $J$ the Jacobian of the normalized residuals with the row norms frozen at
the linearization (norm 1 for zero-gradient rows), the active set $\mathcal A$ growing from $\{t_i>0\}$ by the rows the
step would violate (three rounds), $\lambda$ from $10^{-6}\max\operatorname{diag}(J_{\mathcal A}^\top J_{\mathcal A})$
and $\Delta$ adapted to accepted steps, and the length the first $\alpha\in\{1,\dots,2^{-7}\}$ (one batch) that
decreases $\sum_i([f_i]_+/\|\nabla f_i\|)^2$ or $\max_i[f_i]_+$.  At the root, still infeasible domains then get
250 normalized descent steps on $\max_i[f_i]_+$.

\paragraph{SQP polish with an interior-point QP (GPU).}  ALM stops at the first feasible points, under a large
penalty, so up to 32 distinct points (half the best feasible, half the best of the $5\%$ least violated) and
the first 32 starts are polished (50 iterations, 5\,s; in the tree 8 points, 1\,s) by solving the always
feasible elastic QP
\begin{equation}
\min_{d,\,v}\ \nabla f_0^\top d+\tfrac12d^\top Hd+\nu\,\mathbf 1^\top v\quad\text{s.t.}\quad
f(z)-0.9\cdot10^{-6}+Jd\le v,\quad v\ge0,\quad0\le z+d\le1,
\label{eq:appb_qp}
\end{equation}
with $H$ a Powell-damped BFGS Hessian of the Lagrangian.  A Mehrotra interior-point method solves it with one
$n\times n$ Cholesky factorization of $H+J^\top\operatorname{diag}(\sigma)J+\operatorname{diag}(\omega_{\mathrm{lo}}
+\omega_{\mathrm{hi}})$ per iteration ($v$ eliminated; $\sigma,\omega$ from the multiplier-to-slack ratios), at most
20 iterations replayed from a CUDA graph; steps pass an Armijo test on $f_0+\nu\sum_i[f_i-0.9\cdot10^{-6}]_+$ with
one second-order correction.  Aiming at $0.9\cdot10^{-6}$ uses the tolerance the lower bounds also allow.

\begin{figure}[h]
\centering
\includegraphics[width=\linewidth]{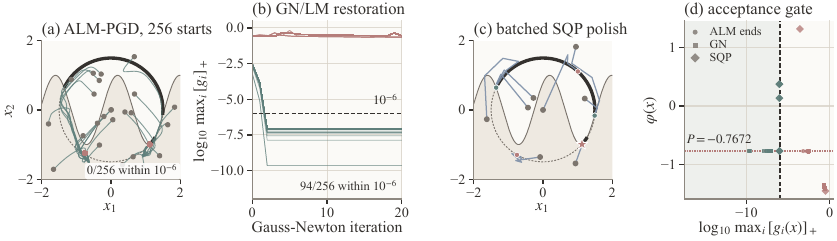}
\caption{\textbf{Root primal pipeline on a toy}: $\min x_2+0.2x_1$ subject to $x_1^2+x_2^2=2.25$ (two rows) and
$x_2\ge\sin(3.5x_1+0.9)$ on $[-2,2]^2$ (feasible set: thick arcs; local minima $-0.767$, $0.134$, $0.373$).
(a)~ALM-PGD from 256 starts (20$\times$200 steps; every 11th trajectory drawn): no end point is within $10^{-6}$.
(b)~Restoration from 256 seeds: 94 reach $10^{-6}$, most in two iterations; the others did not reach the feasible
arcs.  (c)~SQP from 10 points: 8 pass the $10^{-6}$ test (the target is $0.9\cdot10^{-6}$).
(d)~The gate accepts 102 of the 522 points; $P\approx-0.767211$ for the tolerance-relaxed problem (the quoted local
minima refer to zero-slack constraints).}
\label{fig:appb_primal_pipeline}
\end{figure}

\paragraph{IPOPT at the root (CPU).}  IPOPT (MUMPS) uses the exact sparse first and second derivatives of the
structure walk (\cref{appsec:comp_compile}), checked for non-finite values, and starts from the box centre and,
while more than $20\%$ of its effective budget (the shared root budget of 60\,s, see Budgets) remains, from the minimizer of the root CROWN planes
$\arg\min_{x\in B_0}\{a_0^\top x:a_i^\top x+b_i\le0,\ i\ge1\}$.  It targets $f_i\le\tfrac12\cdot10^{-6}$,
polishes an accepted point at $f_i\le0.99\cdot10^{-6}$, and finishes a rejected one in the box at $f_i\le0$.

\begin{figure}[h]
\centering
\includegraphics[width=\linewidth]{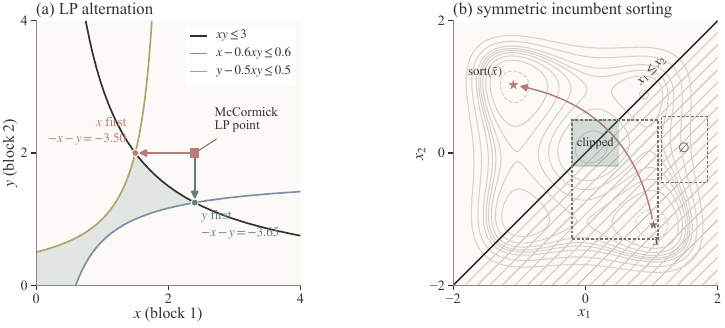}
\caption{(a)~\textbf{LP alternation} on $\min -x-y$ over $[0,4]^2$ with $xy\le3$, $x-0.6xy\le0.6$,
$y-0.5xy\le0.5$.  From the same infeasible McCormick point $(2.4,2)$, a run with $x$ as the first block reaches the vertex $(1.5,2)$
($-3.5$) and a run with $y$ first the global optimum $(2.4,1.25)$ ($-3.65$); these are separate runs, and at both end
points the trust-region LP predicts no decrease.  (b)~\textbf{Symmetric sorting}: the row $x_1\le x_2$ removes the hatched half, $\bar x$ becomes
$\mathrm{sort}(\bar x)$ with the same value $-0.540263$, and clipping shrinks a box crossing the diagonal and
empties one below it.}
\label{fig:appb_lpalt_symmetry}
\end{figure}

\paragraph{LP alternation (CPU process).}  The blocks $x_{B_1},\dots,x_{B_K}$ are the colour classes of the
product graph together with product-free variables, while variables with a unary nonlinear dependence (powers,
squares, atoms) are held fixed, so the objective and every row are affine in one block when the others are fixed, and a
block step is an LP solved on the CPU: from an infeasible point the elastic LP
\[
\min_{x_{B_k},\,s^\pm\ge0}\ \sum_{r:\|A_r\|_\infty>0}\frac{s^+_r+s^-_r}{\|A_r\|_\infty}\quad\text{s.t.}\quad
c_L\le A\,x_{B_k}+c+s^+-s^-\le c_U,
\]
then the objective LP $\min\{\nabla_{B_k}f_0^\top x_{B_k}:c_L\le Ax_{B_k}+c\le c_U\}$, which does not increase $f_0$ from a feasible current point (in exact arithmetic; numerically the acceptance rule
rejects worse incumbents), with inconsistent constant rows handled separately.  Blocks cycle while they progress; then an LP on all rows linearized at the point, in a trust
region, moves all blocks and block LPs restore feasibility.  Seeds are the McCormick LP point in every block
order, ratio fits and McCormick vertices of perturbed costs; feasible points also seed sparse IPOPT.  The LP steps tend to return basic solutions, i.e.\ vertices (\cref{fig:appb_lpalt_symmetry}a), where degenerate
network optima lie; every candidate still passes the feasibility gate.

\paragraph{Symmetric incumbent sorting (CPU).}  For each chain of interchangeable blocks
(\cref{appsec:comp_compile}), the root incumbent and warm start become $\mathrm{sort}(\bar x)$, the blocks permuted
so that $x_{B_1[1]}\le\dots\le x_{B_m[1]}$; block swaps preserve $f_0$ and the multiset of rows (checked exactly),
so the point keeps its value and satisfies the ordering rows used by clipping (\cref{fig:appb_lpalt_symmetry}b).

\paragraph{In the tree (GPU, host schedule).}  At scheduled iterations the ALM runs on every box of the batch
before bounding, with about 1{,}000 starts in total, 20$\times$200 steps, 10\,s and the default parameters
($\eta_0=1$, $\rho_0=10$, gated multiplier updates, no row scaling).  Start 0 of a box is the incumbent clamped
into it (else the dual minimizer $x^{\mathrm{best}}$ of \cref{eq:dual_lb}), $\mu^{(0)}$ is the batch mean of
$\lambda^{\mathrm{best}}$, and up to 10 inherited constraint planes that can be active join the rows.
Restoration follows, plain and then with the steering row $f_0(x)\le P-10^{-4}\max(1,|P|)$, and an SQP polish
ends the call; points still pass the gate on the original rows.  The period $T$ starts at 5, doubles (up to 20)
after three calls without improvement and resets on improvement; after four consecutive failures a call runs
only if $(P-\min_{\mathcal B_t}\mathrm{LB})/(1+|P|)<0.01$.  The lifted LP tree also proposes its leaf points
(\cref{appsec:comp_lifted}).

\subsection{Numerics}
\label{appsec:comp_numerics}
The numerical tensor kernels use double precision; symbolic coefficients use exact rationals and selected
certificates use multiprecision interval arithmetic.  Every LP-multiplier bound (the PDHG and CPU LP solver multipliers of the
plane LP and of the lifted LP, the root dual LPs and the transferred certificates) and every LP infeasibility verdict
is evaluated by one weak-duality routine with outward rounding, whereas interval propagation and constrained
concretization use separate conservative tests.  In that routine, under finite arithmetic and without
underflow or overflow, library errors being assumed negligible: inner products carry the Higham bound
$\gamma_k$ of their summation error, a double-double evaluation is used when the error term could change a
decision, and casts to lower precision round down.  The lifted-LP rows carry outward cushions, and the
posynomial and relative-entropy certificates are evaluated in interval arithmetic.  The concretizations of
the CROWN and interval bounds add a forward-error margin to their inner products and bias sums.  The operator
relaxations themselves, the separable floors and the projected Adam bound of the dual ascent are evaluated in
ordinary double precision without outward rounding.  Incumbents are verified on the original model; for
scoring, each returned point is replayed on the exact double-precision model at $10^{-6}$, and the 22 points
of the instances that only \GEAR{} certifies among \GEAR{}, Gurobi and SCIP (\cref{appsec:win_loss}) were
also re-evaluated in AMPL, with a largest violation of $9.9999994\cdot10^{-7}$.

\subsection{Budgets}
\label{appsec:comp_budgets}
Several steps have wall-clock budgets: bound tightening (5\,s on the McCormick LP and up to 10\,s in total; in
the tree, 5\,s per call and 20\,s in total), the LP branch-and-bound (30\,s), the exact LPs of one GPU batch on the
lifted relaxation (0.5\,s), the quadratic-tangent KKT solve (10\,s), each certificate LP (20\,s), the root ALM
search (60\,s) and IPOPT (30\,s configured, raised to the shared 60\,s root budget in parallel root mode, subject to the
run deadline), the SQP polish (5\,s at the root, 1\,s in the tree), each in-tree ALM call
(10\,s) and its Gauss-Newton search (2\,s).  The solved set therefore depends on the speed of the host.

\section{Benchmark and Experimental Protocol}
\label{appsec:dataset}

\subsection{Source Pools}

The 505-instance benchmark is built from two continuous NLP source pools: 265 \texttt{gams\_world} instances and 240 \texttt{minlplib\_linear} instances.
\GEAR{} and every baseline run on the same instance keys.

\subsection{Canonicalization and Graph Construction}

Each retained instance is converted to
\[
\min_x f_0(x)
\quad
\mathrm{s.t.}
\quad
f_i(x)\le0,\ i=1,\ldots,m,
\qquad
x\in[\ell,u].
\]
Maximization objectives are negated.
Two-sided inequalities are split into paired one-sided inequalities.
Equalities are encoded as paired inequalities with a fixed tolerance.
Objective offsets, source variable bounds, and model scaling are preserved unless the source parser already normalizes them.

\paragraph{PyTorch graph.}

The graph uses one vector input $x\in\mathbb R^n$ and emits
$f_0$ together with the shifted constraints, as in
\eqref{eq:nlp_graph}. Linear constraints are also graph outputs, reproduced exactly by LBP; when the
parser exposes them as an explicit matrix, that matrix is also kept for clipping and the dual. Residual infeasibility uses individual
constraint rows or their nonnegative mixtures and does not require a
separate hard-max graph output. Shared algebraic subexpressions can be
reused when exposed by the parser.

\subsection{Inclusion and Exclusion Policy}

An instance is retained only if it is continuous, has finite variable bounds after canonical preprocessing, uses supported graph primitives, has valid domain guards for singular operators, and can be parsed from its AMPL model into the \GEAR{} graph representation.

\begin{table}[h]
\centering
\small
\begin{tabular}{lr}
\toprule
\textbf{Source / filter stage} & \textbf{Retained count} \\
\midrule
\texttt{gams\_world} continuous computation-graph NLPs & 265 \\
\texttt{minlplib\_linear} continuous computation-graph NLPs & 240 \\
\midrule
Total & 505 \\
\bottomrule
\end{tabular}
\caption{Retained instances by source.}
\label{tab:dataset_sources}
\end{table}

\paragraph{Exclusions.}

Excluded cases fall into structural categories: integer or discrete variables, missing finite bounds, unsupported operators, domain singularities that cannot be split or rewritten, parser failures, inconsistent source metadata, or root memory requirements beyond the configured GPU budget.

\subsection{Translation Precision and Primal-Dual Soundness}
\label{appsec:translation_soundness}

\GEAR{} parses each instance into a PyTorch computation graph and exact
affine rows. Translation can introduce floating-point discrepancies
relative to the original NLP model through evaluation order, intermediate
precision, and transcendental implementations. These discrepancies must
be accounted for before graph-level bounds are interpreted as
original-model certificates.

On the primal side, the source-model validation protocol checks every
accepted incumbent against the original constraints under their stated
tolerances. The reported benchmark solutions were re-substituted into the
original model and checked for feasibility with finite objective values.
In the notation of \eqref{eq:feasibility_gate}, those tolerances have
already been absorbed into the constraint functions, so acceptance is
against zero rather than an additional positive threshold.
On the dual side, the conservative row and scalar evaluations described
in \cref{app:admissibility,app:concretization} give original-model bounds
only when their allowances also cover relevant translation and
propagation errors. Under these conditions the corrected bounds err in
the conservative direction.

\subsection{Baseline Protocol}
\label{appsec:scip_controlled}
\label{app:controlled}
\label{app:baseline}

Every baseline in \cref{sec:exp} is run on the same workstation and the same 505 instance keys as \GEAR{}, with a 180\,s limit.

Gurobi and Xpress parallelize the spatial branch-and-bound tree itself, so each instance is solved with 16 threads, one instance at a time.

SCIP is run in its concurrent mode with 16 threads. In this mode several differently configured solver instances race on the same problem and the first one to terminate determines the result, while each LP is solved single-threaded; instances are again processed one at a time.

The remaining solvers have no native parallel tree search on our benchmark. BARON's thread option parallelizes only its MIP subproblems, which do not arise in these purely continuous models. Couenne, IPOPT, SNOPT, MINOS, and LGO are single-threaded. LindoGlobal does provide a multi-threaded mode, but in our tests it returned objective values inconsistent with the returned points, so we run it single-threaded.

For all of these solvers we saturate the CPU at the level of throughput: 16 instances run concurrently, each with one thread, and each instance keeps its own 180\,s wall-clock limit.

Knitro supports both layouts. On the tuning subset, parallel multistart with 16 threads and 16 concurrent single-thread runs with serial multistart solved the same instances, and we keep the latter.

For every baseline we select a single global configuration, without per-instance tuning. The candidates are drawn from each solver's documented options, covering presolve and bound-tightening intensity, cutting planes, branching, local-search heuristics, and, for the natively parallel solvers, the thread layouts $16\times1$, $8\times2$, and $4\times4$. We evaluate them on a stratified subset of 40 instances, half drawn from instances that are hard for at least one solver and half sampled uniformly at random, and keep the configuration with the most certified solutions.  For SCIP we additionally keep powers of sums unexpanded (\texttt{expr/pow/expandmaxexponent} $=1$); this setting was found on the full benchmark, where it certifies 8 more instances and loses none.
All solvers use the same termination criteria: a feasibility tolerance of $10^{-6}$, a relative optimality gap of $10^{-4}$, and an absolute gap of $10^{-5}$.
Every returned point is re-evaluated by AMPL against all original constraints and variable bounds at $10^{-6}$.
A global solver is credited with an instance only if its own primal and dual bounds satisfy the gap criterion and no other run found a feasible point that contradicts its certificate.
A certificate is also discarded when the solver itself flags it as unreliable (AMPL status \texttt{solved?}, e.g.\ BARON's ``globality is not guaranteed'').
A global solver's \texttt{infeasible} verdict counts as solved when the instance is confirmed infeasible: no run of any solver returned a point feasible at $10^{-6}$ and at least two global solvers report infeasibility within the limit (one instance, \texttt{mosel\_start\_folioqpgraph}, for which \GEAR{} also certifies infeasibility).
A local solver is credited only if its verified objective matches a bound certified by a global solver within the same tolerance.

\subsection{\GEAR{} Logs and Scoring}

\GEAR{} logs provide the last parsed primal incumbent $P$, valid lower bound $D$, final status, wall time, and node count.
Missing final summary lines are counted as unsolved.
A \GEAR{} run is counted solved only when it satisfies the finite-gap test of \cref{eq:exp_solver_criterion} or terminates after all active domains have been pruned.
Baseline runs are scored as in \cref{app:baseline}.

\subsection{Released Tables}

The release lists every instance with its source, name, numbers of variables and constraints, sense, and the functions it uses.
For \GEAR{} and each baseline it gives one row per instance with the final status, primal and dual values, and time; the baseline rows also give the largest constraint violation found when the returned point is re-evaluated in AMPL.

\section{Analysis of the Results and Limitations}
\label{appsec:analysis}

\subsection{Which Problems Benefit Most from the GPU}
\label{appsec:gpu_benefit}
\label{subsec:exp_gpu_benefit}

\begin{table}[h]
\centering
\scriptsize
\setlength{\tabcolsep}{2.6pt}
\renewcommand{\arraystretch}{0.92}
\begin{tabular}{@{}lrrrrrrr@{}}
\toprule
\textbf{Problem class} & $n$ & \GEAR{} & SCIP & Gurobi & Union & \GEAR{}-only & Union-only \\
\midrule
$\le 5$ variables        & 183 & \textbf{182} & 171 & 177 & 182 & 0 & 0 \\
6-20 variables          & 151 & \textbf{145} & 118 & 127 & 140 & 8 & 3 \\
21-100 variables        &  80 & 66 & 58 & 66 & \textbf{71} & 3 & 8 \\
101-1000 variables      &  65 & 28 & 17 & 29 & \textbf{32} & 5 & 9 \\
$>1000$ variables        &  26 & 0 & 0 & 2 & \textbf{2} & 0 & 2 \\
\midrule
box constraints only     &  95 & \textbf{92} & 81 & 81 & 89 & 3 & 0 \\
1-20 constraints        & 235 & \textbf{230} & 197 & 217 & 226 & 6 & 2 \\
$>20$ constraints        & 175 & 99 & 86 & 103 & \textbf{112} & 7 & 20 \\
\midrule
transcendental terms     & 147 & \textbf{131} & 100 & 110 & 118 & 16 & 3 \\
polynomial/rational only & 358 & 290 & 264 & 291 & \textbf{309} & 0 & 19 \\
\bottomrule
\end{tabular}
\caption{Solved instances by problem class. Union: the six global solvers. Transcendental: \texttt{exp}, \texttt{log}, \texttt{sqrt} or trigonometric terms.}
\label{tab:exp_gpu_benefit}
\end{table}

\Cref{tab:exp_gpu_benefit} breaks the comparison of \cref{sec:exp} down by problem class.  Three patterns emerge.

\paragraph{\GEAR{} leads on low- to medium-dimensional problems with few constraints and
composite nonlinearities.}
On instances with at most 20 constraints, \GEAR{} certifies 322/330
(97.6\%), against 278 for SCIP and 298 for Gurobi.  On instances with transcendental terms, it
certifies 131/147 (89.1\%), against 100 and 110; all 16 instances that \GEAR{} solves and no
global solver certifies belong to this class (e.g.\ \texttt{Paviani}, \texttt{hs101},
\texttt{hs102}, \texttt{hs110}, \texttt{schitt\_s380}, \texttt{kriging\_peaks}, \texttt{pricing050},
and \texttt{polygon25} to \texttt{100}, which a special-structure certificate closes at the root;
\cref{appsec:components}).
On these problems, the difficulty is the nonconvexity of composed exponential, logarithmic and
power terms over a small box.  Factorable relaxations of such compositions tend to be loose, so a
CPU spatial branch-and-bound can process many weak nodes and still end with a gap.  \GEAR{} closes most of
these instances at the root (\cref{appsec:comp_dual,appsec:win_loss}); on those that need a tree, it evaluates a clipped
linear relaxation for many input domains in one GPU batch, which spreads the cost of a step over its boxes.

\paragraph{The GPU does not replace cutting planes on large structured models.}
On instances with more
than 20 constraints, \GEAR{} certifies 99/175 (56.6\%) against 103 for Gurobi, and 20 of the 22
instances that only the global union solves fall into this class: pooling, wastewater and
water-network models, and circle-packing (\texttt{kall}) models with tens to thousands of bilinear
constraints.  Their certificates rely on LP-based branch-and-bound with RLT and other cutting planes
and on bound propagation at every node, and their dimension (up to 2,399 variables) makes splitting
the input box expensive; no instance with more than 1,000 variables is solved by \GEAR{}.

\paragraph{The benefit is coverage, not latency.}
On the 405 instances that both \GEAR{} and at least
one global solver certify, the fastest global solver needs a median of 0.03-0.06\,s across size
classes, against 4-29\,s for \GEAR{}, which is faster on only 8 of them.  The GPU therefore
extends the set of instances that can be certified within the budget; on instances that CPU solvers
already close quickly, its fixed start-up cost dominates.

\subsection{Where \GEAR{} Wins and Where It Loses}
\label{appsec:win_loss}

We compare \GEAR{} instance by instance with the two strongest baselines, Gurobi and SCIP (the run of
\cref{sec:exp}).  Of the 505 instances, 399 are certified by \GEAR{} and by at least one of them, 22 only by
\GEAR{}, 18 only by Gurobi or SCIP, and 66 by none.  For each of the 40 disagreements we read the logs of
all three solvers and assign one cause (\cref{apptab:win_loss}); the numbers below are read from those
logs.  ``Required gap'' is the tolerance of \cref{eq:exp_solver_criterion}.

\begin{table}[h]
\centering
\scriptsize
\setlength{\tabcolsep}{3pt}
\renewcommand{\arraystretch}{1.08}
\begin{tabular}{@{}p{0.27\linewidth}p{0.33\linewidth}rp{0.30\linewidth}@{}}
\toprule
\textbf{Cause} & \textbf{Instances} & \textbf{n} & \textbf{Evidence from the logs} \\
\midrule
\multicolumn{4}{@{}l}{\emph{\textup{\GEAR{}} certifies, Gurobi and SCIP do not (22)}} \\
Posynomial objective: certificate in logarithmic coordinates & \texttt{hs101}, \texttt{hs102}, \texttt{cute\_hs101}, \texttt{cute\_hs102}, \texttt{schitt\_s380} & 5 & baseline gaps end at 34 to $10^{4}\times$ the required gap after $8.7\cdot10^{3}$ to $2.7\cdot10^{6}$ nodes \\
Polygon area under a diameter bound: special-structure certificate & \texttt{polygon25}, \texttt{50}, \texttt{75}, \texttt{100} & 4 & baseline gaps end at $5.0\cdot10^{3}$ to $9.8\cdot10^{3}\times$ the required gap \\
Symmetric log-barrier objective: bound tightening collapses the box & \texttt{Paviani}, \texttt{hs110}, \texttt{cute\_hs110} & 3 & baseline gaps end at $1.9\cdot10^{3}$ to $3.0\cdot10^{3}\times$ the required gap after $6.0\cdot10^{4}$ to $1.4\cdot10^{6}$ nodes \\
\GEAR{}'s own bounds where the baselines' relaxations stay loose & \texttt{pricing050}, \texttt{cute\_spanhyd}, \texttt{cute\_probpenl}, \texttt{Shekelfox10} & 4 & baseline gaps end at $1.7\cdot10^{3}$ to $10^{4}\times$ the required gap \\
A baseline's point fails the $10^{-6}$ re-check & \texttt{kriging\_peaks-full050}, \texttt{100}, \texttt{200}, \texttt{ex8\_4\_4}, \texttt{elena\_hs54}, \texttt{cute\_avion2} & 6 & baseline violations $1.4\cdot10^{-6}$ to $1.7\cdot10^{-5}$, one of 0.16 (\texttt{elena\_hs54}); \GEAR{}'s at most $9.9999994\cdot10^{-7}$ \\
\midrule
\multicolumn{4}{@{}l}{\emph{Gurobi or SCIP certifies, \textup{\GEAR{}} does not (18)}} \\
Same root relaxation; the baselines' tree search is far faster & 6 \texttt{wastewater}, 3 \texttt{waterund}, 5 \texttt{kall}, \texttt{pooling\_digabel16}, \texttt{ex5\_3\_3} & 16 & Gurobi: $2.9\cdot10^{3}$ to $4.0\cdot10^{4}$ nodes/s; \GEAR{}'s CPU lifted-LP tree: 1.1-207 nodes/s (median 27) \\
No lifted LP applies; \GEAR{}'s box bounds stay loose & \texttt{eq6\_1} & 1 & after $3.0\cdot10^{6}$ boxes $D=-165$ against $P=670.69$; Gurobi closes in $8.3\cdot10^{5}$ nodes (70\,s) \\
Beyond \GEAR{}'s lifted-LP size limit & \texttt{pooling\_sppa9tp} & 1 & \\
\bottomrule
\end{tabular}
\caption{The 40 instances on which \GEAR{} and the union of Gurobi and SCIP disagree, by cause.}
\label{apptab:win_loss}
\end{table}

\paragraph{Wins.}
Sixteen of the 22 wins come from a stronger bound.  On the posynomial programs, the substitution
$x=e^{y}$ makes the problem convex and \GEAR{}'s certificate is tight at the root, while the baselines
relax each product of fractional powers separately and branch on the original variables.  On the
polygon instances, relaxing every term $r_i r_j\sin(\theta_i-\theta_j)$ separately leaves Gurobi with the
trivial bound; \GEAR{} recognizes the structure and verifies a bound with a dedicated certificate.  On the symmetric log-barrier objectives, ordering the
symmetric variables and tightening the bounds of the lifted relaxation against the incumbent collapses
the box around the optimum.  On the remaining six, a baseline's returned point violates a constraint by more than
$10^{-6}$ when re-evaluated in AMPL: by at most $1.7\cdot10^{-5}$ on five, and by 0.16 on \texttt{elena\_hs54} (Gurobi).
\GEAR{} solves its 22 wins in a median of 9.8\,s (at most 114\,s).

\paragraph{Losses.}
Sixteen of the 18 losses are polynomial QCQPs with many bilinear constraints: water and wastewater
networks, pooling and circle packing.  On each of them \GEAR{}'s lifted McCormick LP and Gurobi's root
relaxation have the same or nearly the same value, so neither starts with a better relaxation; the
difference is the tree search.  Gurobi re-solves a warm-started LP with propagation and cuts at every
node, at $2.9\cdot10^{3}$ to $4.0\cdot10^{4}$ nodes per second.  \GEAR{}'s CPU lifted-LP tree processes
1.1-207 nodes per second, and after its open leaves are handed to the GPU, the GPU phase closes a median
of 0.74\% of the gap that remains (9 of 14 instances close less than 1\%): the batched GPU bounds are
weaker than the parent's exact LP.  On \texttt{eq6\_1} (square roots) no lifted LP applies, and the GPU
bounds alone stay far from the incumbent.  The baselines solve their 18 wins in a median of 1.2\,s.  This is the
regime discussed in \cref{subsec:exp_gpu_benefit}, where per-node LP throughput, not batched bounding,
decides the outcome.

\subsection{Limitations}
\label{subsec:limits}

\paragraph{Scope.}
\GEAR{} targets bounded continuous computation-graph NLPs with supported factorable operators and domain guards for singular primitives.
It does not cover unbounded variables, unsupported operators, integer variables, or exact equality handling beyond tolerance-paired inequalities.

\paragraph{Large bilinear models.}
\GEAR{} is weakest on large models with many bilinear constraints (\cref{appsec:win_loss}).
On these models its lifted LP and Gurobi's root relaxation have the same or nearly the same value, its CPU tree processes far fewer nodes per second than Gurobi's, and the batched GPU bounds are weaker than the parent's exact LP.

\paragraph{Conditional convergence.}
The convergence result in \cref{app:termination} is dimension-exponential and conditional.
\Cref{thm:eps_global} relies on Assumption~\ref{ass:oracle_and_primal}, including a $\delta$-complete primal-recovery hypothesis on ALM-PGD and a small-scale infeasibility-pruning hypothesis.
Neither hypothesis is needed for lower-bound or pruning soundness; both are needed for finite termination.

\paragraph{Operator domain restrictions.}
Fused division requires a denominator interval bounded away from zero.
Non-integer powers require guarded nonnegative bases.
Trigonometric relaxations require period localization before convex or concave envelopes are applied.
Outside these domains, the implementation must split, rewrite, or reject the operation.

\section{More Related Works}
\label{appsec:more_related}

\paragraph{Linear bound propagation and neural network verification.}
Linear bound propagation was introduced for ReLU networks by CROWN~\citep{zhang2018crown} and
DeepPoly~\citep{singh2019abstract} and generalized to arbitrary computation graphs by
auto\_LiRPA~\citep{xu2020automatic}.  $\alpha$-CROWN~\citep{xu2021fast} optimizes the slopes of the
operator relaxations by gradient ascent, and $\beta$-CROWN~\citep{wang2021betaCROWN} adds Lagrangian
multipliers for split constraints inside bound propagation.  \GEAR{} differs in what it optimizes: it
keeps the operator relaxations fixed (apart from the closed-form McCormick plane selection of
\cref{sec:operators}) and optimizes only nonnegative multipliers $\lambda$ on the constraints' linear
bounds after propagation, so an ascent step costs $O(mn)$ per box instead of a new propagation.
GCP-CROWN~\citep{zhang2022general} and BICCOS~\citep{zhou2024biccos} add cutting planes to the bound
computation, and Clip-and-Verify~\citep{zhou2025clip} uses linear constraints to clip input and
intermediate domains; \GEAR{} adopts relaxed and complete clipping from the latter and replaces its
coordinate descent by gradient steps when there are many constraints (\cref{subsec:dual}).  Other
verifiers use MILP and SMT formulations~\citep{tjeng2017evaluating,katz2017reluplex,anderson2020strong},
semidefinite relaxations~\citep{raghunathan2018semidefinite,dathathri2020enabling,lan2022tight,chiu2025sdpcrown},
and branch-and-bound over activations or inputs~\citep{bunel2018unified,bunel2020branch,de2021improved,ferrari2022complete,shi2025genbab}.
These methods certify a fixed property; none of them maintains a primal incumbent that moves the
pruning threshold.

\paragraph{Relaxations of factorable programs.}
Factorable programming builds a relaxation of a nonlinear program from relaxations of its elementary
operators~\citep{tawarmalani2002convexification,liberti2006reformulation,khajavirad2018hybrid}.  McCormick
envelopes~\citep{McCormick1976} give the convex and concave envelopes of a bilinear term on a box, and
RLT~\citep{sherali2013reformulation} multiplies constraints to tighten them further.  An LP relaxation
keeps all McCormick planes at once and couples them through auxiliary variables, while backward
linear bound propagation uses one plane per term, chosen by the sign of its coefficient.  \GEAR{}'s
plane selection by concretized value and its direct quotient relaxation (\cref{sec:operators}) narrow
this gap while keeping one linear bound per output.

\paragraph{Bound tightening and Lagrangian bounds.}
Optimization-based bound tightening solves one LP per variable bound; Gleixner
et al.~\citep{gleixner2017obbt} derive cheaper Lagrangian variable bounds from LP duals and reuse them
during the search.  The Lagrangian bound of \cref{eq:dual_lb} is the same weak-duality construction
applied to linear bounds from propagation instead of an LP relaxation, clipping tightens variable and
intermediate bounds without LP solves, and the infeasibility test of \cref{eq:infeas_cert} is a
Farkas-type certificate over the same linear bounds.

\paragraph{Branch-and-bound search.}
Spatial branch-and-bound descends from the branch-and-bound method of Land and
Doig~\citep{doig1960automaticmethodofsolvingdiscrete}; Morrison et al.~\citep{morrison2016babsurvery} survey
strategies for searching, branching and pruning.  Global NLP solvers usually branch on variables of
nonconvex terms that are violated at the LP solution, while verification branch-and-bound scores
splits by their estimated effect on the bound~\citep{bunel2020branch,lu2019neural}.  \SBReweight{}
(\cref{subsec:branch}) is closer to the latter: it weights the coefficient magnitudes of the linear
bounds by box width and by the Lagrangian multipliers, without solving an LP.

\paragraph{Primal recovery and local NLP solvers.}
Augmented Lagrangian methods go back to Hestenes and Powell; we use the inequality form analyzed by
Rockafellar~\citep{rockafellar1973augmented,rockafellar1973dual}.  Local solvers such as
IPOPT~\citep{wachter2006implementation}, SNOPT~\citep{gill2005snopt}, MINOS~\citep{murtagh2003minos},
CONOPT~\citep{drud1994conopt} and Knitro~\citep{byrd2006knitro} solve KKT systems with factorizations.
\GEAR{}'s projected ALM uses only graph evaluations, automatic differentiation and Adam
steps~\citep{kingma2014adam}, so thousands of searches run in one batch, and the Lagrangian
multipliers $\lambda$ of the bound only initialize its own multipliers $\mu$ (\cref{subsec:primal}).

\paragraph{Optimization over trained networks and GPU solvers.}
A related line of work optimizes over trained neural networks with mixed-integer
formulations~\citep{fischetti2017deepneuralnetworksasmilp,anderson2020strong,huchette2023deep},
bringing optimization techniques to verification; \GEAR{} goes in the opposite direction and brings
verification bounds to optimization.  GPU LP solvers such as PDLP~\citep{applegate2021practical,applegate2025pdlp}
accelerate one large LP with first-order iterations, whereas \GEAR{} parallelizes across many small
bounding problems in the branch-and-bound tree.


\newpage
\section*{NeurIPS Paper Checklist}

\begin{enumerate}

\item {\bf Claims}
    \item[] Question: Do the main claims made in the abstract and introduction accurately reflect the paper's contributions and scope?
    \item[] Answer: \answerYes{}.
    \item[] Justification: The abstract and \cref{sec:intro} state the bounded computation-graph NLP scope, the main contributions, and the comparison with six global solvers; \cref{sec:limits_related} states the main limitations.

\item {\bf Limitations}
    \item[] Question: Does the paper discuss the limitations of the work performed by the authors?
    \item[] Answer: \answerYes{}.
    \item[] Justification: \Cref{subsec:limits}  discusses unsupported operators, bounded-domain scope, conditional convergence, and large bilinear models.

\item {\bf Theory assumptions and proofs}
    \item[] Question: For each theoretical result, does the paper provide the full set of assumptions and a complete (and correct) proof?
    \item[] Answer: \answerYes{}.
    \item[] Justification: \cref{appsec:proofs} gives full proofs for the main theoretical claims.

    \item {\bf Experimental result reproducibility}
    \item[] Question: Does the paper fully disclose all the information needed to reproduce the main experimental results of the paper to the extent that it affects the main claims and/or conclusions of the paper (regardless of whether the code and data are provided or not)?
    \item[] Answer: \answerYes{}.
    \item[] Justification: \Cref{sec:exp} specifies the 505 instances, solver criteria, hardware, time limit, and parsing policy; \cref{appsec:dataset} gives canonicalization and filtering details.

\item {\bf Open access to data and code}
    \item[] Question: Does the paper provide open access to the data and code, with sufficient instructions to faithfully reproduce the main experimental results, as described in supplemental material?
    \item[] Answer: \answerYes{}.
    \item[] Justification: The current draft describes the benchmark construction and audit fields and includes them in the supplementary files.

\item {\bf Experimental setting/details}
    \item[] Question: Does the paper specify all the training and test details (e.g., data splits, hyperparameters, how they were chosen, type of optimizer) necessary to understand the results?
    \item[] Answer: \answerYes{}.
    \item[] Justification: The experimental protocol, the \GEAR{} and baseline configurations, time limits, and solved criteria are specified in \cref{sec:exp} and \cref{app:baseline}; dataset filtering is specified in \cref{appsec:dataset}.

\item {\bf Experiment statistical significance}
    \item[] Question: Does the paper report error bars suitably and correctly defined or other appropriate information about the statistical significance of the experiments?
    \item[] Answer: \answerYes{}.
    \item[] Justification: The main experiment is a deterministic solver benchmark over the full retained instance set, so \cref{fig:exp_cactus_morandi} reports exact coverage counts rather than sampling-based error bars.

\item {\bf Experiments compute resources}
    \item[] Question: For each experiment, does the paper provide sufficient information on the computer resources (type of compute workers, memory, time of execution) needed to reproduce the experiments?
    \item[] Answer: \answerYes{}.
    \item[] Justification: \Cref{sec:exp} reports the 180 s wall-clock limit, solver pairing, and the NVIDIA RTX 5090 GPU used for \GEAR{} runs.
    
\item {\bf Code of ethics}
    \item[] Question: Does the research conducted in the paper conform, in every respect, with the NeurIPS Code of Ethics \url{https://neurips.cc/public/EthicsGuidelines}?
    \item[] Answer: \answerYes{}.
    \item[] Justification: The work is an algorithmic solver study using public benchmark instances and does not involve human subjects, private data, or deployed decision systems.

\item {\bf Broader impacts}
    \item[] Question: Does the paper discuss both potential positive societal impacts and negative societal impacts of the work performed?
    \item[] Answer: \answerNA{}.
    \item[] Justification: The contribution is a global optimization algorithm for computation-graph NLPs; the paper does not introduce a model or dataset with direct societal deployment claims.
    
\item {\bf Safeguards}
    \item[] Question: Does the paper describe safeguards that have been put in place for responsible release of data or models that have a high risk for misuse (e.g., pre-trained language models, image generators, or scraped datasets)?
    \item[] Answer: \answerNA{}.
    \item[] Justification: The paper does not release high-risk generative models, scraped datasets, or other assets requiring misuse safeguards.

\item {\bf Licenses for existing assets}
    \item[] Question: Are the creators or original owners of assets (e.g., code, data, models), used in the paper, properly credited and are the license and terms of use explicitly mentioned and properly respected?
    \item[] Answer: \answerNo{}.
    \item[] Justification: The benchmark sources and solver families are cited, but the draft does not yet enumerate all source licenses and terms of use.

\item {\bf New assets}
    \item[] Question: Are new assets introduced in the paper well documented and is the documentation provided alongside the assets?
    \item[] Answer: \answerYes{}.
    \item[] Justification: The paper introduces a 505-instance benchmark; \cref{appsec:dataset} documents its construction, filters, and released tables.

\item {\bf Crowdsourcing and research with human subjects}
    \item[] Question: For crowdsourcing experiments and research with human subjects, does the paper include the full text of instructions given to participants and screenshots, if applicable, as well as details about compensation (if any)? 
    \item[] Answer: \answerNA{}.
    \item[] Justification: The paper does not involve crowdsourcing or research with human subjects.

\item {\bf Institutional review board (IRB) approvals or equivalent for research with human subjects}
    \item[] Question: Does the paper describe potential risks incurred by study participants, whether such risks were disclosed to the subjects, and whether Institutional Review Board (IRB) approvals (or an equivalent approval/review based on the requirements of your country or institution) were obtained?
    \item[] Answer: \answerNA{}.
    \item[] Justification: The paper does not involve human subjects, so IRB or equivalent review is not applicable.

\item {\bf Declaration of LLM usage}
    \item[] Question: Does the paper describe the usage of LLMs if it is an important, original, or non-standard component of the core methods in this research? Note that if the LLM is used only for writing, editing, or formatting purposes and does \emph{not} impact the core methodology, scientific rigor, or originality of the research, declaration is not required.
    \item[] Answer: \answerNA{}.
    \item[] Justification: LLMs are not an important, original, or non-standard component of the proposed solver methodology.

\end{enumerate}

\end{document}